\documentclass[11pt]{article}

\usepackage{tikz} 
\usepackage{amsmath}
\usepackage{mathrsfs} 
\usepackage{amsfonts} 
\usepackage{pgfplots,diagbox,multirow}
\usepackage{amssymb}

\usepackage[colorlinks, linkcolor=blue, anchorcolor=blue,urlcolor=blue, citecolor=blue]{hyperref}
\usepackage{subfigure}

\usepackage{appendix}

\def\ar{\!\!\!&} 

\def\proof{\noindent{\it Proof.~~}} 
\def\qed{\hfill$\Box$\medskip} 
\usepackage[top=2.5cm,bottom=2.5cm, outer=2.5cm, inner=2.5cm]{geometry}

 \usepackage{caption}

\newtheorem{theorem}{Theorem}[section]

\newtheorem{lemma}[theorem]{Lemma}
\newtheorem{corollary}[theorem]{Corollary}
\newtheorem{proposition}[theorem]{Proposition}
\newtheorem{remark}[theorem]{Remark}
\newtheorem{condition}[theorem]{Condition}

\newtheorem{assumption}[theorem]{Assumption}

\def\red{\color{red}}

\def\beqlb{\begin{eqnarray}}\def\eeqlb{\end{eqnarray}} 
\def\beqnn{\begin{eqnarray*}}\def\eeqnn{\end{eqnarray*}} 

\DeclareMathOperator*{\arginf}{arg\,inf}

\def\ar{\!\!\!&}

\def\proof{\noindent{\it Proof.~~}}\def\qed{\hfill$\Box$\medskip}

\begin{document}

	\title{\Large\bf  ASYMPTOTICS FOR EXPONENTIAL FUNCTIONALS OF RANDOM WALKS}
	
	\author{\sc\normalsize By Wei Xu\footnote{Department of Mathematics, Humboldt-Universit\"at zu Berlin, Germany. E-mail: \href{mailto:xuwei@math.hu-berlin.de}{xuwei.math@gmail.com}}}
	\date{}
	\maketitle

	\begin{abstract}
		This paper provides a detailed description for the asymptotics of exponential functionals of random walks with light/heavy tails.
	    We give the convergence rate based on the key observation that the asymptotics depends on the sample paths with either slowly decreasing local minimum or final value below a low level.
	    Also, our thoughtful analysis of the interrelationship between the local minimum and the final value provides the exact expression for the limiting coefficients in terms of some transformations of the random walk.
	\end{abstract}

    {\small	\textit{MSC2020 subject classifications:} Primary 60G50, 60J55; secondary 60F10,  60E05.
	
     }
	
	\renewcommand{\baselinestretch}{1.15}

  \section{Introduction}
 \setcounter{equation}{0}
 
 Exponential functionals of random walks have been studied deeply in the past decades because of their wide and important applications in various fields such as mathematical finance, physics and population evolution.
 In the general setting, for a one-dimension random walk $S$ generated by a sequence of independent copies $X_1,X_2,\cdots$ of a random variable $X$, its  exponential functional is usually defined by
 \beqnn
 I_n:= \sum_{k=1}^n e^{-S_k},\quad  n=1,2,\cdots
 \eeqnn
 with the convention $I_0=0$. It increases to the limit $I_\infty \in(0,\infty]$ almost surely as $n\to\infty$.
 According to the Blumenthal zero-one law, the probability that $I_\infty$ is finite, $\mathbf{P}(I_\infty<\infty)$, necessarily equals to $0$ or $1$.
 Furthermore, the long-term behavior of $S$ shows that this probability equals to $0$ if and only if $S$ does not drift to infinity.
 In this case, much attention has been drawn to the speed at which $I_n$ increases to infinity, especially the decay rate of the following expectation
 \beqlb\label{ExpRW}
 \mathbf{E}\big[F(I_n)\big]= \mathbf{E}\Big[F\Big(\sum_{k=1}^n e^{-S_k}\Big)\Big]
 \eeqlb
 with $F$ be a positive function on $(0,\infty)$ that vanishes at infinity.
 The main aim of this work is to provide a detailed description for the decay rate of this expectation as $n\to\infty$.
 
 The study of the expectation $\mathbf{E}[F(I_n)]$ has been initiated to explore the asymptotic properties of stochastic systems in random environment.
 For instance, let $\tau$ be the extinction time of a linear fractional Galton-Watson process in an i.i.d. random environment  with geometric offspring distribution\footnote{The asymptotic results in the following references are established for Galton-Watson processes in i.i.d. random environment with  general branching mechanism and the linear fractional case is considered as a typical example.}. The survival probability of the population at time $n$ can be written as
 \beqnn
 \mathbf{P}(\tau>n)=  \mathbf{E}\big[(1+I_n)^{-1}\big]
 \eeqnn
 in which the random walk $S$ is generated by the  offspring distribution that varies randomly as time goes; see Application~1 in  \cite{Hirano1998}.
 Under the assumption that $\mathbf{E}[X]=0$ and $\mathbf{E}[|X|^2]<\infty$, Kozlov \cite{Kozlov1976} first showed that $ \mathbf{P}(\tau>n)\sim C\cdot n^{-1/2}$ as $n\to\infty$ for some $C>0$. 
 When $\mathbf{E}[X]<0$, Liu \cite{Liu1993} proved a rough asymptotic result
 \beqnn
 \big(\mathbf{P}(\tau>n)\big)^{1/n}\sim C\cdot \zeta
 \quad \mbox{with}\quad \zeta :=\inf_{\lambda \in[0,1]}\mathbf{E}[e^{\lambda X}].
 \eeqnn
 More precise results were provided later in \cite{AfanasyevBoinghoffKerstingVatutin2014,AfanasyevGeigerKerstingVatutin2005b,GuivarchLiu2001}, i.e,  $ \mathbf{P}(\tau>n)\sim C\cdot \zeta^n $, $C\cdot n^{-1/2}\zeta^n$ or $C\cdot n^{-3/2}\zeta^n$ if $\mathbf{E}[|X|^2]<\infty$ and $\mathbf{E}[Xe^X]<0$, $=0$ or $>0$.
 As the second example, we write $\mathcal{M}$ for the maximum of a sample transient random walk moving in an i.i.d. random environment. The tail probability of $\mathcal{M}$ can be given by
 \beqnn
 \mathbf{P}(\mathcal{M}\geq n)= \mathbf{E}\big[ \hat{I}'_\infty( \hat{I}'_\infty+ I_n)^{-1}\big], \quad n=1,2,\cdots,
 \eeqnn
 where $ \hat{I}'_\infty:= 1+ \sum_{i=-1}^{-\infty}e^{S'_i}$ and the two random walks $\{S'_k:k=-1,-2,\cdots\}$ and $\{S_k:k=1,2,\cdots\}$ are generated by the random environment at negative and positive integer points respectively; see (8) in \cite{Afanasev1990}.
 The asymptotics of $\mathbf{P}(\mathcal{M}\geq n)$ was considered in \cite{Afanasev1990} with  $\mathbf{E}[X]<0$, i.e., it is asymptotically equivalent to  $C\cdot\zeta^n $, $C\cdot n^{-1/2}\zeta^n$ or $C\cdot n^{-3/2}\zeta^n$ when $\mathbf{E}[Xe^X]<0$, $=0$ or $>0$ respectively.
 
 For general function $F$, to our best knowledge the asymptotics of the expectation (\ref{ExpRW}) has only been studied in Hirano \cite{Hirano1998} under the assumption that $\mathbf{E}[Xe^{\alpha X}]=0$ and $F$ is completely monotone satisfying $ F(x)\leq C x^{-\beta} $ for  two constants $\beta>\alpha >0$.
 His results show that $ \mathbf{E}[F(I_n)] \sim C\cdot n^{-3/2} (\mathbf{E}[e^{\alpha X}])^n$.
 In this work, we provide under natural assumptions, a more complete and accurate description for the asymptotic behavior of the expectation (\ref{ExpRW}) with general function $F$ and random walk $S$ whose generic step $X$ owns light- or heavy-tailed distribution.
 Our results state that beside of the optimal polynomial decay rate of $F$, i.e.,  $\theta_F:=\sup\{ \theta>0: \sup_{x>0}x^\theta F(x)<\infty\}$, the long-term properties of the expectation (\ref{ExpRW}) also heavily rely on the speed of the random walk $S$ decreasing to $-\infty$, which can be quantified by 
 \beqlb\label{Lambda}
 \varrho :=\inf_{\lambda \in[0,\theta_F]}\mathbf{E}[e^{\lambda X}] , \quad \Lambda:= \arginf_{\lambda \in[0,\theta_F]}\mathbf{E}[e^{\lambda X}]
 \eeqlb
 and the behavior of $S_n$ under $\mathbf{P}^{(\Lambda)}$ as $n\to\infty$. Here the Laplace transform is allowed to be infinite for some $\lambda \in [0,\theta_F]$ and $\arginf$ stands for argument of the infimum and the probability law $\mathbf{P}^{(\Lambda)}$ is the Esscher transform of $\mathbf{P}$ associated to the martingale $\{\varrho^{-n}\cdot e^{\Lambda S_n}:n=0,2,\cdots  \}$. 
 Consequently, we see that the convergence rate of the expectation (\ref{ExpRW}) changes dramatically in the following six disjoint cases:
 \begin{center}
 	\begin{tabular}{|c|c|c|c|c|c|}
 		\hline
 		\multicolumn{3}{|c|}{$S$ oscillates under $\mathbf{P}^{(\Lambda)}$} & \multicolumn{3}{|c|}{$S$ drifts to $-\infty$ under $\mathbf{P}^{(\Lambda)}$} \\
 		\hline
 		$0=\Lambda<\theta_F$ & $0<\Lambda<\theta_F$  & $0<\Lambda=\theta_F$ & $0=\Lambda<\theta_F$ & $0<\Lambda<\theta_F$&  $0<\Lambda=\theta_F$\\
 		\hline
 		Figure~\ref{Figure1}(a) & Figure~\ref{Figure1}(b) & Figure~\ref{Figure1}(c) & Figure~\ref{Figure1}(d) & Figure~\ref{Figure1}(e) & Figure~\ref{Figure1}(f)\\
 		\hline
 		Case~I &Case~II&Case~III&Case~IV& Case~V&Case~VI\\
 		\hline
 	\end{tabular}
 	
 \end{center}
 The feature of this work is that both the convergence rate and the   expression of the limiting coefficient in each regime are given.
 Roughly speaking, the expectation (\ref{ExpRW}) decays at some exponential rate with different regularly varying modifying factor, i.e.,
 \beqnn
 \mathbf{E}[F(I_n)] \sim \psi_n \cdot \varrho^n , 
 \eeqnn
 where  $\psi_n$ is a regularly varying sequence vanishing at infinity and determined by the fluctuation of the random walk $S$ under the probability law $\mathbf{P}^{(\Lambda)}$. 
 Particularly, when $0\leq \Lambda<\theta_F$ the comparison with the asymptotics of the first entrance time $\tau_0^-$ of the random walk $S$ in $(-\infty,0)$ shows that for some constant $C>0$ depending on $F$,
 \beqnn
 \mathbf{E}[F(I_n)]  \sim C \cdot \mathbf{P}(\tau_0^->n).
 \eeqnn
 
 Our methodology is based on the fluctuation theory for random walks.
 There are two major challenges, which are also two key steps, in the asymptotic analysis of the expectation (\ref{ExpRW}).
 The first one is to find out the sample paths of random walk that make the main contribution to the expectation (\ref{ExpRW}).
 Enlightened by the asymptotic analysis of survival probabilities of Galton-Watson processes in an i.i.d. random environment; see \cite{AfanasyevGeigerKerstingVatutin2005,AfanasyevBoinghoffKerstingVatutin2012,AfanasyevBoinghoffKerstingVatutin2014,BansayeVatutin2017,VatutinZheng2012}, we  observe that the characteristics of these sample paths vary dramatically in different cases.
 More precisely, sample paths with slowly decreasing local minimum make the main contribution in Case~I and III; the expectation (\ref{ExpRW}) in Case~II  mainly replies on sample paths with either the local minimum attained at the beginning of the time interval or the final value below a low level. In Case~IV, the main contribution is made by sample paths with early large step. 
 The key sample paths in Case~V not only have early large step and low final value but also  decay slowly before the large step.
 Different to other cases, the contribution of any sample path can not be asymptotically ignored in Case~VI.
 The second challenge is to provide an exact expression of the limiting coefficients.
 In order to achieve it, we decompose each key sample path into several parts at its local minimum and then seek out the subdivisions that make the main contribution to the expectation (\ref{ExpRW}). 
 For instance, the contribution of the key sample paths after large time $K$ is negligible in Case I, III and IV, i.e., the impact of $I_n$ can be well approximated by that of $I_K$. 
 However, in Case II and V the key sample paths make the main contribution to the expectation (\ref{ExpRW}) at both the beginning and the end of the time interval, i.e., the impact of $I_n$ can be well approximated by that of $I_K+ (I_n-I_{n-K})$ for large $K$. 
 Finally, we give the representations for the limiting coefficients in terms of some transformations of these subdivisions.

 The second purpose of this work is to offer assistance in our future study of the asymptotic behavior of exponential functional of a L\'evy process $\{\xi_t:t\geq 0\}$ defined by
 \beqnn
 I_t(\xi):= \int_0^t e^{-\xi_s}ds, \quad t\geq 0.
 \eeqnn
 Similarly, we have $ I_t(\xi) \to I_\infty(
 \xi) \in(0,\infty]$ a.s. as $t\to\infty$  and $I_\infty(\xi)<\infty$ a.s. if and only if $\xi$ drifts to infinity. 
 Readers may refer to \cite{BarkerSavov2021,CarmonaPetitYor1997, PardoPatieSavov2012,PatieSavov2018,Vechambre2019} and references therein for many interesting results related to $I_t(\xi)$ and $I_\infty(\xi)$.
 In the case $I_\infty(\xi)=\infty$ a.s., we are usually interested in the decay rate of the expectation $\mathbf{E}[F(I_t(\xi))]$ defined as in (\ref{ExpRW}), because of its close connection to the long-term properties of random processes in random environment, e.g. continuous-state branching processes in L\'evy random environment and diffusion in L\'evy random environment; see \cite{BansayeMillanSmadi2013,  HeLiXu2018, LiXu2018, PalauPardoSmadi2016}.
 Under the assumption that $F(x)\leq C x^{\theta_F}$ and $\theta_F<\lambda_*:=\sup\{\lambda\geq 0: \mathbf{E}[e^{\lambda\xi_1}]<\infty\}$,   four different
 regimes for the convergence rate of the expectation (1.2) were provided in \cite{BansayeMillanSmadi2013, PalauPardoSmadi2016, LiXu2018}.
 When $\lambda_*=0$, Patie and Savov \cite{PatieSavov2018} and Xu \cite{Xu2021a} proved the polynomial decay rate for $\mathbf{E}[F(I_t(\xi))]$ with  $\xi$  oscillating and satisfying the Spitzer's condition or $\mathbf{E}[\xi_1]<0$ and $\mathbf{P}(\xi_1>x)$ regularly varying at infinity.
 For the case $0<\lambda_*\leq\theta_F$,  to the best of our knowledge, the asymptotics of $\mathbf{E}[F(I_t(\xi))]$ is still an open problem.
 Notice that $\xi$ can be well approximated by the random walk $\xi^\delta:=\{ \xi_n^\delta:n=0,1,\cdots \}$ with $\xi_n^\delta:=\xi_{n\delta}$ for small $\delta>0$. Let $I^{\delta}$ be the exponential functional of $\xi^\delta$, we may conjecture that there exists a positive, regularly varying function $\Psi$ on $\mathbb{R}_+$ satisfying $\Psi(n)\sim \psi_n$ as $n\to\infty$ and
 \beqnn
 \mathbf{E}[F(I_t(\xi))] \sim \mathbf{E}[F(\delta\cdot I^{\delta}_{[t/\delta]})] \sim \Psi(t)\cdot |\varrho_\xi|^t
 \quad \mbox{and}\quad 
 \varrho_\xi := \inf_{\lambda\in[0,\theta_F]} \mathbf{E}\big[e^{\lambda \xi_1}\big].
 \eeqnn
 

 The remainder of this paper is organized as follows. In Section~\ref{Sec.MainR}, we recall some basic
 elements of fluctuation theory for random walks and then provide the accurate asymptotic results for the expectation (\ref{ExpRW}). We give in Section~\ref{Sec.AuxiliaryR} some auxiliary asymptotic results for random walks, which will be used in the proofs for our main results. The asymptotic results for the expectation (\ref{ExpRW}) with $S$  oscillating or drifting to $-\infty$ under $\mathbf{P}^{(\Lambda)}$ are proved separately in Section~\ref{Sec.Oscilliate} and \ref{Sec.NegD}.
 
  \section{Preliminaries and main results} \label{Sec.MainR}
 \setcounter{equation}{0}
 
 In this section we first introduce some basic notation and elements of fluctuation theory for random walks.
 We then provide the main results in this paper about the asymptotic behavior of exponential functionals of random walks.
 
 \subsection{Random walks}
 
 Suppose that the random walk $S$ is defined on a complete probability space $(\Omega, \mathscr{F},\mathbf{P})$.
 Let $\mathscr{F}^S$ and $(\mathscr{F}^S_n)_{n\geq 0}$ denote the $\sigma$-algebra and the filtration generated by $S$, i.e., $\mathscr{F}^S:= \sigma(S_k:k=0,1,\cdots)$ and  $\mathscr{F}^S_n:=\sigma(S_k:k=0,1,\cdots,n)$.
 For any probability measure $\mu$ on $\mathbb{R}$, we denote by $\mathbf{P}_\mu$ and $\mathbf{E}_\mu$ the law and expectation of the random walk $S$ with initial state $S_0$ distributed as $\mu$.
 When $\mu=\delta_x$ is a Dirac measure at point $x\in\mathbb{R}$, we write $\mathbf{P}_x$ for $\mathbf{P}_{\delta_x}$ and $\mathbf{E}_x$ for $\mathbf{E}_{\delta_x}$.
 For simplicity, we also write $\mathbf{P}$ for  $\mathbf{P}_0$ and $\mathbf{E}$ for $\mathbf{E}_0$. 
 
 According to its behavior as $n\to\infty$, the random walk $S$ can be exactly classified into three types:
 (i) drifts to $\infty$: $S_n\to \infty$ a.s.; (ii) drifts to $-\infty$: $S_n\to -\infty$ a.s.; (iii) oscillation: $\limsup_{n\to\infty} S_n=-\liminf_{n\to\infty} S_n=\infty$ a.s.
 We write $M:=\{M_n:n=1,2,\cdots  \}$ and $L:=\{L_n:n=1,2,\cdots  \}$ for the {\it running maximum} and {\it minimum processes} respectively,
 \beqnn
 M_n:= \max_{1\leq i\leq n} S_i
 \quad \mbox{and}\quad
 L_n:=  \min_{1\leq i\leq n} S_i .
 \eeqnn
 In addition to $M_n$ and $L_n$, we will also use the following two random times:
 \beqnn
 \sigma^+_n:= \inf \big\{0\leq i\leq n: S_i=S_0\vee M_n \big\}
 \quad \mbox{and}\quad
 \sigma^-_n:= \inf\big\{ 0\leq i\leq n: S_i= S_0\wedge  L_n \big\}
 \eeqnn
 the {\it first times} up to time $n$ that the maximum and minimum are attained.
 For $x\in\mathbb{R}$, denote by $\tau_x^+$ and $\tau_x^-$ the {\it first entrance times} of $S$ in $(x,\infty)$ and $(-\infty,x)$ respectively, i.e.,
 \beqnn
 \tau_x^+ :=\inf\{n>0: S_n>x \}
 \quad \mbox{and}\quad
 \tau_x^- :=\inf\{n>0: S_n<x \}.
 \eeqnn
 In the sequel, we always denote by $\tilde{S}$ an independent copy of $S$.  For the quantities introduced to $S$, the corresponding ones for  $\tilde{S}$ are denoted by tildes, for instance,  $\tilde{I}$ and $\tilde{\tau}_x^-$.

 The {\it renewal function} $V$  associated with the  strict descending ladder height process is defined by $V(x)=0$ if $x<0$ and
 \beqnn
 V(x):= 1+ \sum_{n=1}^\infty \mathbf{P}(S_{\gamma_n^-} \geq -x )  ,\quad x\geq 0,
 \eeqnn
 where $\{ \gamma_i^-: i=0,1,\cdots \}$ are {\it strict descending ladder epochs} of $S$ with $\gamma_0^-=0$ and
 \beqnn
 \gamma_i^- := \inf \big\{ n>\gamma_{i-1}^- : S_n<S_{\gamma_{i-1}^-} \big\},
 \quad i\geq 1.
 \eeqnn
 It is obvious that $V$ is a non-decreasing right-continuous function with $V(0)=1$.
 Using  the duality lemma, we also have
 \beqnn
 V(x) = 1+ \sum_{n=1}^\infty \mathbf{P}( S_n\geq -x , \sigma_n^-=n )
 = 1+ \sum_{n=1}^\infty \mathbf{P}(  S_n\geq -x  , M_n< 0 ).
 \eeqnn
 When the random walk $S$ drifts to $-\infty$, we have  $ \mathbf{E}_x[\tau_0^-]<\infty$ and $V(x)=\mathbf{E}_x[\tau_0^-]/\mathbf{E}[\tau_0^-]$ for any $x\geq 0$.
 
 The {\it dual process} of $S$ is denoted by $\hat{S}$, that is $\hat{S}=-S$ when the starting point is $0$.
 Let $\hat{\mathbf{P}}_x$ be the law of $\hat{S}$ under $\mathbf{P}_{-x}$.
 For the quantities introduced to $S$, the corresponding ones for $\hat{S}$ are denoted by hats, for instance $\hat{I}$, $\hat{\tau}_x^-$ and so on. Specially, the renewal process $\hat{V}$ associated with the  strict descending ladder height process of $\hat{S}$ is equal to the renewal process associated with the  strict ascending ladder height process of $S$, i.e.
 \beqnn
 \hat{V}(x) = 1+ \sum_{n=1}^\infty \mathbf{P}( S_n\leq x , \sigma_n^+=n )
 = 1+ \sum_{n=1}^\infty \mathbf{P}(  S_n\leq x  , L_n> 0 ).
 \eeqnn
 When  $S$ drifts to  $\infty$,   we also have  $ \mathbf{E}_x[\tau_0^+]<\infty$ and $\hat{V}(x)=\mathbf{E}_x[\tau_0^+]/\mathbf{E}[\tau_0^+]$  for any $x\geq 0$.

 If the random walk $S$ does not drift to $-\infty$, the process $\{ V(S_n)\mathbf{1}_{\{ \tau_0^->n \}} :n=0,1,\cdots \}$ is a $\mathbf{P}_x$-martingale for any $x\geq 0$.
 In this case, we introduce a probability $\mathbf{P}^\uparrow_x$ on $(\Omega,\mathscr{F},\mathscr{F}^S)$ defined by
 \beqnn
 \mathbf{P}^\uparrow_x(A):= \int_A V(S_n)\mathbf{1}_{\{ \tau_0^->n \}}  \frac{d\mathbf{P}_x}{V(x)},\quad A\in \mathscr{F}^S_n,n\geq 0.
 \eeqnn
 It is usually well-known as {\it Doob's $h$-transform} of $\mathbf{P}_x$. Particularly, under $ \mathbf{P}^\uparrow$ the process $S$ turns to be a homogeneous Markov process on  $[0,\infty)$ with transition function
 \beqnn
 p^\uparrow(x,dy) := \frac{V(y)}{V(x)} \mathbf{P}(x+X\in dy),\quad x,y\geq 0.
 \eeqnn
 Similarly, in the case that $\hat{S}$ does not drift to $-\infty$, for any $x\geq 0$ we can also introduce Doob's $h$-transform  of $\hat{\mathbf{P}}_x$  on $(\Omega,\mathscr{F},\mathscr{F}^{\hat{S}}_n)$ defined by the renewal function $\hat{V}$, i.e.,
 \beqnn
 \mathbf{P}^\downarrow_x(A):= \int_A \hat{V}(\hat{S}_n)\mathbf{1}_{\{ \hat\tau_0^->n \}} \frac{d\red \mathbf{P}_{x}}{\hat{V}(x)},\quad A\in \mathscr{F}^{\hat{S}}_n,n\geq 0.
 \eeqnn
 Under $ \mathbf{P}^\downarrow$ the process $\hat{S}$ turns to be a homogeneous Markov process on  $[0,\infty)$ with transition function
 \beqnn
 p^\downarrow(x,dy) := \frac{\hat{V}(y)}{\hat{V}(x)} \mathbf{P}(x-X\in dy),\quad x,y\geq 0.
 \eeqnn
 %
 In order to make our following statements much easier to be understood, we write $S^{\uparrow}$ and $S^{\downarrow}$ for the two independent Markov processes with transition functions $p^\uparrow$ and $p^\downarrow$ respectively.
 Their exponential functionals are denoted as $I^{\uparrow}$ and $I^{\downarrow}$.
 Moreover, we also write $\hat{S}^{\uparrow}$ and $\hat{S}^{\downarrow}$ for the dual processes of $S^{\uparrow}$ and $S^{\downarrow}$ respectively and also write  $\hat{I}^{\uparrow}$ and $\hat{I}^{\downarrow}$ for their exponential functionals.
 Repeating the preceding argument, we see that $\hat{S}^{\uparrow}$ and $\hat{S}^{\downarrow}$ are two homogeneous Markov processes  with transition functions
 \beqnn \red 
 \hat{p}^\uparrow(x,dy) := \frac{V(-y)}{V(-x)} \mathbf{P}(x-X\in dy)
 \quad \mbox{and}\quad
 \hat{p}^\downarrow(x,dy) := \frac{\hat{V}(-y)}{\hat{V}(-x)} \mathbf{P}(x+X\in dy), \quad x,y\leq 0.
 \eeqnn
 
 As we have mentioned before, the asymptotics of the tail distribution of generic step $X$ plays a crucial role in the following classification and asymptotic analysis of exponential functionals of random walks.
 Thus we need the Laplace transform of $X$
 \beqnn
 \mathcal{L}_X(\lambda):=  \mathbf{E} \big[e^{\lambda X}\big],\quad \lambda \in\mathbb{R}.
 \eeqnn
 Certainly, it may happen that $  \mathcal{L}_X(\lambda) =\infty$ for some $\lambda\in\mathbb{R} $.
 Let $\mathcal{D}_{\mathcal{L}_X}:= \{ \lambda\in\mathbb{R}: \mathcal{L}_X(\lambda)<\infty \}$, $\mathcal{D}^+_{\mathcal{L}_X}:= \mathcal{D}_{\mathcal{L}_X}\cap [0,\infty)$ and $\mathcal{D}^-_{\mathcal{L}_X}:= \mathcal{D}_{\mathcal{L}_X}\cap (-\infty,0]$.
 For each $\lambda_0 \in\mathcal{D}_{\mathcal{L}_X}$,  the process $\{ \exp\{ \lambda_0 S_n\} \cdot | \mathcal{L}_X(\lambda_0)|^{-n}:n=0,1,\cdots  \}$ is a $\mathbf{P}$ martingale, which allows us to define a probability measure $\mathbf{P}^{(\lambda_0)}$ on $(\Omega,\mathscr{F},\mathscr{F}^S)$
 \beqlb\label{MeasureChan}
 \mathbf{P}^{(\lambda_0)}(A):= \int_A  \exp\{ \lambda_0 S_n\}\cdot | \mathcal{L}_X(\lambda_0)|^{-n} d\mathbf{P}, \quad A \in\mathscr{F}^S_n, n\geq 0.
 \eeqlb
 It is known that the process $S$ under $\mathbf{P}^{(\lambda_0)}$ is still a random walk and the generic step $X$ has Laplace transform
 \beqnn
 \mathbf{E}^{(\lambda_0)}[e^{\lambda X}]=
 \mathcal{L}_X(\lambda_0+\lambda)/\mathcal{L}_X(\lambda_0),\quad \lambda\in\mathbb{R}.
 \eeqnn
 Let $V^{(\lambda_0)}$ and $\hat{V}^{(\lambda_0)}$ be the renewal functions associated with the  strict descending ladder height processes of $S$ and $\hat{S}$ under $\mathbf{P}^{(\lambda_0)}$.
 

 \subsection{Main results}
 
 We now provide the asymptotic results for the expectation $ \mathbf{E}[F(I_n)]$ in which $F$ is a {\it positive, bounded} function on $(0,\infty)$ and always satisfies the following assumption:
 \begin{assumption}\label{Con.F.UpperBound}
 	The set $\mathcal{D}_F:=\{\theta> 0 : \sup_{x>0} x^\theta F(x)<\infty  \}$ is not null and the supremum is denoted as $\theta_F\in(0,\infty]$.
 \end{assumption}
 To simplify the representation of our main results, let us list the following conditions:
 \begin{condition}\label{Con.F.LipCon}
 	For each $\delta>0$, there exists a constant $C_\delta>0$ such that $|F(x)-F(y)|\leq C_\delta|x-y|$ for any $x,y\geq \delta$.
 \end{condition}
 \begin{condition}\label{Con.AsymPower}
 	There exist a constant $K_0>0$ such that $F(x)\sim K_0x^{-\theta_F}$ as $x\to\infty$.
 \end{condition}
 
 %
 
 It is known that when $S$ does not drift to infinity, with probability one it visits the negative half-line infinite times and hence $I_\infty=\infty$ a.s.
 For the converse,  when $S$ drifts to $\infty$,  Erickson's theorem in \cite{Erickson1973} shows that $\lim_{n\to\infty}S_n/n\in(0,\infty]$ a.s. and hence $I_\infty<\infty$ a.s.
 The first asymptotic result for exponential functionals of random walks is summarized as follows.
 
 \begin{lemma}\label{FirstLemma}
 	The following assertions are equivalent:\\
 	{\rm(i)}  $I_\infty<\infty$ a.s.;\quad {\rm(ii)}  $\mathbf{P}( I_\infty<\infty )>0$;\quad {\rm(iii)}  $S$ drifts to $\infty$, i.e., $S_n\to \infty$ a.s. as $n\to\infty$.
 \end{lemma}
 
\begin{figure}
 	\centering
 	
 	%
 	\subfigure{
 		\begin{minipage}[c]{0.31\linewidth}
 			\centering
 			\includegraphics[scale=0.53]{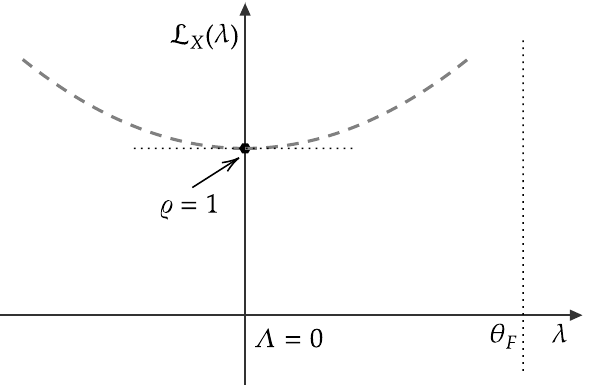}
 			{\scriptsize {\bf(a)} $0=\Lambda<\theta_F$,
 				$S$ is $\mathbf{P}$-oscillating}
 		\end{minipage}
 	}
 	\subfigure{
 		\begin{minipage}[c]{0.31\linewidth}
 			\centering
 			\includegraphics[scale=0.53]{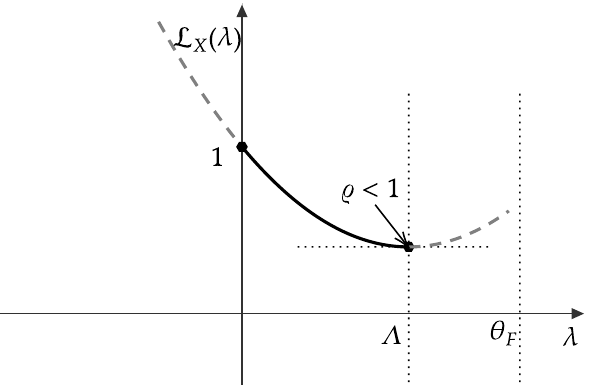}
 			{\scriptsize {\bf(b)} $0<\Lambda<\theta_F$,
 				$S$ is $\mathbf{P}^{(\Lambda)}$-oscillating}
 		\end{minipage}
 	}
 	\subfigure{
 		\begin{minipage}[c]{0.31\linewidth}
 			\centering
 			\includegraphics[scale=0.53]{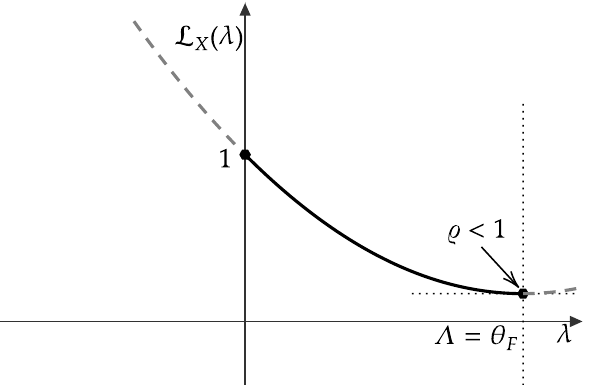}
 			{\scriptsize {\bf(c)} $0<\Lambda=\theta_F$,
 				$S$ is $\mathbf{P}^{(\Lambda)}$-oscillating}
 		\end{minipage}
 	}\\
 	
 	\centering
 	\subfigure{
 		\begin{minipage}[c]{0.31\linewidth}
 			\centering
 			\includegraphics[scale=0.53]{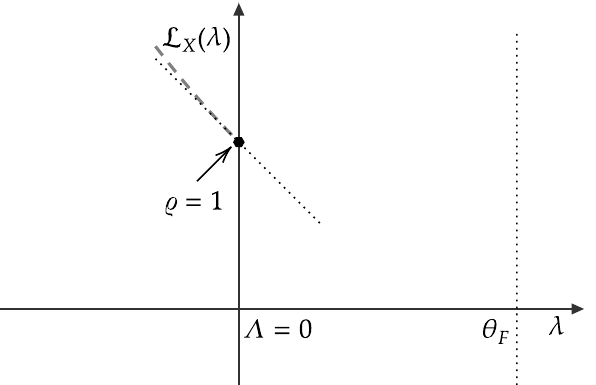}
 			{\scriptsize {\bf(d)} $0=\Lambda<\theta_F$, $\mathbf{E}[X]<0$}
 		\end{minipage}
 	}
 	\subfigure{
 		\begin{minipage}[c]{0.31\linewidth}
 			\centering
 			\includegraphics[scale=0.53]{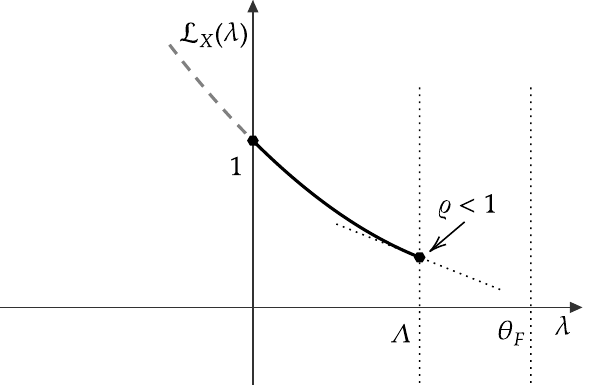}
 			{\scriptsize {\bf(e)} $0<\Lambda<\theta_F$, $\mathbf{E}^{(\Lambda)}[X]<0$}
 		\end{minipage}
 	}
 	\subfigure{
 		\begin{minipage}[c]{0.31\linewidth}
 			\centering
 			\includegraphics[scale=0.53]{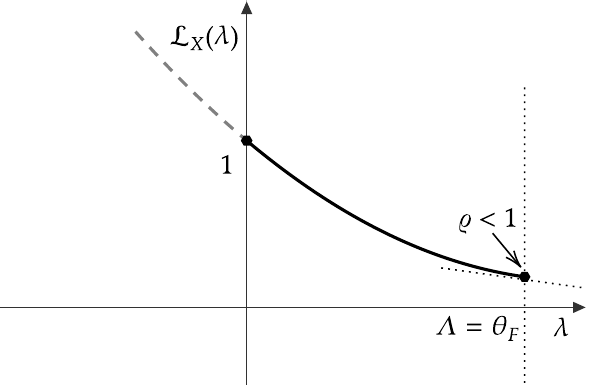}
 			{\scriptsize {\bf(f)} $0<\Lambda=\theta_F$, $\mathbf{E}^{(\Lambda)}[X]<0$}
 		\end{minipage}
 	}
 	\caption{\footnotesize 
 		Figures (a)-(c) draw the three possibilities of the Laplace transform of  random walk $S$ oscillating under $\mathbf{P}^{(\Lambda)}$, in which $\mathcal{D}^+_{\mathcal{L}_X} \supset [0,\Lambda]$. 
 		Figures (d)-(f) draw the three possibilities  of the Laplace transform of random walk $S$ with negative drift under $\mathbf{P}^{(\Lambda)}$,  in which $\mathcal{D}^+_{\mathcal{L}_X}=[0,\Lambda]$. 
 		The thick dotted lines in each figure represent the set besides $[0,\Lambda]$ on which the Laplace transform may be finite.}
 	\label{Figure1} 
 \end{figure}

 In this paper, we are mainly interested in the case $I_\infty=\infty$ a.s., that is the random walk $S$  will not drift to $\infty$.
 As we have mentioned before, beside of the decay rate of $F$ the long-term behavior of the expectation $ \mathbf{E}[F(I_n)]$ also heavily depends on the speed at which the random walk $S$ approaches to $-\infty$.
 The fluctuation of $S$ is closely related to the following crucial quantity:
 \beqnn
 \varrho:= \inf_{\lambda\in[0,\theta_F]} \mathcal{L}_X(\lambda).
 \eeqnn
 It is obvious that  $\varrho\leq 1$ and the equality holds if and only if $\mathcal{D}_{\mathcal{L}_X}^+=\{ 0 \}$ or $\mathbf{E}[X]=0$; see Figure~\ref{Figure1}.
 In the sequel of this paper we always make the following assumption:
 \begin{assumption}
 	The infimum of $\mathcal{L}_X(\lambda)$ over $[0,\theta_F]$ can be attained, i.e. $ \Lambda \in \mathcal{D}_{\mathcal{L}_X}$ in which $\Lambda$ is defined in (\ref{Lambda}) and satisfies $\mathcal{L}_X(\Lambda)=\varrho$. 
 \end{assumption}
 
 According to the location of $\Lambda$ in $[0,\theta_F]$ and the asymptotic behavior of random walk $S$ under $\mathbf{P}^{(\Lambda)}$, six regimes arise for the expectation $\mathbf{E}[F(I_n)]$; see Figure~\ref{Figure1}. 
 In the next two subsections, we provide an explicit description for the asymptotic behavior of $\mathbf{E}[F(I_n)]$ in each regime.

 \subsubsection{The oscillating cases}
 
 We now provide the asymptotic results for the expectation $ \mathbf{E}[F(I_n)]$ with the random walk $S$ being oscillating under $\mathbf{P}^{(\Lambda)}$. We first consider the case in which $\Lambda=0$ and the well-known {\it Spitzer's condition} holds for $S$; see Figure~\ref{Figure1}(a).
 \begin{condition}\label{Con.Spitzer}
 	There exists a constant $\rho\in(0,1)$ such that as $n\to\infty$, 
 	\beqnn
 	\frac{1}{n}\sum_{k=1}^n\mathbf{P}( S_k>0 )  \to \rho.
 	\eeqnn
 \end{condition}
 
 Spitzer's condition is of key importance in fluctuation theory of random walks and it is equivalent to the convergence $\mathbf{P}( S_n>0 ) \to \rho$; see Theorem 1 in \cite{Doney1995}.
 Especially, all symmetric random walks satisfy Spitzer's condition with $\rho=1/2$.
 To show the exact decay rate of the expectation $ \mathbf{E}[F(I_n)]$, we need the positive, slowly varying function $\ell_1$  defined by
 \beqlb\label{ell.1}
 \ell_1(x) :=\frac{1}{\Gamma(\rho)}\exp\Big\{ \sum_{k=1}^\infty \frac{(1-1/x)^k}{k}\big(\mathbf{P}(S_k\geq 0) -\rho\big) \Big\},\quad x\geq 1,
 \eeqlb
 where $\Gamma $ is the Gamma function.

 \begin{theorem}\label{MainThm.01}
 	If $\theta_F> \Lambda =0$, under Condition~\ref{Con.F.LipCon} and~\ref{Con.Spitzer},  we have as $n\to\infty$,
 	\beqlb\label{MainThm.01.01}
 	\mathbf{E}[F(I_n)]\sim  C_{F,1} \cdot	\mathbf{P}(\tau_0^->n)\sim C_{F,1} \cdot n^{\rho-1} \ell_1(n),
 	\eeqlb
 	where the limit coefficient $C_{F,1}\in (0,\infty)$ is given by
 	\beqlb\label{MainThm.01.02}
 	C_{F,1}=\sum_{k=0}^\infty \mathbf{E}\big[  F\big(I_k + e^{-S_k} \cdot I^\uparrow_\infty) ; \sigma_k^- =k \big].
 	\eeqlb
 	
 \end{theorem}

 \begin{remark}\label{ConRateMeasureChange} 
 By the inequality above (2.22) in \cite{AfanasyevGeigerKerstingVatutin2005} with $\delta=\rho/2>0$, there exists a constant $C>0$ such that $ \liminf_{n\to\infty}n^{-\rho/2}S^{\uparrow}_n\geq C $ a.s., which directly induces $I^\uparrow_\infty<\infty$ a.s. and hence  $C_{F,1}>0$. 
 	Similarly, the other three random variables $I^\downarrow_\infty$, $\hat{I}^\uparrow_\infty$ and $\hat{I}^\downarrow_\infty$ are also finite almost surely.
 \end{remark}
 
 %

 We now turn to consider the asymptotics of  the expectation $ \mathbf{E}[F(I_n)]$ with $0< \Lambda<\theta_F$; see Figure~\ref{Figure1}(b).
 Let $\mathbf{P}^{(\Lambda)}$ be a probability measure defined in (\ref{MeasureChan}) with $\lambda_0=\Lambda$.
 Using the change of measure, we have
 \beqlb\label{MeasureChange}
 \mathbf{E}\big[F(I_n)\big]
 =   \varrho^n \cdot   \mathbf{E}^{(\Lambda)}
 \big[e^{-\Lambda S_n}F(I_n)\big].
 \eeqlb
 In contrast to the previous case, the long-term behavior of the expectation $\mathbf{E}\big[F(I_n)\big]$ depends not only on the sample paths with slowly decreasing local minimum but also on the sample paths with final value $S_n$ below a low level. In order to get an exact description for the distribution of $S_n$, we need the following condition.
 
 \begin{condition}\label{Con.AttractionDomain}
 	Under $\mathbf{P}^{(\Lambda)}$ the random walk $S$ is in the domain of attraction of a strictly stable law without centering with index $\alpha \in(0,2]$ and positivity parameter $\rho \in(0,1)$, we write $S\in \mathcal{D}_{\alpha,\rho}^{(\Lambda)}$.
 \end{condition}
 
 This condition is usually stronger than Spitzer's condition. Indeed, if $S\in \mathcal{D}_{\alpha,\rho}^{(\Lambda)}$ we have $\mathbf{P}^{(\Lambda)}(S_n>0)\to \rho$. For the converse,   Spitzer's condition usually does not imply a domain of attraction; readers may refer to \cite[p.380]{BinghamGoldieTeugels1987} for more details.  In particular,  $S\in\mathcal{D}_{\alpha,\rho}^{(\Lambda)}$ with $\alpha\in(1,2]$ if and only if Spitzer's condition holds for $S$ under $\mathbf{P}^{(\Lambda)}$, and in this case we always have $\rho=1-1/\alpha$; see Proposition~8.9.16 in \cite[p.384]{BinghamGoldieTeugels1987}.
 Moreover, Condition~\ref{Con.AttractionDomain} holds if and only if there exists a  function $\ell_2$ that is slowly varying at $\infty$ such that under $\mathbf{P}^{(\Lambda)}$,
 \beqnn
 \frac{S_n}{n^{1/\alpha}\ell_2(n)} \to Y_\alpha
 \eeqnn
 in distribution, where $Y_\alpha$ is a strictly stable random variable of parameter $\alpha$ and positivity parameter $\rho $.
 To simplify the following statements, we define
 \beqlb\label{SeqAn}
 A_n :=n^{-1-1/\alpha}/\ell_2(n),\quad n=1,2,\cdots.
 \eeqlb
 The following quantities are necessary to give an exact decay rate for the expectation $\mathbf{E}[F(I_n)]$. We write $\{g_\alpha(x):x\in\mathbb{R}\}$ for the probability density function of $Y_\alpha$.
 The fact that $V^{(\Lambda )}(x)=O(x)$ and $\hat{V}^{(\Lambda )}(x)=O(x)$ as $x\to\infty$ induces that the following two integrals are finite:
 \beqnn
 \mathcal{L}_{V^{(\Lambda )}}^{(\lambda)}(y):= \int_0^y e^{-\lambda z} V^{(\Lambda )}(z) dz \quad \mbox{and}\quad  \mathcal{L}_{\hat{V}^{(\Lambda )}}^{(\lambda)}(y):= \int_0^y e^{-\lambda z}\hat{V}^{(\Lambda )}(z) dz,\quad \lambda,y >0.
 \eeqnn
 Let $\mu_{V^{(\Lambda )}}^{(\lambda)}(dx)$ and $\mu_{\hat{V}^{(\Lambda )}}^{(\lambda)}(dx)$ be two probability measures on $\mathbb{R}_+$ defined by
 \beqlb\label{measureMu}
 \mu_{V^{(\Lambda )}}^{(\lambda)}(dx) := \frac{e^{-\lambda x} V^{(\Lambda )}(x) }{  \mathcal{L}_{V^{(\Lambda )}}^{(\lambda)}(\infty)}dx
 \quad \mbox{and}\quad
 \mu_{\hat{V}^{(\Lambda )}}^{(\lambda)}(dx)
 :=  \frac{e^{-\lambda x}\hat{V}^{(\Lambda )}(x)}{ \mathcal{L}_{\hat{V}^{(\Lambda )}}^{(\lambda)}(\infty) } dx.
 \eeqlb
 \begin{theorem} \label{MainThm.02}
 	If $\theta_F>\Lambda >0$, under Condition~\ref{Con.F.LipCon} and~\ref{Con.AttractionDomain},  we have as $n\to\infty$,
 	\beqlb\label{eqn.MainThm.02}
 	\mathbf{E}[F(I_n)]
 	\sim C_{F,2} \cdot \mathbf{P}(\tau_0^->n)
 	\sim  C_{F,2}\cdot g_\alpha(0) \mathcal{L}_{\hat{V}^{(\Lambda )}}^{(\Lambda )}(\infty)   \cdot  \varrho^n \cdot A_n ,
 	\eeqlb
 	where the limit coefficient\,\footnote{The positiveness of $C_{F,2}$ is a direct consequence of the finiteness of $I_\infty^{ \uparrow}$ and $ I_\infty^{ \downarrow}$; see Remark~\ref{ConRateMeasureChange}. } $C_{F,2}\in(0,\infty)$ is given by\,\footnote{$\mathbf{E}_{(x,y,z)}^{(\Lambda)}$ is the expectation of $(S,S^{\uparrow},S^{\downarrow})$ under $\mathbf{P}^{(\Lambda)}$ with initial state $(x,y,z)$. }
 	\beqnn
 	C_{F,2}	
 	\ar=\ar \sum_{k=0}^\infty  \Big[ \int_0^\infty  \mathbf{E}_{(0,0,y)}^{(\Lambda)} \big[e^{-\Lambda S_k} F\big(I_k+e^{-S_k}(I^{\uparrow}_\infty+e^{-y}+I^{\downarrow}_\infty)\big); \sigma_k^-=k \big] \mu_{\hat{V}^{(\Lambda )}}^{(\Lambda )}(dy)\cr
 	\ar\ar +  \int_0^\infty  \mathbf{E}_{(0,0,z)}^{(\Lambda)} \big[ e^{-\Lambda S_k} F\big(e^{z}(1+   I_k+ I_\infty^{\uparrow}  + I_\infty^{\downarrow})\big) ;\tau_0^->k  \big]  \frac{ e^{\Lambda z}  V^{(\Lambda)}(z) }{  \mathcal{L}_{\hat{V}^{(\Lambda)}}^{( \Lambda )}(\infty)   }dz   \Big]. 
 	\eeqnn 
 \end{theorem}

 We now continue to consider the asymptotics of  the expectation $ \mathbf{E}[F(I_n)]$  with $\theta_F=\Lambda>0$; see Figure~\ref{Figure1}(c).
 In this case, we find that the impact of sample paths with final value below a low level on the expectation on the right side of (\ref{MeasureChange}) is more complicated than that in the preceding case.
 In more detail,  if the random walk ends up below a low level at time $n$, we observe that its local minimum over the time interval $[0,n]$ is not only approximately equal to the final value but also prefers to be attained at the end of the time interval, and consequently its contribution to the expectation $ \mathbf{E}[F(I_n)]$ can be roughly represented as
 \beqnn
 e^{-\Lambda S_n} F(I_n) =O\big(  e^{-\Lambda S_n} ( I_n-I_{\sigma_n^--1} )^{-\Lambda}  \big) = O\big((n-\sigma_n^-)^{-\Lambda}\big).
 \eeqnn
 These make it difficult to provide an exact description for the decay rate of the expectation $ \mathbf{E}[F(I_n)]$.
 To keep away from these difficulties, we consider the aysmptotics of the expectation $\mathbf{E}[F(I_n)]$ with  $F$ satisfying Condition~\ref{Con.AsymPower}.
 Let $\hat\ell_1^{(\Lambda)}$ be a positive and slowly varying function defined by
 \beqlb\label{hatell.1}
 \hat\ell_1^{(\Lambda)}(x)\ar=\ar   \frac{1}{\Gamma(1-\rho)}\exp\Big\{ -\sum_{k=1}^\infty \frac{(1-1/x)^k}{k}\big( \mathbf{P}^{(\Lambda)}(S_k> 0) -\rho \big) \Big\},\quad x\geq 1 .
 \eeqlb
 \begin{theorem}\label{MainThm.03}
 	If $\theta_F=\Lambda>0$, under Condition~\ref{Con.AsymPower} and~\ref{Con.AttractionDomain},
 	we have as $n\to\infty$,
 	\beqnn
 	\mathbf{E}[F(I_n)]\sim K_0\cdot C_{F,3}\cdot  \varrho^n \cdot \mathbf{P}^{(\Lambda)}(\tau_0^+>n)
 	\sim K_0\cdot C_{F,3}\cdot  \varrho^n \cdot n^{-\rho}  \hat\ell_1^{(\Lambda)}(n),
 	\eeqnn
 	where the coefficient $C_{F,3} \in(0,\infty)$ is given by
 	\beqnn
 	C_{F,3}= \sum_{k=0}^\infty \mathbf{E}^{(\Lambda)}\big[  \big(1+\hat{I}_k + e^{-\hat{S}_k} \cdot \hat{I}^\uparrow_\infty)^{-\Lambda} ; \hat\sigma_k^- =k \big].
 	\eeqnn
 	
 \end{theorem}

 \subsubsection{The negative drift cases}
 
 We now provide the asymptotic results for the expectation $\mathbf{E}[F(I_n)]$ in which under $\mathbf{P}^{(\Lambda)}$ the random walk $S$ has negative drift, i.e.,
 \beqlb\label{Constant.a}
 a:=-\mathbf{E}^{(\Lambda)}[X]>0.
 \eeqlb
 We first consider the case with $\Lambda= \theta_F >0$ and $F$ satisfying Condition~\ref{Con.AsymPower}; see Figure~\ref{Figure1}(f).
 We notice that the expectation $\mathbf{E}[F(I_n)]$ is asymptotically equivalent to $K_0\cdot\mathbf{E}[I_n^{-\theta_F}]$.
 Using the change of measure and then the duality lemma, we have
 \beqnn
 \mathbf{E}[I_n^{-\Lambda}] \sim \mathbf{E}^{(\Lambda)}[(1+\hat{I}_{n-1})^{-\Lambda}] \cdot   \varrho^n.
 \eeqnn
 By Lemma~\ref{FirstLemma}, as $n\to\infty$ we have $\hat{I}_{n-1}\to \hat{I}_\infty<\infty$ a.s. under $\mathbf{P}^{(\Lambda)}$ and hence $\mathbf{E}^{(\Lambda)}[(1+\hat{I}_\infty)^{-\Lambda}] >0$.
 \begin{theorem}\label{MainThm.06} If $ \Lambda= \theta_F>0$, under Condition~\ref{Con.AsymPower},  we have as $n\to\infty$,
 	\beqnn
 	\mathbf{E}[F(I_n)] \sim  K_0\cdot\mathbf{E}^{(\Lambda)}[(1+\hat{I}_\infty)^{-\Lambda}] \cdot  \varrho^n .
 	\eeqnn
 	
 \end{theorem}
 
 We now provide the asymptotic results for the expectation $\mathbf{E}[F(I_n)]$ with $\theta_F>\Lambda\geq 0$ and
 the distribution of generic step $X$ always satisfying the following two regular variation assumptions:
 \begin{assumption}\label{Con.RegularV}
 	For some constant $\beta>1$, the tail-probability $	\mathbf{P}^{(\Lambda)}(X>x)$ is regularly varying with index $-\beta$, i.e., there exists a slowly varying function $\ell_3(x)$ at $\infty$ such that $	\mathbf{P}^{(\Lambda)}(X>x) \sim x^{-\beta}\ell_3(x)$ as $x\to\infty$,
 	
 \end{assumption}
 \begin{assumption}\label{Con.LocalRegularV}
 	For every $\delta>0$, we have $\mathbf{P}^{(\Lambda)}(X\in(x,x+\delta]) \sim \beta x^{-1-\beta}\ell_3(x)\cdot \delta$ as $x\to\infty$.
 \end{assumption}
 
 When $\Lambda=0$; see Figure~\ref{Figure1}(d), we first observe that the random walk prefers to attain the local minimum around  the first large step.
 Because of the negative drift, its local minimum will decrease slowly if there is a large step occurring at the beginning,  otherwise it drifts to $-\infty$ very fast and its contribution to the expectation $\mathbf{E}[F(I_n)]$ can be asymptotically ignored.
 By analyzing the contribution of sample paths before and after the first large step separately, we give in the next theorem, the exact decay rate of the expectation $\mathbf{E}[F(I_n)]$.

 \begin{theorem}\label{MainThm.04}
 	If $\theta_F>\Lambda=0$, assume  $F$ is non-increasing\,\footnote{Theorem~\ref{MainThm.04} still holds if  $F:=F_+-F_-$ with $F_+,F_-$ being bounded, non-increasing and satisfying Condition~\ref{Con.F.LipCon}.} and satisfies Condition~\ref{Con.F.LipCon}, we have as $n\to\infty$,
 	\beqlb\label{Eqn.MainThm.04.01}
 	\mathbf{E}[F(I_n)] \sim \frac{C_{F,4}}{\mathbf{E}[\tau_0^-]}\cdot \mathbf{P}(\tau_0^->n) \sim C_{F,4} \cdot \mathbf{P}(X\geq an),
 	\eeqlb
 	where the coefficient $C_{F,4}:= \sum_{k=1}^\infty \mathbf{E}[C_{F,4}(k)] \in (0,\infty)$ with
 	\beqnn
 	C_{F,4}(k):= \lim_{n\to\infty} \mathbf{E}\big[F(I_{k-1}+ e^{-S_{k-1}-X} \tilde{I}_{n})  \, \big|\, X\geq an,\mathscr{F}^S_{k-1}\big].
 	\eeqnn
 \end{theorem}	
 
 We now start to consider the case $\theta_F>\Lambda>0$; see Figure~\ref{Figure1}(e).
 Similarly as in the asymptotic analysis for the previous case, we see that the slow decreasing of its local minimum also stems from the early arrival of a large step.
 Moreover, (\ref{MeasureChange}) shows that the expectation $\mathbf{E}[F(I_n)]$ is also heavily effected by the final value $S_n$.
 However, we observe that the random walk would stay above a high level for a long time after the large step and hence its final value is more likely above a high level, which gives rise to its meager contribution to the expectation $\mathbf{E}[F(I_n)]$.
 Consequently, the main contribution to the expectation $\mathbf{E}[F(I_n)]$ is made by the sample paths with not only  an early large step but also small final value.
 In order to get an exact relationship between the local minimum and the final value, we need the following additional technic condition.
 
 \begin{condition}\label{MomentCon}
 	Assume that $\beta\neq 2$ and $\mathbf{E}[|X|^\kappa]<\infty$ for some $\kappa\in(1,2)$.
 \end{condition}
 Under this condition, 
 we have $n^{-1/\kappa}(S_n+an)$ converges to $0$ in distribution as $n\to\infty$ by the Kolmogorov-Marcinkiewicz-Zygmund law of large numbers; see Theorem 10.3 in \cite[p.311]{Gut2013}.
 These tell that the final value $S_n$ falls into a bounded interval around $0$
 if and only if the size of the early large step is about $an$.
 In order to simplify the notation, we define the following sequence
 \beqlb\label{SeqB}
 B_n=\frac{\beta}{an} \mathbf{P}^{(\Lambda)}(X\geq an),\quad n=1,2,\cdots.
 \eeqlb
 
 \begin{theorem}\label{MainThm.05}
 	If $\theta_F>\Lambda>0$, under Condition~\ref{Con.F.LipCon} and \ref{MomentCon}, we have as $n\to\infty$,
 	\beqnn
 	\mathbf{E}[F(I_n)] \sim \frac{C_{F,5}}{\mathcal{L}_{\hat{V}^{(\Lambda )}}^{(\Lambda )}(\infty)} \mathbf{P}(\tau_0^->n)  \sim  C_{F,5} \cdot \varrho^n \cdot B_n,
 	\eeqnn
 	where the coefficient  $C_{F,5}:=\sum_{k=1}^\infty C_{F,5}(k) \in(0,\infty)$ with
 	\beqnn
 	C_{F,5}(k):= \int_{-\infty}^\infty   \mathbf{E}^{(\Lambda)} \big[e^{-\Lambda S_{k}-\Lambda  z} F\big(I_{k}+ e^{-S_{k}-z}(1+\tilde{\hat{I}}_\infty)  \big)   \big]dz
 	\eeqnn
 	and $\tilde{\hat{I}}_\infty$ being an independent copy of $\hat{I}_\infty$.
 \end{theorem}


  \section{Auxiliary results for random walks}\label{Sec.AuxiliaryR}
 \setcounter{equation}{0}
 
 In this section we list as well as generalize some well-known asymptotic results for random walks under conditions and assumptions introduced in Section~\ref{Sec.MainR} with $\Lambda =0$, e.g., asymptotics of the first passage times, local probabilities conditioned to stay positive and conditional limit theorems. 
 
 \begin{remark}\label{Remark.PassageT.01}
 \red Notice that $V$ and $\hat{V}$ are c\`adl\`ag, all the following claims still hold  with $\tau_0^->n$, $\tau_0^+>n $, $V(x)$ and $\hat{V}(x)$ replaced by $L_n>0$, $M_n<0 $, $V(x-)$ and $\hat{V}(x-)$ respectively. 
 \end{remark}
 We first recall a useful asymptotic result for regularly varying sequences; readers can find it in \cite{AsmussenFossKorshunov2003}.
 
 \begin{lemma}\label{AuxiliaryLemma}
 	
 	Let $\{  b_n  \}$ be a regularly varying sequence. Consider two summable sequences $\{  f_n :n=0,1,\cdots \}$ and $\{ g_n:n=0,1,\cdots \}$ satisfying that  $f_n \sim c_1\cdot b_n$ and $g_n\sim c_2\cdot b_n$ with $c_1,c_2\geq 0$. We have as $n\to\infty$,
 	\beqnn
 	\sum_{k=0}^n f_{n-k}g_k \sim \Big(  c_1 \sum_{k=0}^\infty g_k+ c_2 \sum_{k=0}^\infty f_k \Big) \cdot b_n. 
 	\eeqnn

 	%
 	%

 \end{lemma}
 
 \subsection{Oscillating random walk}  
 
 The fluctuation theory for oscillating random walks satisfying Spitzer's condition has been well developed and abundant results have been gotten. Here we list a part of them that will be used in our following proofs.  
 Recall the two slowly varying functions $\ell_1$ and $\hat\ell_1$ defined in (\ref{ell.1}) and (\ref{hatell.1}) with $\Lambda=0$.  
 A simple calculation shows that $\ell_1(x)\hat\ell_1(x)$ converges to a positive constant as $x\to \infty$. 
 The next lemma comes from Theorem~8.9.12 in \cite[p.381]{BinghamGoldieTeugels1987} and Lemma~2.1 in  \cite{AfanasyevGeigerKerstingVatutin2005}.
 \begin{lemma}\label{Lemma.PassageT.01}
 	Under Condition~\ref{Con.Spitzer},  for every $x\geq 0$ we have as $n\to\infty$,
 	\beqnn
 	\mathbf{P}_x(\tau_0^->n) \sim V(x) n^{\rho-1} \ell_1(n)
 	\quad \mbox{and}\quad 
 	\mathbf{P}_{-x}(\tau_0^+>n) \sim \hat{V}(x) n^{-\rho} \hat\ell_1(n).
 	\eeqnn
 	Moreover, there exists a constant $C>0$ such that for all $x\geq 0$ and $n\geq 0$,
 	\beqnn
 	\mathbf{P}_x(\tau_0^->n) \leq C V(x) n^{\rho-1} \ell_1(n)
 	\quad \mbox{and}\quad 
 	\mathbf{P}_{-x}(\tau_0^+>n) \leq C \hat{V}(x) n^{-\rho} \hat\ell_1(n).
 	\eeqnn
 \end{lemma}
 The following result is a direct consequence of the proof for Lemma~2.2 in \cite{AfanasyevGeigerKerstingVatutin2005} with $u(x)=e^{-\lambda x}$. 
 \begin{lemma} \label{Lemma.Upper.01}
 	Under Condition~\ref{Con.Spitzer}, we have for every  $\lambda >0$, both of
 	the two sequences 
 	\beqnn
 	\mathbf{E}[e^{\lambda S_n}:\sigma_n^-=n]=\mathbf{E}[e^{\lambda S_n}:M_n< 0] 
 	\quad \mbox{and}\quad 
 	\mathbf{E}[e^{-\lambda S_n}:\sigma_n^+=n] =\mathbf{E}[e^{-\lambda S_n}:L_n>0]
 	\eeqnn
 	are summable and can be uniformly bounded by $c_0/n$ for some constant $c_0>0$.
 \end{lemma}
 
 The first conditional limit theorem for random walks in the next lemma was proved by Bertoin and Doney \cite{BertoinDoney1994} and the second one can be gotten immediately by using the duality lemma. 
 
 \begin{lemma}\label{Lemma.MeasureChange.02}
 	For $k\geq 1$, let $f$ be a bounded function on $\mathbb{R}^{k+1}$ and $f(S):=f(S_0,S_1,\cdots,S_k)$. 
 	Under Condition~\ref{Con.Spitzer}, for any $x\geq 0$ we have as $n\to\infty$, 
 	\beqnn
 	\mathbf{E}_x[f(S)\,|\, \tau_0^->n] 
 	\to\mathbf{E}_x [f(S^\uparrow )]
 	\quad\mbox{and}\quad 
 	\mathbf{E}_{-x}[f(S )\,|\, \tau_0^+>n] \to \mathbf{E}_{-x} [f(\hat{S}^\uparrow )],
 	\eeqnn
 	where $f(S^\uparrow):=f(S^\uparrow_0,S^\uparrow_1, \cdots, S^\uparrow_k)$  and $f(\hat{S}^\uparrow):= f(\hat{S}^\uparrow_0,\hat{S}^\uparrow_1, \cdots, \hat{S}^\uparrow_k)$.
 	
 \end{lemma}

 By Remark~\ref{ConRateMeasureChange}, we see that conditioned to stay positive the random walk $S$ drifts to $\infty$ a.s.
 To meet the needs of the following proofs for our main theorems, we provide in the following lemma, some large deviation estimates and uniform upper estimates for the final value of $S$ conditioned to stay positive. It can be proved by slightly extending and modifying the proofs for Proposition~2.1 and Corollary~2.4 in \cite{AfanasyevBoinghoffKerstingVatutin2012}, who considers the case with $\alpha\in(1,2]$. 
 
 \begin{lemma}\label{Lemma.CondLaplace.01}
 	Suppose $S\in \mathcal{D}_{\alpha,\rho}$ with $\alpha\in (0,2]$ and $\rho\in(0,1)$. 	For any $\lambda>0$, $x\geq 0$ and $y\in[0,\infty]$, we have as $n\to\infty$,
 	\beqnn
 	\mathbf{E}_x[e^{-\lambda S_n}; S_n\leq y, \tau_0^->n] \sim g_\alpha(0)  V(x)  \mathcal{L}_{\hat{V}}^{(\lambda)}(y) \cdot A_n ,\cr
 	\mathbf{E}_{-x}[e^{\lambda S_n};S_n\geq -y, \tau_0^+>n] \sim  g_\alpha(0)   \hat{V}(x) \mathcal{L}_{V}^{(\lambda)}(y) \cdot A_n.
 	\eeqnn
 	Moreover,  there exists a constant $C>0$ such that uniformly in $x,y\geq 0$ and $n\geq 1$,
 	\beqnn
 	\mathbf{E}_x[e^{-\lambda S_n}; S_n\leq y, \tau_0^->n] 
 	\ar\leq \ar C\cdot V(x) \mathcal{L}_{\hat{V}}^{(\lambda)}(y) \cdot A_n,\cr
 	\mathbf{E}_{-x}[e^{\lambda S_n};S_n\geq -y, \tau_0^+>n] 
 	\ar\leq \ar C\cdot  \hat{V}(x) \mathcal{L}_{V}^{(\lambda)}(y) \cdot A_n.
 	\eeqnn
 \end{lemma} 
 \proof Repeating the proof of  Proposition~2.1 in  \cite{AfanasyevBoinghoffKerstingVatutin2012} with
 	$a_n$ and $b_n$ replaced by $1/(nA_n) $ and $A_n$ respectively, we have 
 	\beqnn
 	\mathbf{E}_x[e^{-\lambda S_n}; \tau_0^->n] = \mathbf{E}_x[e^{-\lambda S_n};  L_n \geq 0] \sim g_\alpha(0)  V(x)  \mathcal{L}_{\hat{V}}^{(\lambda)}(\infty) \cdot A_n,
 	\eeqnn
 which induces that the finite measure $A_n^{-1} \cdot e^{-\lambda y}\cdot \mathbf{P}_x( S_n\in dy;\tau_0^->n )$ converges weakly to a finite measure on $\mathbb{R}_+$ with density $g_\alpha(0)  \cdot e^{-\lambda y}\cdot V(x)  \hat{V}(y)$ as $n\to\infty$. 
 Thus the first desired asymptotic result holds. 
   The second one can be proved similarly. 
   Similarly as in the proof of  Corollary~2.4 in \cite{AfanasyevBoinghoffKerstingVatutin2012} with $b_n$ replaced by $A_n$, we also can get the two upper bound estimates.
 	\qed 
 
 Let $C_0(\mathbb{R}^2)$ be the space of continuous functions on $\mathbb{R}^2$ vanishing at infinity. 
 The following lemma generalizes the conditional limit theorem  from Lemma~10 in \cite{Hirano1998}, which  considered the case of   $G(f(S),g_n(S))=f(S)g_n(S)$; see the following proof.
 
 \begin{lemma}\label{Lemma.MeasureChange.03}
 	Suppose $S\in \mathcal{D}_{\alpha,\rho}$ with $\alpha\in (0,2]$ and $\rho\in(0,1)$. 
 	Let $G\in C_0(\mathbb{R}^2)$ and $f, g$ be two bounded, continuous functions on $\mathbb{R}^{k+1}$ and $\mathbb{R}^{m+1}$ respectively for some $k,m\geq 1$. 
 	Let $f(S):=f(S_0,S_1,\cdots,S_k)$  and $g_n(S):=g(S_n,S_{n-1},\cdots,S_{n-m})$ for $n\geq m$.
 	For any $x\geq 0$ and $\lambda>0$ we have as $n\to\infty$,
 	\beqnn
 	\frac{\mathbf{E}_x [G(f(S),g_n(S))e^{-\lambda S_n} ; \tau_0^->n]}{\mathbf{E}_x[e^{-\lambda S_n} ; \tau_0^->n]}  \ar\to\ar \int_0^\infty\mathbf{E}_{(x,y)}[G(f(S^\uparrow),g(S^\downarrow))] \mu_{\hat{V}}^{(\lambda)}(dy),\cr
 	\frac{\mathbf{E}_{-x} [ G(f(S),g_n(S)) e^{\lambda S_n} ; \tau_0^+>n]}{\mathbf{E}_{-x}[e^{\lambda S_n} ; \tau_0^+>n]}  \ar\to\ar \int_0^\infty \mathbf{E}_{(-x,-y)}[G(f(\hat{S}^\uparrow),g(\hat{S}^\downarrow))]\mu_{V}^{(\lambda)}(dy),
 	\eeqnn
 	where $g(S^\downarrow):=g(S^\downarrow_0,S^\downarrow_1,\cdots,S^\downarrow_m)$, $g(\hat{S}^\downarrow):=g(\hat{S}^\downarrow_0,\hat{S}^\downarrow_1,\cdots,\hat{S}^\downarrow_m)$,  $f(S^\uparrow)$, $f(\hat{S}^\uparrow)$ are defined as in Lemma~\ref{Lemma.MeasureChange.02}
 	and $\mathbf{E}_{(x,y)}$ is the expectations of $(S^{\uparrow},S^\downarrow)$ or $(\hat{S}^{\uparrow},\hat{S}^\downarrow)$ with initial state $(x,y)$.
 \end{lemma} 
 \proof If $G(f(S),g_n(S))=f(S)g_n(S)$, the desired results can be gotten immediately  by repeating the proof for Lemma~10 in \cite{Hirano1998}. 
 The general results follow by the Stone-Weierstrass theorem.
 \qed

 %
 %

 \subsection{ Random walk with negative drift}\label{AuxiliaryNegD}
 
 We now list the asymptotic results for random walks with negative drift satisfying Assumption~\ref{Con.RegularV} and \ref{Con.LocalRegularV}.
 For every $x\geq 0 $, we  need the following important stopping time
 \beqnn 
 \mathcal{T}^x:= \inf\{k\geq 1: X_k> x  \} 
 \eeqnn  
 the first arrival of step with size larger than $x$.
 For each $k\geq 1$, a simple calculation together with Assumption~\ref{Con.RegularV} induces that the two events $\mathcal{T}^x=k$ and $X_k\geq x$ are asymptotically equivalent as $x\to\infty$, i.e., 
 \beqlb\label{Equvi.LargeJ}
 \mathbf{P}(\mathcal{T}^x=k | X_k\geq x)\sim \mathbf{P}(X_k\geq x| \mathcal{T}^x=k  )\to 1,
 \eeqlb
 and hence $ \mathbf{P}(  \mathcal{T}^x \leq k ) \sim k  \mathbf{P}( X\geq x  )$. The following well-known asymptotic results for the tail-probabilities of  the final value $S_n$ and the first passage time $\tau_0^-$ can be found in many literature, e.g., Theorem 5 in \cite{AsmussenFossKorshunov2003} and Theorem~2.2 in \cite{DenisovShneer2013}. 
 \begin{lemma}\label{Lemma.Passage.NegD}
 	For every $x\geq 0$, we have as $n\to\infty$,
 	\beqnn
 	\mathbf{P}( S_n\geq x )\sim n  \mathbf{P}( X\geq an  )
 	\quad \mbox{and}\quad 
 	\mathbf{P}_x(\tau_0^- \geq n  )   \sim \mathbf{E}_x[\tau_0^-] \mathbf{P}( X\geq an ) . 
 	\eeqnn
 \end{lemma}
 By this lemma and the Markov property, for each integer $k\geq 0$ we  have  as $n\to\infty$,
 \beqnn
 \mathbf{P}(\sigma_n^-=k)=   \mathbf{P}(\sigma_k^-=k) \cdot \mathbf{P}(\tau_0^-> n-k) \sim   \mathbf{P}(\sigma_k^-=k) \mathbf{E}[\tau_0^-] \cdot \mathbf{P}( X\geq an ) . 
 \eeqnn
 The next lemma comes from Remark 3.7 in \cite{Xu2021a}, which extends Theorem~3.2 in \cite{Durrett1980}. 
 It shows that an early large step is necessary to keep the random walk staying positive for a long time.
 \begin{lemma}\label{Lemma.CondBigStep}
 	For every integer $k\geq 1$, we have as $n\to\infty$,
 	\beqnn
 	\mathbf{P}(\mathcal{T}^{an}>n|\tau_0^->n)\to 0
 	\quad \mbox{and}\quad 
 	\mathbf{P}(\mathcal{T}^{an}\leq k|\tau_0^->n) \to  \mathbf{E}[\tau_0^- \wedge k]/\mathbf{E}[\tau_0^- ] .  
 	\eeqnn
 \end{lemma}
 
 Moreover, we also observe that the random walk will stay above a high level for a long time after the first large step.  
 This implies that the local minimum of the random walk should be attained near the early large step; see the next lemma and it can be proved by slightly modifying the proof for Lemma~4.6 in \cite{Xu2021a}. 
 \begin{lemma}\label{Lemma.BigStepMini}
 	Let $b\in(0,a]$ and pick an integer $k\geq 0$. For every $\epsilon>0$, there exist two integers  $n_0,t_0\geq 1$ such that for any $T\geq t_0$ and $n\geq n_0$,
 	\beqnn
 	\mathbf{P}(\mathcal{T}^{bn}\geq T, \sigma_n^- =k) \leq \epsilon \cdot \mathbf{P}(X\geq an). 
 	\eeqnn
 \end{lemma}
 \proof By the duality lemma, for $n>T>k$ we have 
 \beqnn
 \mathbf{P}(\mathcal{T}^{bn}\geq T, \sigma_n^- =k) 
 \ar\leq\ar  \mathbf{P}(\mathcal{T}^{an}\geq T-k, \tau_0^- >n-k )\cdot \mathbf{P}(  \sigma_k^- =k)\cr
 \ar\leq\ar \mathbf{P}(\mathcal{T}^{an}\geq T-k \,|\,   \tau_0^- >n-k )\cdot \mathbf{P}( \tau_0^- >n-k). 
 \eeqnn
 Using Lemma~\ref{Lemma.Passage.NegD} and then Lemma~\ref{Lemma.CondBigStep}, we have for some $C>0$,
 \beqnn
 \limsup_{n\to \infty} \frac{\mathbf{P}(\mathcal{T}^{bn}\geq T, \sigma_n^- =k)}{\mathbf{P}(X\geq an)} \leq C \cdot \big( 1-\mathbf{E}[\tau_0^- \wedge (T-k)]/\mathbf{E}[\tau_0^- ]  \big),
 \eeqnn
 which goes to $0$ as $T\to\infty$.
 \qed

 Recall the sequence $B_n$ defined in (\ref{SeqB}). 
 The next lemma provides some asymptotic results for joint local probabilities of the local maximum/minimum ($M_n/L_n$) and the final value ($S_n$). Its proof is similar to that of Lemma 3.13 and 3.14 in \cite{Xu2021a}. 
 \begin{lemma}\label{Lemma.LaplaceTrans}
 	For any $\lambda_1,\lambda_2>0$ and $x,y\in[0,\infty]$, we have as $n\to\infty$,
 	\beqnn
 	\mathbf{E}\big[e^{-\lambda_1 M_n-\lambda_2(M_n-S_n)}; M_n\leq x, M_n-S_n\leq y\big]  
 	\ar\sim\ar \mathcal{L}_V^{(\lambda_1)}(x) \mathcal{L}_{\hat{V}}^{(\lambda_2)}(y)  \cdot B_n , \cr
 	\mathbf{E}\big[e^{\lambda_1 L_n+\lambda_2(L_n-S_n)}; -L_n\leq x, S_n-L_n\leq y \big] 
 	\ar\sim\ar   \mathcal{L}_{\hat{V}}^{(\lambda_1)}(x)\mathcal{L}_V^{(\lambda_2)}(y) \cdot B_n .
 	\eeqnn
 \end{lemma}
 
 \proof Here we just provide a brief proof for the first desired result.  
 Similarly as in the proof for (12) in \cite{VatutinZheng2012}, we have as $n\to\infty$,
 \beqnn
 \mathbf{E}[e^{-\lambda_1 S_n}; S_n>0] \sim  \mathbf{E}[e^{\lambda_2 S_n}; S_n\leq 0] \sim \frac{\beta}{a} \mathbf{P}(X>an).
 \eeqnn 
 By the factorization identity of random walk; see Theorem 8.9.1 and 8.9.3 in \cite[p.376-377]{BinghamGoldieTeugels1987},  
 \beqnn
 1+\sum_{n=1}^\infty s^n \mathbf{E}\big[e^{-\lambda_1 M_n-\lambda_2(M_n-S_n)}\big] = \exp\Big\{ \sum_{n=1}^\infty \frac{s^n}{n} \big( \mathbf{E}[e^{-\lambda_1 S_n}; S_n>0] + \mathbf{E}[e^{\lambda_2 S_n}; S_n\leq 0]\big) \Big\} 
 \eeqnn
 for $ s\in(0,1)$ and 
 \beqnn
 \sum_{n=0}^\infty \mathbf{E}\big[e^{-\lambda_1 M_n-\lambda_2(M_n-S_n)}\big] 
 \ar=\ar \sum_{n=0}^\infty \mathbf{E}[e^{-\lambda_1 S_n}; \sigma_n^-=n] \cdot \sum_{n=0}^\infty \mathbf{E}[e^{-\lambda_2 S_n}; \sigma_n^-=0] = \mathcal{L}_V^{(\lambda_1)}(\infty) \mathcal{L}_{\hat{V}}^{(\lambda_2)}(\infty). 
 \eeqnn
 Applying Lemma 2.2.(2) in \cite{DenisovVatutinWachtel2014} to the foregoing equations,  we have as $n\to\infty$,
 \beqnn
 \mathbf{E}\big[e^{-\lambda_1 M_n-\lambda_2(M_n-S_n)}\big]  
 \ar\sim\ar \sum_{k=0}^\infty \mathbf{E}\big[e^{-\lambda_1 M_k-\lambda_2(M_k-S_k)}\big]  \cdot  B_n \sim  \mathcal{L}_V^{(\lambda_1)}(\infty) \mathcal{L}_{\hat{V}}^{(\lambda_2)}(\infty)\cdot B_n, 
 \eeqnn
 which induces that   for any $x,y\in(0,\infty)$,
 \beqnn
 \mathbf{P}(M_n\leq x,M_n-S_n\leq y) \sim   \int_0^x V(z_1)dz_1\int_0^y \hat{V}(z_2)dz_2  \cdot B_n
 \eeqnn
 as $n\to\infty$ and the desired result follows.
 \qed

 The next lemma gives several large deviation estimates for the random walk $S$ conditioned to stay positive/negative. 
 It comes from Remark~3.16 in \cite{Xu2021a} and can be proved like the previous lemma with the help of the Baxter identity
 \beqnn
 1+\sum_{n=1}^\infty s^n \mathbf{E}[e^{\lambda S_n}; \tau_0^+>n]= \exp\Big\{ \sum_{n=1}^\infty \frac{s^n}{n} \mathbf{E}[e^{\lambda S_n}; S_n<0]  \Big\}, \quad s\in(0,1),
 \eeqnn
 which can be found in Chapter XVIII.3 in \cite{Feller1971} or Chapter 8.9 in \cite{BinghamGoldieTeugels1987}. Here we omit the detailed proof.
 \begin{lemma}\label{Lemma.LargeDev.NegD} 	
 	For any $x, \lambda>0$ and $y\in[0,\infty]$, we have as $n\to\infty$,
 	\beqnn
 	\mathbf{E}_x[e^{-\lambda S_n}; S_n\leq y, \tau_0^->n] 
 	\ar\sim \ar  V(x) \mathcal{L}_{\hat{V}}^{(\lambda)}(y) \cdot B_n,\cr
 	\mathbf{E}_{-x}[e^{\lambda S_n};S_n\geq -y,\tau_0^+>n] 
 	\ar\sim \ar  \hat{V}(x)\mathcal{L}_{V}^{(\lambda)}(y) \cdot B_n.
 	\eeqnn
 \end{lemma}


 \section{Proof for Theorems~\ref{MainThm.01}-\ref{MainThm.03}}\label{Sec.Oscilliate}
 \setcounter{equation}{0}

 In this section we prove the asymptotic results for the expectation $\mathbf{E}[F(I_n)]$  with the random walk $S$ being oscillating under $\mathbf{P}^{(\Lambda )}$.
 Although the technical difficulties in the following proofs vary in different cases, one may find that it is the common key point to identify that sample paths with slowly decreasing local infimum make the main contribution to the expectation $\mathbf{E}[F(I_n)]$.
 For simplicity, we may always assume $\theta_F\in\mathcal{D}_F$. All the following proofs still work with  $\theta_F$ replaced by any $\theta \in \mathcal{D}_F\cap(\Lambda ,\infty)$.
 

 \subsection{Proof for Theorem~\ref{MainThm.01}}
 Under Condition~\ref{Con.Spitzer}, we first notice that sample paths with slowly decreasing local infimum would attain the local minimum at the beginning of the time interval.
 And then, we prove in the next proposition, that the contribution of sample paths with the local minimum being late attained to the expectation $\mathbf{E}[F(I_n)]$ can be asymptotically ignored.
 \begin{proposition}\label{Prop.MainThm01.01}
 	For every $\epsilon>0$, there exist two integers $k_0,n_0\geq 1$ such that for any $K\geq k_0$ and $n\geq n_0$,
 	\beqnn
 	\mathbf{E}\big[F(I_n); \sigma_n^-\geq K\big]\leq \epsilon \cdot \mathbf{P}( \tau_0^- >n ) .
 	\eeqnn
 	Moreover,  there exists a constant $C>0$ such that  for any $n\geq 0$,
 	\beqnn
 	\mathbf{E}\big[F(I_n) \big]\leq C \cdot \mathbf{P}( \tau_0^- >n ).
 	\eeqnn
 	
 	
 \end{proposition}
 \proof By Assumption~\ref{Con.F.UpperBound}, we first have $F(I_n) \leq C\cdot \exp\{\theta_F  L_n \}$ for some $C>0$. Using the Markov property of $S$, we have
 \beqnn
 \mathbf{E}\big[F(I_n); \sigma_n^-\geq K\big]  \ar\leq\ar  C\sum_{k=K}^n\mathbf{E}\big[e^{\theta_F S_k}; \sigma_k^-= k\big] \cdot \mathbf{P}( \tau_0^-> n-k ).
 \eeqnn
 From Lemma~\ref{Lemma.PassageT.01} and~\ref{Lemma.Upper.01}, we have $\mathbf{E}\big[e^{\theta_F S_n}; \sigma_n^-= n\big] = o(\mathbf{P}( \tau_0^-> n ))$ as $n\to\infty$. 
 	Applying Lemma~\ref{AuxiliaryLemma} to the foregoing partial sum, we have for large $n$, 
 \beqnn
 \mathbf{E}\big[F(I_n); \sigma_n^-\geq K \big]
 \ar\leq\ar C\sum_{k=K}^\infty \mathbf{E} \big[e^{\theta_F S_k}; \sigma_k^- =k \big] \cdot \mathbf{P}( \tau_0^-> n) .
 \eeqnn
 The desired two claims follow directly from the summability of the sequence $\mathbf{E}[e^{\theta_F S_k}; \sigma_k^- =k]$; see Lemma~\ref{Lemma.Upper.01}, and the fact that
 $\mathbf{E}\big[F(I_n); \sigma_n^- =0 \big] \leq C \mathbf{P}(\sigma_n^- =0) \sim C \mathbf{P}( \tau_0^-> n) $.
 \qed
 
 We now consider the contribution of sample paths with the local minimum attained at the beginning of the time interval, i.e., $\mathbf{E}[F(I_n); \sigma_n^-=k]$ for each fixed $k\geq 0$. Using the Markov property of $S$ and then the duality lemma,  we have
 \beqlb\label{Proof.MainThm.02}
 \mathbf{E}\big[F(I_n);\sigma_n^- =k\big]
 \ar=\ar \mathbf{E}\big[ \mathbf{E}\big[ F\big(I_k + e^{-S_k} \cdot \tilde{I}_{n-k} \big);\tilde\tau_0^- >n-k\,\big|\, \mathscr{F}^S_k \big]; \sigma_k^-=k \big].
 \eeqlb
 Thus it is a crucial step to analyze the asymptotics of the foregoing conditional expectation.
 
 
 \begin{proposition}\label{Prop.MainThm01.02}
 	Let $H$ be a bounded and Lipschitz continuous function on $(0,\infty)$. For every $\epsilon>0$, there exist two integers $k_0,n_0\geq 1$ such that for any $K\geq k_0$ and $n\geq n_0$,
 	\beqnn
 	\mathbf{E}\big[|H(I_n)-H(I_K)|;\tau_0^->n\big] =\mathbf{E}\big[|H(I_n)-H(I_K)|;L_n\geq 0\big] \leq \epsilon \cdot \mathbf{P}(\tau_0^->n).
 	\eeqnn
 	
 \end{proposition}
 \proof 
 By the duality lemma, we have $  \mathbf{E}\big[|H(I_n)-H(I_K)|;\tau_0^->n\big]=   \mathbf{E}\big[|H(I_n)-H(I_K)|; L_n= 0\big]$. By Remark~\ref{Remark.PassageT.01} and Lemma~\ref{Lemma.Passage.NegD}, the boundedness of $H$ induces that $  \mathbf{E} [|H(I_n)-H(I_K)|;L_n= 0 ] \leq C \mathbf{P}(L_n= 0)=o(\mathbf{P}(\tau_0^->n)) $ uniformly in $K$ as $n\to\infty$
 and hence it suffices to prove
 \beqnn
 \mathbf{E}\big[|H(I_n)-H(I_K)|\,;L_n>0\big] \leq \epsilon \cdot \mathbf{P}(\tau_0^->n).
 \eeqnn
 By the Lipschitz continuity of $H$,  there exists a constant $C>0$ such that
 \beqlb\label{Prop4.2.03}
 \mathbf{E}\big[ |H(I_n)-H(I_K)|\,; L_n> 0\big]
 \ar\leq\ar 
 C\sum_{k=K+1}^n \mathbf{E}\big[ e^{- S_k};L_n> 0\big].
 \eeqlb
 By the Markov property of $S$, we see that  $\mathbf{E}[ e^{- S_k}; L_n> 0 ]$ can be bounded by
 \beqnn
 \mathbf{E} \big[ e^{- S_k} \mathbf{P}_{S_k} (\tilde{L}_{n-k}>0)\,; L_k> 0 \big] .
 \eeqnn
 By the second result in Lemma~\ref{Lemma.PassageT.01}, we have $\mathbf{P}_{S_k} (\tilde{L}_{n-k}>0)\leq C V(S_k) \mathbf{P}(L_{n-k}>0) $ and hence 
 \beqlb\label{Prop4.2.04}
 \mathbf{E}[ e^{- S_k}; L_n> 0 ] \ar\leq\ar C \mathbf{E}\big[ e^{- S_k} V(S_k)\,; L_k> 0 \big]\cdot  \mathbf{P}( L_{n-k}>0 ).
 \eeqlb 
 The fact that $V(x)= O(x)$ as $x\to\infty$ shows $V(x)e^{- x} \leq C e^{- x/2}$ for all $x\geq 0$.
 From Lemma~\ref{Lemma.PassageT.01} and~\ref{Lemma.Upper.01}, we have for large $k$,
 	\beqnn
 	\mathbf{E}\big[ e^{- S_k} V(S_k)\,;L_k> 0\big]
 	\leq  C \cdot \mathbf{E}[ e^{- S_k/2}; L_k> 0]  
 	=o \big(\mathbf{P}( \tau_0^->k )\big).
 	\eeqnn
 Taking this and (\ref{Prop4.2.04}) back into (\ref{Prop4.2.03}) and then using Lemma~\ref{AuxiliaryLemma}, we have for large $n$,
 \beqnn
 \mathbf{E}\big[\left|H(I_n)-H(I_K)\right|; L_n>0\big]
 \ar\leq\ar C\sum_{k=K}^\infty \mathbf{E}\big[ e^{-S_k/2};   L_k> 0 \big] \cdot \mathbf{P}( \tau_0^-> n ),
 \eeqnn
 Here the summation above vanishes as $K\to\infty$; see Lemma~\ref{Lemma.Upper.01}, and the desired result follows immediately.
 \qed
 
 \begin{proposition}\label{Prop.MainThm01.03}
 	Let $H$ be a bounded and continuous function on $(0,\infty)$. For $K\geq 1$, we have
 	\beqnn
 	\lim_{n\to\infty} \mathbf{E}\big[ H(I_K)\,|\, \tau_0^->n \big] = \mathbf{E} \big[ H(I_K^\uparrow) \big]>0.
 	\eeqnn
 	Moreover,  the sequence $\mathbf{E}[ H(I_K^\uparrow)]$ converges to $\mathbf{E} [ H(I_\infty^\uparrow)]>0$ as $K\to\infty$. 
 \end{proposition}
 \proof The first claim follows from Lemma~\ref{Lemma.MeasureChange.02}. 
 For the second one, by  the dominated convergence theorem we have $\lim_{K\to\infty}\mathbf{E}[ H(I_K^\uparrow)]=\mathbf{E} [\lim_{K\to\infty} H(I_K^\uparrow) ]$. 	The continuity of $H$ and the fact that $I_K^\uparrow$ increases to $ I_\infty^\uparrow<\infty$ a.s.; see Remark~\ref{ConRateMeasureChange}, induce that $\mathbf{E} [\lim_{K\to\infty} H(I_K^\uparrow) ]=\mathbf{E} [ H(I_\infty^\uparrow) ]$.
 \qed

 %

 %
 
 {\it Proof for Theorem~\ref{MainThm.01}.}
 We first prove this theorem with $F$ being globally Lipschitz continuous on $(0,\infty)$.
 By Proposition~\ref{Prop.MainThm01.01},
 \beqlb\label{Proof.MainThm.01}
 \lim_{n\to\infty} \frac{\mathbf{E}[F(I_n)]}{\mathbf{P}(\tau_0^-\geq n)} \ar=\ar \lim_{K\to\infty} \lim_{n\to\infty} \frac{\mathbf{E}[F(I_n);\sigma_n^- \leq K]}{\mathbf{P}(\tau_0^-\geq n)}  =   \sum_{k=0}^\infty \lim_{n\to\infty} \frac{\mathbf{E}[F(I_n);\sigma_n^- =k]}{\mathbf{P}(\tau_0^-\geq n)}  .
 \quad
 \eeqlb
 Applying Proposition~\ref{Prop.MainThm01.02} with $H(x)= \mathbf{E} [ F (I_k + e^{-S_k} \cdot x ) \, |\,  \mathscr{F}^S_k   ]$ to the conditional expectation on the right side of (\ref{Proof.MainThm.02}) and then using the dominated convergence theorem, we have 
 \beqlb\label{eqn.5.63}
 \lefteqn{\lim_{n\to\infty} \frac{\mathbf{E}[F(I_n);\sigma_n^- =k]}{\mathbf{P}(\tau_0^-\geq n)}
 	= 
 	\lim_{K\to\infty} \lim_{n\to\infty} \frac{ \mathbf{E}\big[ \mathbf{E}\big[ F\big(I_k + e^{-S_k} \cdot \tilde{I}_{K} \big);\tilde\tau_0^- >n-k\,\big|\, \mathscr{F}^S_k \big]; \sigma_k^-=k \big]}{\mathbf{P}(\tau_0^-\geq n)}}\ar\ar\cr
 \ar=\ar \lim_{K\to\infty}\lim_{n\to\infty}\frac{\mathbf{P}(\tilde{\tau}_0^-\geq n-k)}{\mathbf{P}(\tau_0^-\geq n)}  \cdot \mathbf{E}\Big[ \mathbf{E}\big[ F\big(I_k + e^{-S_k} \cdot \tilde{I}_{K} \big)\,\big|\,\tilde\tau_0^- >n-k, \mathscr{F}^S_k \big]; \sigma_k^-=k \Big]\cr
 \ar=\ar  \lim_{n\to\infty}\frac{\mathbf{P}(\tilde{\tau}_0^-\geq n-k)}{\mathbf{P}(\tau_0^-\geq n)}
 \cdot \mathbf{E}\Big[ \lim_{K\to\infty}  \lim_{n\to\infty} \mathbf{E}\big[ F\big(I_k + e^{-S_k} \cdot \tilde{I}_{K} \big)\,\big|\,\tilde\tau_0^- >n-k, \mathscr{F}^S_k \big]; \sigma_k^-=k \Big] . \quad
 \eeqlb
 From Lemma~\ref{Lemma.PassageT.01}, we see the first limit on the right side of the second equality equals to $1$. 
 Applying Proposition~\ref{Prop.MainThm01.03} with $H(x)= F(I_k + e^{-S_k} \cdot x )$ to the second limit, we also have 
 \beqnn
 \lim_{K\to\infty} \lim_{n\to\infty} \mathbf{E}\big[ F\big(I_k + e^{-S_k} \cdot \tilde{I}_{K} \big)\,\big|\,\tilde\tau_0^- >n-k, \mathscr{F}^S \big]= \mathbf{E}\big[ F\big(I_k + e^{-S_k} \cdot I^\uparrow_\infty \big)\,\big|\,   \mathscr{F}^S \big].
 \eeqnn
 Taking this back into (\ref{eqn.5.63}) and then (\ref{Proof.MainThm.01}), we have
 \beqnn
 \lim_{n\to\infty} \frac{\mathbf{E}[F(I_n);\sigma_n^-=k] }{\mathbf{P}( \tau_0^- >n)}
 =  \mathbf{E} \big[ F\big(I_k + e^{-S_k} \cdot I^\uparrow_\infty \big); \sigma_k^-=k \big]>0.
 \eeqnn
 Hence the desired asymptotic equivalences in (\ref{MainThm.01.01}) hold and the limit coefficient $C_{F,1}$ is finite because of Proposition~\ref{Prop.MainThm01.01}. 
 For general $F$ satisfying  Condition~\ref{Con.F.LipCon} and $\delta>0$,
 we define $F_\delta(x)= F(x\vee \delta)$, which is  globally Lipschitz continuous.
 By Chebyshev's inequality,
 \beqnn
 \mathbf{E}\big[|F(I_n)-F_\delta(I_n)|\big] \leq  2C\cdot \mathbf{P}( I_n\leq \delta)
 \leq 2C\cdot \mathbf{P}(\sigma_n^-=0) + 2C \delta \cdot \mathbf{E}[I_n^{-1};\sigma_n^-\geq 1].
 \eeqnn
 By the duality lemma, we have $ \mathbf{P}(\sigma_n^-=0) \sim\mathbf{P}(\tau_0^-\geq n) $ for large $n$.
 Similarly as in Proposition~\ref{Prop.MainThm01.01}, we also can prove  $\mathbf{E}[I_n^{-1};\sigma_n^-\geq 1]\leq C \cdot \mathbf{P}(\tau_0^-\geq n)$ for any $n\geq 1$ and hence
 \beqnn
 \lim_{n\to\infty}\frac{\mathbf{E}[F(I_n)]}{\mathbf{P}(\tau_0^-\geq n)}  =\lim_{\delta\to 0}\lim_{n\to\infty}\frac{\mathbf{E}[F_\delta(I_n)]}{\mathbf{P}(\tau_0^-\geq n)}.
 \eeqnn
 The preceding result shows that  $\mathbf{E}[F_\delta(I_n)]/\mathbf{P}(\tau_0^-\geq n)\to C_{F_\delta,1}\in(0,\infty)$ as $n\to\infty $, where $C_{F_\delta,1}$ can be represented as the summation in (\ref{MainThm.01.02}) with $F$ replaced by $F_\delta$.
 Specially, when $F(x)=Cx^{-\theta_F}$ for any $C>0$, the monotone convergence theorem induces that  $ C_{F_\delta,1}\to C_{F,1} $ as $\delta\to 0+$.
 This can be extended to the general $F$ satisfying Assumption~\ref{Con.F.UpperBound} by the dominated convergence theorem.
 \qed
 
 %

 \subsection{Proof for Theorem~\ref{MainThm.02}}
 By (\ref{MeasureChange}), we observe that the asymptotics of the expectation  $\mathbf{E}^{(\Lambda )}[F(I_n)e^{-\Lambda S_n}]$ is mainly determined by sample paths with either the local minimum attained at the beginning of the time interval or the final value below a low level.
 Moreover, we also observe that if the random walk ends up at time $n$ below a low level, it will tend to attain the local minimum at the end of the time interval $[0,n]$.
 Precisely, the next proposition proves that the contribution of sample paths with the local minimum attained at neither the beginning nor the end of the time interval $[0,n]$ to the expectation $\mathbf{E}[F(I_n)]$ can be asymptotically ignored.

 \begin{proposition}\label{Proof.Prop.4.04}
 	For every $\epsilon>0$, there exist two integers $k_0,n_0\geq 1$ such that for any $K\geq k_0$ and $n\geq n_0$,
 	\beqnn
 	\mathbf{E}\big[F(I_n) ;\sigma_n^-\in[K,n-K] \big]\leq \epsilon \cdot \mathbf{P} (\tau_0^->n  ).
 	\eeqnn
 	Moreover, there exists a constant $C>0$ such that for any $n\geq 1$,
 	\beqnn
 	\mathbf{E}\big[F(I_n)\big] \leq C\cdot \mathbf{P} (\tau_0^->n  ).
 	\eeqnn
 \end{proposition}
 \proof  By the change of measure given in (\ref{MeasureChange}),  it suffices to prove
 \beqnn
 \mathbf{E}^{(\Lambda )}\big[F(I_n)e^{-\Lambda  S_n};\sigma_n^-\in[K,n-K] \big]\leq \epsilon \cdot \mathbf{E}^{(\Lambda )}\big[e^{-\Lambda  S_n};\tau_0^->n \big].
 \eeqnn
 By Assumption~\ref{Con.F.UpperBound} we have $\mathbf{E}^{(\Lambda )}[F(I_n)e^{-\Lambda  S_n};\sigma_n^-=k] \leq C\cdot \mathbf{E}^{(\Lambda )}[e^{\theta_{F}S_k-\Lambda S_n}; \sigma_n^-=k]$ for any $k\geq 1$.
 Using the Markov property, the independent increments of $S$ and then the duality lemma, we have
 \beqnn
 \mathbf{E}^{(\Lambda )}[e^{\theta_FS_k-\Lambda  S_n}; \sigma_n^-=k]
 \ar=\ar \mathbf{E}^{(\Lambda )}[e^{(\theta_F-\Lambda )S_k }; \sigma_k^-=k] \cdot \mathbf{E}^{(\Lambda )}[ e^{-\Lambda  S_{n-k}};\tau_0^->n-k ]\cr
 \ar=\ar \mathbf{E}^{(\Lambda )}[e^{(\theta_F-\Lambda )S_k }; \tau_0^+>k] \cdot \mathbf{E}^{(\Lambda )}[ e^{-\Lambda  S_{n-k}};\tau_0^->n-k ]
 \eeqnn
 and hence $ \mathbf{E}^{(\Lambda )}\big[F(I_n)e^{-\Lambda  S_n}; \sigma_n^-\in[K,n-K]\big]$ can be bounded by
 \beqnn
 C \sum_{k=K}^{n-K}  \mathbf{E}^{(\Lambda )}[e^{(\theta_F-\Lambda )S_k };\tau_0^+>k] \cdot \mathbf{E}^{(\Lambda )}[ e^{-\Lambda  S_{n-k}};\tau_0^->n-k ].
 \eeqnn
 Recall the sequence $A_n$ defined in (\ref{SeqAn}). 
 Using the first asymptotic equivalence in Lemma~\ref{Lemma.CondLaplace.01} with $\lambda = \Lambda$, $y=\infty$, the second one with $\lambda=\theta_F-\Lambda$, $y=\infty$ and then Lemma~\ref{AuxiliaryLemma} to the foregoing summation, we have for some $C>0$ and  large $n$, 
 \beqnn
 \lefteqn{\mathbf{E}^{(\Lambda )}\big[F(I_n)e^{-\Lambda  S_n}; \sigma_n^-\in[K,n-K]\big]}\ar\ar\cr
 \ar\leq \ar C \cdot A_n \cdot \sum_{k=K}^\infty \Big( \mathbf{E}^{(\Lambda )}[e^{(\theta_F-\Lambda )S_k };\tau_0^+>k]+ \mathbf{E}^{(\Lambda )}[ e^{-\Lambda  S_k};\tau_0^->k ] \Big) .
 \eeqnn 
 %
 Using Lemma~\ref{Lemma.CondLaplace.01} again, the last summation vanishes as $K\to\infty$ and then the first claim follows. The second claim follows from the first one and the fact that
 $\mathbf{E}^{(\Lambda )}[F(I_n)e^{-\Lambda  S_n};\sigma_n^-=0] \leq C \mathbf{E}^{(\Lambda )}[e^{-\Lambda  S_n};\tau_0^->n]$.
 \qed
 

 We now start to consider the contribution of sample paths with the local minimum attained early to the expectation $ \mathbf{E} [F(I_n)]$, i.e. $ \mathbf{E} [F(I_n) ;\sigma_n^- =k] $ for $k\geq 0$.
 As in the proof for Theorem~\ref{MainThm.01}, we first consider the conditional expectation on the right side of (\ref{Proof.MainThm.02}).
 Different to Proposition~\ref{Prop.MainThm01.02},
 the next proposition shows that conditioned on $\tau_0^->n$ for large $n$, the exponential functional $I_n$ can be well approximated by the sum of $I_J$ and $I_{n-J,n}:=I_n-I_{n-J}$ for large $J<n/2$.
 
 \begin{proposition}\label{Proof.Main02.Prop02}
 	Let $H$ be a bounded and Lipschitz continuous function. For every $\epsilon>0$, there exist two integers $j_0,n_0\geq 0$ such that for any $J\geq j_0$ and $n\geq n_0$,
 	%
 	\beqnn
 	\mathbf{E} \big[ |H(I_n)-H(I_{J}+I_{n-J,n})|; \tau_0^->n\big] \leq \epsilon \cdot \mathbf{P}(\tau_0^->n).
 	\eeqnn
 \end{proposition}
 \proof  By the Lipschitz continuity of $H$, we first have
 %
 \beqlb\label{Prop.MainThm.02.eqn07}
 \mathbf{E} \big[ |H(I_n)-H(I_{J}+I_{n-J,n})|\,\big|\, \tau_0^->n\big] 	\ar\leq\ar \sum_{j=J+1}^{n-J}	\mathbf{E} \big[ e^{-S_j }\,\big|\, \tau_0^->n \big] .
 \eeqlb
 By the change of measure, we have
 \beqlb\label{Prop.MainThm.02.eqn08}
 \mathbf{E} \big[ e^{-S_j } \,\big|\, \tau_0^->n \big]
 \ar=\ar \frac{\mathbf{E} [ e^{-S_j } ;\tau_0^->n ] }{\mathbf{P}(\tau_0^->n)  }
 = \frac{\mathbf{E}^{(\Lambda  )} [ e^{-S_j- \Lambda   S_n } ;\tau_0^->n ] }{\mathbf{E}^{(\Lambda)}[e^{-\Lambda  S_n};\tau_0^->n] }.
 \eeqlb
 Using the Markov property and the independent increments of $S$ to the numerator of the last fraction, we have
 \beqnn
 \mathbf{E}^{(\Lambda  )} [ e^{-S_j- \Lambda   S_n } ;\tau_0^->n ]
 = \mathbf{E}^{(\Lambda )}\big[ e^{- S_j } \mathbf{E}^{(\Lambda )}_{S_j}\big[ e^{-\Lambda  \tilde{S}_{n-j}}; \tilde{\tau}_0^->n-j\big] ; \tau_0^->j\big].
 \eeqnn
 Together with the fact that $V(x)=O(x)$ as $x\to\infty$, the second claim in Lemma~\ref{Lemma.CondLaplace.01}  implies that there exists a constant $C>0$ such that
 \beqnn
 \mathbf{E}^{(\Lambda )}_{S_j}\big[ e^{-\Lambda  \tilde{S}_{n-j}}; \tilde{\tau}_0^->n-j\big]
 \ar\leq\ar C \cdot S_j \cdot A_{n-j} 
 \eeqnn
 and hence by the fact that $xe^{- x} \leq Ce^{- x/2}$ for some $C>0$ and any $x\geq 0$,
 \beqnn
 \mathbf{E}^{(\Lambda  )} [ e^{-S_j- \Lambda   S_n } ;\tau_0^->n ]
 \ar\leq\ar C\mathbf{E}^{(\Lambda  )} \big[ e^{- S_j/2 } ;\tau_0^->j \big] \cdot A_{n-j}
 \leq C\cdot A_jA_{n-j}.
 \eeqnn
 Taking this back into (\ref{Prop.MainThm.02.eqn08}) and then using the first claim in Lemma~\ref{Lemma.CondLaplace.01} with $y=\infty$, we have for large $n$,
 \beqnn
 \mathbf{E} \big[ e^{-S_j } \,\big|\, \tau_0^->n \big] \ar\leq\ar
 \frac{C\cdot A_jA_{n-j}}{ A_n}.
 \eeqnn
 Taking this back into (\ref{Prop.MainThm.02.eqn07}) and then using Lemma~\ref{AuxiliaryLemma},  we have for some $C>0$ and large $n$,
 \beqnn
 \mathbf{E} \big[ |H(I_n)-H(I_{J}+I_{n-J,n})|\,\big|\, \tau_0^->n\big]
 \leq   C\sum_{j=J+1}^\infty A_j,
 \eeqnn
 which goes to $0$ as $J\to\infty$ and the desired result follows.
 \qed

 \begin{proposition}\label{Proof.Main02.Prop03}
 	%
 	Recall $\mu_{\hat{V}^{(\Lambda )}}^{(\Lambda )}(dy)$ defined in (\ref{measureMu}) with $\lambda=\Lambda$.	 Let $H$ be a bounded and continuous function on $(0,\infty)$ vanishing at infinity.
 	For each $J\geq 1$, we have as $n\to\infty$,
 	\beqnn
 	\mathbf{E} \big[ H(I_{J}+I_{n-J,n})\,\big|\, \tau_0^->n\big]
 	\ar\to\ar  \int_0^\infty	\mathbf{E}_{(0,y)}^{(\Lambda )} \big[ H(I^{\uparrow}_{J}+e^{-y}+I^{\downarrow}_{J-1})\big] \mu_{\hat{V}^{(\Lambda )}}^{(\Lambda )}(dy).
 	\eeqnn
 	Moreover, the limit coefficient converges as $J\to\infty$ to a finite limit given by
 	\beqnn
 	\int_0^\infty	\mathbf{E}_{(0,y)}^{(\Lambda )} \big[ H(I^{\uparrow}_\infty+e^{-y}+I^{\downarrow}_\infty)\big] \mu_{\hat{V}^{(\Lambda )}}^{(\Lambda )}(dy).
 	\eeqnn
 \end{proposition}
 \proof By the change of measure
 \beqnn
 \mathbf{E} \big[ H(I_{J}+I_{n-J,n})\,\big|\, \tau_0^->n\big]   \ar=\ar \frac{	\mathbf{E}^{(\Lambda )} \big[ H(I_{J}+I_{n-J,n})e^{-\Lambda  S_n}; \tau_0^->n\big]   }{	\mathbf{E}^{(\Lambda )} \big[ e^{-\Lambda  S_n}; \tau_0^->n\big]   }.
 \eeqnn
 Applying Lemma~\ref{Lemma.MeasureChange.03} to the foregoing expectations with $G(x_1,x_2)=H(x_1+x_2)$, $f(S)= I_J$ and $g_n(S)=I_{n-J,n}$, we can get the desired convergence immediately.
 The second claim follows by the dominated convergence theorem and the fact that $I^{\uparrow}_\infty, I^{\downarrow}_\infty<\infty$ a.s. under $\mathbf{P}^{(\Lambda )} $; see Remark~\ref{ConRateMeasureChange}.
 \qed

 \begin{corollary}\label{Proof.Corollary.01}
 	For each $k\geq 0$, we have as $n\to\infty$,
 	\beqnn
 	\frac{\mathbf{E} [F(I_n) ;\sigma_n^-=k]}{\mathbf{P}(\tau_0^->n)}
 	\ar\to\ar \int_0^\infty  	\mathbf{E}_{(0,0,y)}^{(\Lambda )} \big[ e^{-\Lambda S_k} F\big(I_k+e^{-S_k}(I^{\uparrow}_\infty+e^{-y}+I^{\downarrow}_\infty)\big); \sigma_k^-=k \big] \mu_{\hat{V}^{(\Lambda )}}^{(\Lambda )}(dy).
 	\eeqnn
 \end{corollary}
 \proof Applying Proposition~\ref{Proof.Main02.Prop02} together with  the continuity of $F$ on $(0,\infty)$ to the conditional expectation on the right side of (\ref{Proof.MainThm.02}), we have 
 \beqnn
 \lefteqn{\lim_{n\to\infty} \frac{\mathbf{E}\big[ F\big(I_k + e^{-S_k} \cdot \tilde{I}_{n-k} \big);\tilde\tau_0^- >n-k\,\big|\, \mathscr{F}^S_k \big]}{\mathbf{P}(\tau_0^->n)} }\ar\ar\cr
 \ar=\ar  \lim_{J\to\infty} \lim_{n\to\infty}\frac{  \mathbf{E}\big[ F\big(I_k + e^{-S_k} \cdot (\tilde{I}_{J}+\tilde{I}_{n-k-J,n-k}) \big);\tilde\tau_0^- >n-k\,\big|\, \mathscr{F}^S_k \big]}{\mathbf{P}(\tau_0^->n)}  \cr
 \ar=\ar \lim_{n\to\infty} \frac{\mathbf{P}(\tilde{\tau}_0^->n-k)}{\mathbf{P}(\tau_0^->n)}\cdot   \lim_{J\to\infty} \lim_{n\to\infty} \mathbf{E}\big[ F\big(I_k + e^{-S_k} \cdot (\tilde{I}_{J}+\tilde{I}_{n-k-J,n-k}) \big)\,\big|\, \tilde\tau_0^- >n-k,\mathscr{F}^S_k \big]  . 
 \eeqnn
 By the change of measure and Lemma~\ref{Lemma.CondLaplace.01} with $y=\infty$, the first limit on the right side of the last equality equals to $ \varrho^{-k}$. 
 By Proposition~\ref{Proof.Main02.Prop03} with $H(x)=F(I_k+ e^{-S_k}x)$, we have the second limit equals to 
 \beqnn
 \int_0^\infty  \mathbf{E}_{(0,y)}^{(\Lambda )}  \big[ F\big(I_k+e^{-S_k}(I^{ \uparrow}_\infty+e^{-y}+I^{ \downarrow}_\infty)\big)\,\big|\, \mathscr{F}^S_k\big] \mu_{\hat{V}^{(\Lambda )}}^{(\Lambda )}(dy).
 \eeqnn
 Taking these back into (\ref{Proof.MainThm.02}) and then using the change of measure, we have
 \beqnn
 \lefteqn{\lim_{n\to\infty}\frac{\mathbf{E} [F(I_n) ;\sigma_n^-=k]}{\mathbf{P}(\tau_0^->n)}
 	=  \mathbf{E}\Big[  \varrho^{-k} \int_0^\infty  \mathbf{E}_{(0,y)}^{(\Lambda )}  \big[ F\big(I_k+e^{-S_k}(I^{ \uparrow}_\infty+e^{-y}+I^{ \downarrow}_\infty)\big)\,\big|\, \mathscr{F}^S_k\big] \mu_{\hat{V}^{(\Lambda )}}^{(\Lambda )}(dy); \sigma_k^-=k   \Big] }\qquad \qquad\qquad\quad\ar\ar\cr
 \ar=\ar \mathbf{E}^{(\Lambda)}\Big[ e^{-\Lambda S_k}  \int_0^\infty  \mathbf{E}_{(0,y)}^{(\Lambda )}  \big[ F\big(I_k+e^{-S_k}(I^{ \uparrow}_\infty+e^{-y}+I^{ \downarrow}_\infty)\big)\,\big|\, \mathscr{F}^S_k\big] \mu_{\hat{V}^{(\Lambda )}}^{(\Lambda )}(dy); \sigma_k^-=k   \Big]
 \eeqnn
 and the desired result follows by Fubini's theorem.
 \qed 

 We now start to consider the impact of sample paths with the local minimum attained at the end of the time interval on the expectation $\mathbf{E}[F(I_n)]$.
 Applying the Markov property to
 $ \mathbf{E}[F(I_n);\sigma^-_n=n-k]$ for $n>k\geq 0$, we see it equals to
 \beqlb\label{eqn.4.66}
 \mathbf{E} \big[  \mathbf{E} \big[  F(I_{n-k}+ e^{-S_{n-k}} \tilde{I}_k);\sigma^-_{n-k}=n-k \,\big|\, \mathscr{F}^{\tilde{S}}_k \big]; \tilde{\tau}_0^->k \big].
 \eeqlb
 Like the previous argument, we first need to consider the aysmptotics of the conditional expectation.
 In the next proposition, we show that it is rarely contributed by sample paths with final value below a very low level.
 
 \begin{proposition}
 	For every $x\geq 0$ and $\epsilon>0$, there exist two integers $y_0,n_0\geq 1$ such that for any $y\geq y_0$ and $n\geq n_0$,
 	\beqnn
 	\mathbf{E}  \big[   F\big(I_{n}+x e^{-S_{n}} \big);S_n<-y,\sigma^-_{n}=n \big]\leq \epsilon\cdot \mathbf{P}\big(   \tau_0^->n  \big).
 	\eeqnn
 \end{proposition}
 \proof By the change of measure, we see that the desired  inequality  holds if and only if
 \beqnn
 \mathbf{E}^{(\Lambda )}  \big[ e^{-\Lambda  S_n}  F\big(I_{n}+x e^{-S_{n}} \big);S_n<-y,\sigma^-_{n}=n \big]
 \leq \epsilon\cdot\mathbf{E}^{(\Lambda )}\big[ e^{-\Lambda  S_n} ; \tau_0^->n  \big].
 \eeqnn
 Applying the duality lemma to the expectation on the left side of this inequality, we see it equals to
 \beqnn
 \mathbf{E}^{(\Lambda )}  \big[ e^{-\Lambda  S_n}  F\big( e^{-S_{n}}(1+\hat{I}_{n-1}+x ) \big);S_n<-y,M_n<0 \big],
 \eeqnn
 which can be bounded by $C	\mathbf{E}^{(\Lambda )}  \big[ e^{( \theta_F-\Lambda )S_n}  ;S_n<-y,M_n<0\big] $ because of the assumption that $F(x)\leq Cx^{-\theta_F}$.
 From Lemma~\ref{Lemma.CondLaplace.01}, there exists a constant $C>0$ such for any $n\geq 1$,
 \beqnn
 \lim_{n\to\infty} \frac{\mathbf{E}^{(\Lambda )}[ e^{( \theta_F-\Lambda )S_n}  ;S_n<-y,M_n<0]}{\mathbf{E}^{(\Lambda )} [ e^{-\Lambda S_n}  ;\tau_0^->n]} 
 \leq C\cdot \big|\mathcal{L}^{( \theta_F-\Lambda)}_V(\infty) - \mathcal{L}^{( \theta_F-\Lambda)}_V(y) \big|,
 \eeqnn
 which goes to $0$ as $y\to\infty$ and hence the desired result follows.
 \qed

 \begin{proposition}
 	Suppose $F$ is globally Lipschitz continuous on $(0,\infty)$. For each $x,y\geq0$ and $\epsilon>0$, there exist two integers $j_0,n_0>1$ such that for any $J\geq j_0$ and $n\geq n_0$,
 	\beqlb\label{eqn.74}
 	\mathbf{E}\big[ \big|F(I_n+x e^{-S_{n}}) -F(I_J+I_{n-J,n}+x e^{-S_{n}}) \big|; S_n\geq -y, \sigma_n^-=n  \big]\leq \epsilon \cdot \mathbf{P}(\tau_0^->n).
 	\eeqlb
 \end{proposition}
 \proof By the Lipschitz continuity of $F$, the expectation on the left side of (\ref{eqn.74}) can be bounded by
 \beqlb\label{eqn.4.72}
 C \mathbf{E}\big[I_{J+1,n-J}; S_n\geq -y, \sigma_n^-=n  \big]\ar=\ar C \sum_{j=J+1}^{n-J} \mathbf{E}\big[e^{-S_j}; S_n\geq -y, \sigma_n^-=n \big] .
 \eeqlb
 Applying the duality lemma and then the change of measure to $\mathbf{E}\big[e^{-S_j}; S_n\geq -y, \sigma_n^-=n \big] $,
 we see that it equals to
 \beqlb\label{eqn.4.73}
 \mathbf{E}\big[e^{S_{n-j}-S_n}; S_n\geq -y, M_n<0 \big]
 \ar=\ar \mathbf{E}^{(\Lambda )}\big[e^{S_{n-j}-(1+\Lambda )S_n}; S_n\geq -y, M_n<0 \big] \cdot \varrho^n,
 \eeqlb
 which can be bounded by $e^{(\kappa+1+\Lambda )y} \mathbf{E}^{(\Lambda )} [e^{S_{n-j}+\kappa S_n};  M_n<0] \cdot \varrho^n$ for any $\kappa>0$.
 By the Markov property and the second claim in Lemma~\ref{Lemma.CondLaplace.01},  there exists a constant $C>0$ such that
 \beqnn
 \mathbf{E}^{(\Lambda )}\big[e^{S_{n-j}+\kappa S_n};  M_n<0 \big] \ar=\ar \mathbf{E}^{(\Lambda )}\big[e^{S_{n-j} } \mathbf{E}^{(\Lambda )}_{S_{n-j}}[e^{\kappa \tilde{S}_j}; \tilde{M}_{j}<0];  M_{n-j}<0 \big] \cr
 \ar\leq\ar C \mathbf{E}^{(\Lambda )}\big[e^{S_{n-j} }S_{n-j};  M_{n-j}<0 \big] \cdot A_j
 \leq  C A_j A_{n-j}.
 \eeqnn
 Taking this and (\ref{eqn.4.73}) back into  (\ref{eqn.4.72}) and then using Lemma~\ref{Lemma.CondLaplace.01} with $\lambda=0$ and $y=\infty$, we have 
 for large $n$,
 \beqnn
 \frac{\mathbf{E}\big[I_{J+1,n-J}; S_n\geq -y, \sigma_n^-=n  \big] }{\mathbf{P}(\tau_0^->n)}
 \ar\leq\ar \frac{C}{A_n}\sum_{j=J+1}^{n-J}  A_jA_{n-j} , 
 \eeqnn
 which is asymptotically equivalent to $2C \sum_{j=J}^\infty  A_j$ as $n\to\infty$; see Lemma~\ref{AuxiliaryLemma}. Hence  the desired result follows as $J\to\infty$ because of the summability of the sequence $A_n$.
 \qed
 
 \begin{proposition}
 	For any  $x,y\geq 0$ and $J\geq 1$, we have as $n\to\infty$,
 	\beqnn
 	\lefteqn{	\frac{\mathbf{E}[F(I_J+I_{n-J,n} +xe^{-S_n} ); S_n\geq -y , \sigma_n^-=n]}{\mathbf{P}(\tau_0^->n)} }\ar\ar\cr
 	\ar\to \ar   \int_0^y \mathbf{E}_{(0,z)}^{(\Lambda )}\big[F\big(e^{z}\big(1+ x+ I_{J-1}^{\uparrow}  + I_{J-1}^{\downarrow}  \big)\big)   \big] \frac{ e^{ \Lambda  z} V^{(\Lambda)}(z)}{  \mathcal{L}_{\hat{V}^{(\Lambda)}}^{( \Lambda )}(\infty)   }dz  .
 	\eeqnn
 	Moreover, the limit coefficient converges as $J\to\infty$ and then $y\to\infty$ to
 	\beqnn 
 	\int_0^\infty \mathbf{E}_{(0,z)}^{(\Lambda )}\big[F\big(e^{z}\big(1+ x+ I_\infty^{\uparrow}  + I_\infty^{\downarrow}  \big)\big)   \big]\frac{ e^{\Lambda z} V^{(\Lambda)}(z) }{  \mathcal{L}_{\hat{V}^{(\Lambda)}}^{( \Lambda )}(\infty)   } dz.
 	\eeqnn
 \end{proposition}
 \proof For $\eta\in (0,\theta_F-\Lambda)$, by the change of measure we have 
 \beqlb\label{eqn.4.34}
 \lefteqn{ \frac{\mathbf{E}[F(I_J+I_{n-J,n} +xe^{-S_n} ); S_n\geq -y , \sigma_n^-=n]}{\mathbf{P}(\tau_0^->n)}}\ar\ar\cr
 \ar=\ar \frac{\mathbf{E}^{(\Lambda)}[F(I_J+I_{n-J,n} +xe^{-S_n} )e^{-\Lambda S_n}; S_n\geq -y , \sigma_n^-=n]}{\mathbf{E}^{(\Lambda)}[e^{-\Lambda S_n};\tau_0^->n]}\cr
 \ar=\ar\frac{\mathbf{E}^{(\Lambda)}[e^{\eta S_n};\tau_0^+>n]}{\mathbf{E}^{(\Lambda)}[e^{-\Lambda S_n};\tau_0^->n]} \cdot \frac{\mathbf{E}^{(\Lambda)}[F(I_J+I_{n-J,n} +xe^{-S_n} )e^{-\Lambda S_n}; S_n\geq -y , \sigma_n^-=n]}{\mathbf{E}^{(\Lambda)}[e^{\eta S_n};\tau_0^+>n]}.
 \eeqlb
 By Lemma~\ref{Lemma.CondLaplace.01} with $x=0$, $y=\infty$ and $\lambda =\Lambda$ or $\eta$, the first fraction on the right side of second equality in (\ref{eqn.4.34}) converges to  $	\mathcal{L}_{V^{(\Lambda )}}^{(\eta)}(\infty) / \mathcal{L}_{\hat{V}^{(\Lambda )}}^{( \Lambda  )}(\infty)  $ as $n\to\infty$. For the numerator of the second fraction, by the duality lemma we have 
 \beqnn
 \lefteqn{	\mathbf{E}^{(\Lambda )} \Big[F\big(I_J+I_{n-J,n} +xe^{-S_n} \big) e^{-\Lambda  S_n}; S_n\geq -y, \sigma_n^-=n\Big] }\ar\ar\cr
 \ar=\ar 	\mathbf{E}^{(\Lambda )} \Big[F\Big(e^{-S_n}\Big(1+x+\sum_{j=1}^{J-1} e^{S_j} +\sum_{i=n-J}^{n-1}e^{S_i}\Big)\Big) e^{-\Lambda  S_n} ; S_n\geq -y , M_n<0\Big] \cr
 \ar=\ar \mathbf{E}^{(\Lambda )}\Big[e^{-(\eta+\Lambda )S_n}F\Big(e^{-S_n}\Big(1+x+\sum_{j=1}^{J-1} e^{S_j} +\sum_{i=n-J}^{n-1}e^{S_i}\Big)\Big) \mathbf{1}_{\{S_n\geq -y\}} \cdot e^{\eta S_n}, M_n<0\Big].
 \eeqnn
 From Remark~\ref{Remark.PassageT.01} and Lemma~\ref{Lemma.CondLaplace.01}, we have $  \{M_n\leq 0\} = \{\tau_0^+>n\}$ and $\mathbf{E}^{(\Lambda)}[e^{\eta S_n};\tau_0^+>n]- \mathbf{E}^{(\Lambda)}[e^{-\eta S_n};M_n<0] =o(\mathbf{E}^{(\Lambda)}[e^{\eta S_n};\tau_0^+>n]) $. Hence as $n\to\infty$,
 	\beqnn
 	\lefteqn{\frac{\mathbf{E}^{(\Lambda)}[F(I_J+I_{n-J,n} +xe^{-S_n} )e^{-\Lambda S_n}; S_n\geq -y , \sigma_n^-=n]}{\mathbf{E}^{(\Lambda)}[e^{\eta S_n};\tau_0^+>n]} }\ar\ar\cr
 	\ar\sim\ar  \frac{\mathbf{E}^{(\Lambda )}\big[e^{-(\eta+\Lambda )S_n}F\big(e^{-S_n}\big(1+x+\sum_{j=1}^{J-1} e^{S_j} +\sum_{i=n-J}^{n-1}e^{S_i}\big)\big) \mathbf{1}_{\{S_n\geq -y\}} \cdot e^{\eta S_n}, \tau_0^+>n\big]}{\mathbf{E}^{(\Lambda)}[e^{\eta S_n};\tau_0^+>n]} . 
 	\eeqnn
 	Applying Lemma~\ref{Lemma.MeasureChange.03} with $\lambda=\eta$, $f(S)=\sum_{j=1}^{J-1} e^{S_j}$, $g_n(S)=(e^{-(\eta+\Lambda )S_n},e^{-S_n},\sum_{i=n-J}^{n-1}e^{S_i}, \mathbf{1}_{\{S_n\geq -y\}})$ and $G(y,(z_1,z_2,z_3,z_4))=z_1\cdot F(z_2(1+x+y+z_3))\cdot z_4$ to the preceding fraction, we have as $n\to\infty$,
 	\beqnn
 	\lefteqn{\frac{\mathbf{E}^{(\Lambda)}[F(I_J+I_{n-J,n} +xe^{-S_n} )e^{-\Lambda S_n}; S_n\geq -y , \sigma_n^-=n]}{\mathbf{E}^{(\Lambda)}[e^{\eta S_n};\tau_0^+>n]} }\ar\ar\cr
 	\ar\to \ar \int_0^y\mathbf{E}_{(0,-z)}^{(\Lambda )}\Big[e^{(\eta+\Lambda )z}F\Big(e^{z}\Big(1+x+ \sum_{j=1}^{J-1} e^{\hat{S}_j^{\uparrow}} + \sum_{i=1}^{J}e^{\hat{S}_i^{\downarrow}}\Big)\Big)   \Big] \mu_{V^{(\Lambda )}}^{(\eta)}(dz)\cr
 	\ar=\ar \int_0^y e^{(\eta+\Lambda )z}\mathbf{E}_{(0,z)}^{(\Lambda )}\Big[F\Big(e^{z}\big(1+ x+ I_{J-1}^{\uparrow}  + I_{J-1}^{\downarrow}  \big)\Big)   \Big] \mu_{V^{(\Lambda )}}^{(\eta)}(dz).
 	\eeqnn
 	Here the last equality follows from the fact that $\hat{S}^{\uparrow}=- S^{\uparrow}$ and $\hat{S}^{\downarrow}=- S^{\downarrow}$.
 	Taking these back into (\ref{eqn.4.34}) and using the definition of $\mu_{V^{(\Lambda )}}^{(\eta)}(dz)$; see (\ref{measureMu}), we have as $n\to\infty$, 
 	\beqnn
 	\lefteqn{ \frac{\mathbf{E}[F(I_J+I_{n-J,n} +xe^{-S_n} ); S_n\geq -y , \sigma_n^-=n]}{\mathbf{P}(\tau_0^->n)}}\ar\ar\cr
 	\ar\to\ar \frac{	\mathcal{L}_{V^{(\Lambda )}}^{(\eta)}(\infty) }{ \mathcal{L}_{\hat{V}^{(\Lambda )}}^{( \Lambda  )}(\infty)} \int_0^y e^{(\eta+\Lambda )z}\mathbf{E}_{(0,z)}^{(\Lambda )}\Big[F\Big(e^{z}\big(1+ x+ I_{J-1}^{\uparrow}  + I_{J-1}^{\downarrow}  \big)\Big)   \Big] \mu_{V^{(\Lambda )}}^{(\eta)}(dz)\cr
 	\ar =\ar  \int_0^y \mathbf{E}_{(0,z)}^{(\Lambda )}\Big[F\Big(e^{z}\big(1+ x+ I_{J-1}^{\uparrow}  + I_{J-1}^{\downarrow}  \big)\Big)   \Big] \frac{e^{\Lambda z} V^{(\Lambda )}(z) }{  \mathcal{L}_{\hat{V}^{(\Lambda )}}^{( \Lambda  )}(\infty)}dz.
 	\eeqnn
 	Here we have proved the first claim.
 The second one can be proved by using the dominated convergence theorem and the continuity of $F$ as $J\to\infty$, and then using the monotone convergence theorem as $y\to\infty$.
 \qed

 Applying the preceding two propositions to (\ref{eqn.4.66}), we can get the following corollary immediately.
 
 \begin{corollary}\label{Proof.Corollary.02}
 	For each $k\geq 0$, we have as $n\to\infty$
 	\beqnn
 	\frac{\mathbf{E} [F(I_n) ;\sigma_n^-=n-k]}{\mathbf{P}(\tau_0^->n)}	
 	\ar \to\ar
 	\int_0^\infty  \mathbf{E}_{(0,0,z)}^{(\Lambda )}\big[e^{-\Lambda S_k}F\big(e^{z}(1+   I_k+ I_\infty^{\uparrow}  + I_\infty^{\downarrow})\big) ;\tau_0^->k  \big] \frac{  e^{\Lambda z} V^{(\Lambda)}(z) }{  \mathcal{L}_{\hat{V}^{(\Lambda)}}^{( \Lambda )}(\infty)   } dz .
 	\eeqnn

 \end{corollary}

 {\it Proof for Theorem~\ref{MainThm.02}.} Here we just prove this theorem with $F$ being globally Lipschitz continuous.
 Proceeding as in the proof of Theorem~\ref{MainThm.01}, we can prove this theorem for general $F$.
 By Proposition~\ref{Proof.Prop.4.04}, we have 
  \beqnn
 	\lim_{n\to\infty} \frac{\mathbf{E}[F(I_n)]}{\mathbf{P}(\tau_0^->n)}
 	\ar=\ar 
 	\lim_{K\to\infty}\lim_{n\to\infty} \frac{\mathbf{E}[F(I_n); \sigma_n^-\in[K,n-K]]}{\mathbf{P}(\tau_0^->n)}\cr
 	\ar\ar +  \lim_{K\to\infty}\lim_{n\to\infty} \frac{\mathbf{E}[F(I_n); \sigma_n^-\in[0,K-1]]}{\mathbf{P}(\tau_0^->n)}\cr \ar\ar +\lim_{K\to\infty}\lim_{n\to\infty} \frac{\mathbf{E}[F(I_n); \sigma_n^-\in[n-K+1,n]]}{\mathbf{P}(\tau_0^->n)} \cr
 	\ar=\ar \sum_{k=0}^{\infty}\lim_{n\to\infty} \frac{\mathbf{E}[F(I_n);\sigma_n^-=k]}{\mathbf{P}(\tau_0^->n)} +  \sum_{k=0}^{\infty}  \lim_{n\to\infty} \frac{\mathbf{E}[F(I_n);\sigma_n^-=n-k]}{\mathbf{P}(\tau_0^->n)} .
 	\eeqnn 
 The desired result (\ref{eqn.MainThm.02}) follows directly from Corollary \ref{Proof.Corollary.01} and \ref{Proof.Corollary.02}. The finiteness of the limit coefficient $C_{F,2}$ can be gotten from Proposition~\ref{Proof.Prop.4.04}.
 \qed

 \subsection{Proof for Theorem~\ref{MainThm.03}} 
 
 %

 Condition~\ref{Con.AsymPower} and the boundedness of $F$ tell us that for every $\epsilon>0$, there exits a constant $x_0>0$ such that $|F(x)- K_0\cdot  (1+x)^{-\Lambda }|<\epsilon\cdot (1+x)^{-\Lambda }$ for any $x\geq x_0$.
 Thus we need to consider the asymptotic behavior of $\mathbf{E}[(1+I_n)^{-\Lambda }]$ at first. 
 By the change of measure and then duality lemma, 
 \beqnn
 \mathbf{E}[(1+I_n)^{-\Lambda }] =\varrho^n\cdot \mathbf{E}^{(\Lambda )}\big[(1+I_n)^{-\Lambda }e^{-\Lambda S_n}\big]
 =\varrho^n\cdot \mathbf{E}^{(\Lambda )}\big[(1+\hat{I}_{n})^{-\Lambda } \big] .
 \eeqnn
 According to our previous argument, the dual process $\hat{S}$ satisfies Spitzer's condition with positivity parameter $1-\rho$.
 Using Theorem~\ref{MainThm.01} with $F(x)= (1+x)^{-\Lambda }$, we have as $n\to\infty$,
 \beqnn
 \mathbf{E}^{(\Lambda )}\big[(1+\hat{I}_{n})^{-\Lambda } \big]
 \ar\sim\ar  C_{F,3} \cdot \mathbf{P}^{(\Lambda )}(\hat{\tau}_0^->n)= C_{F,3}\cdot \mathbf{P}^{(\Lambda )}(\tau_0^+>n).
 \eeqnn
 
 It remains to prove that $ \mathbf{E}^{(\Lambda )}\big[|F(I_n)-K_0(1+I_n)^{-\Lambda }|e^{-\Lambda S_n}\big]=o\big(\mathbf{P}^{(\Lambda )}(\tau_0^+>n)\big)$ as $n\to\infty$. From  the previous result, we first have for large $n$,
 \beqnn
 \mathbf{E}^{(\Lambda )}\big[|F(I_n)-K_0(1+I_n)^{-\Lambda }|e^{-\Lambda S_n}; I_n>x_0\big]
 \leq \epsilon\cdot \mathbf{E}^{(\Lambda )}\big[(1+I_n)^{-\Lambda }e^{-\Lambda S_n}\big]  \leq \epsilon\cdot C \cdot \mathbf{P}^{(\Lambda )}(\tau_0^+>n).
 \eeqnn
 On the other hand, the boundedness of $F$ implies that $F(x)\leq C(1+x)^{-\Lambda}$ for some $C>0$ and $x\geq 0$. By Chebyshev's inequality we also have
 \beqnn
 \mathbf{E}^{(\Lambda )}\big[|F(I_n)-K_0(1+I_n)^{-\Lambda }|e^{-\Lambda S_n}; I_n\leq x_0\big] \ar\leq\ar C(1+x_0)\mathbf{E}^{(\Lambda )}\big[(1+I_n)^{-\Lambda -1}e^{-\Lambda S_n}\big].  
 \eeqnn
 Notice that $(1+I_n)^{-\Lambda -1} \leq e^{(\Lambda+1)S_n}$ a.s. for any $n\geq 1$. 
 For each $k\geq 0$, by the Markov property of $S$ we have as $n\to\infty$, 
 \beqnn
 \mathbf{E}^{(\Lambda )}\big[(1+I_n)^{-\Lambda -1}e^{-\Lambda S_n};\sigma_n^-=k\big]
 \ar\leq \ar \mathbf{E}^{(\Lambda )}\big[e^{(\Lambda +1)S_k-\Lambda S_n};\sigma_n^-=k\big]\cr
 \ar=\ar \mathbf{E}^{(\Lambda )}\big[e^{ S_k};\sigma_k^-=k\big] \mathbf{E}^{(\Lambda )}\big[e^{-\Lambda S_{n-k}};\tau_0^->n-k\big].
 \eeqnn
 By the duality lemma, Remark~\ref{Remark.PassageT.01} and Lemma~\ref{Lemma.CondLaplace.01}, we have $\mathbf{E}^{(\Lambda )}\big[e^{ S_k};\sigma_k^-=k\big] =\mathbf{E}^{(\Lambda )}\big[e^{ S_k};M_k<0\big]\sim \mathbf{E}^{(\Lambda )}\big[e^{ S_k};\tau_0^+>k\big]$ as $k\to\infty$. 
 Applying Lemma~\ref{AuxiliaryLemma} together with Lemma~\ref{Lemma.CondLaplace.01}, we have
 \beqnn
 \mathbf{E}^{(\Lambda )}\big[(1+I_n)^{-\Lambda - 1}e^{-\Lambda S_n} \big]
 \ar=\ar   \sum_{k=0}^n\mathbf{E}^{(\Lambda )}\big[(1+I_n)^{-\Lambda - 1}e^{-\Lambda S_n};\sigma_n^-=k\big]\cr
 \ar\leq\ar \sum_{k=0}^n  \mathbf{E}^{(\Lambda )}\big[e^{ S_k};\tau_0^+>k\big] \mathbf{E}^{(\Lambda )}\big[e^{-\Lambda S_{n-k}};\tau_0^->n-k\big] = O(A_n)  ,
 \eeqnn
 which is $ o(\mathbf{P}^{(\Lambda )}(\tau_0^+>n))$; see (\ref{SeqAn}) and Lemma~\ref{Lemma.PassageT.01}.
 Putting all estimates above together, we can get the desired result immediately.
 \qed
 

 \section{Proof for Theorems~\ref{MainThm.06}-\ref{MainThm.05}} \label{Sec.NegD}
 
 In this section we prove the asymptotic results for the expectation $\mathbf{E}[F(I_n)]$ with $S$ drifting to $-\infty$ under $\mathbf{P}^{(\Lambda )}$.
 Different to the oscillating case, we observe that the slow decreasing of the local minimum of $S$ usually results from an early large step.
 For simplicity, we again assume $\theta_F\in\mathcal{D}_F$.
 
 \subsection{Proof for Theorem~\ref{MainThm.06}}
 
 By the change of measure and then  Condition~\ref{Con.AsymPower}, we have 
 \beqnn 
 \varrho^{-n}\cdot \mathbf{E} [ F(I_n)]= \mathbf{E}^{(\Lambda )}\big[e^{-\Lambda  S_n} \cdot F(I_n)\big]  = \mathbf{E}^{(\Lambda )} \big[e^{-\Lambda  S_n} (F(I_n) -K_0 I_n^{-\Lambda })\big]   + \mathbf{E}^{(\Lambda )}[e^{-\Lambda  S_n} K_0 I_n^{-\Lambda }] . 
 \eeqnn
 By the duality lemma and then Lemma~\ref{FirstLemma},  the last expectation equals to $ K_0 \mathbf{E}^{(\Lambda )}[( 1+ \hat{I}_{\red n-1})^{-\Lambda }] $, which converges to $K_0 \mathbf{E}^{(\Lambda )}[( 1+ \hat{I}_\infty)^{-\Lambda }]>0$ as $n\to\infty$. 
 Thus it suffices to prove $ \mathbf{E}^{(\Lambda )}[e^{-\Lambda  S_n} |F(I_n) -K_0 I_n^{-\Lambda }|]\to 0$ as $n\to\infty$. By Condition~\ref{Con.AsymPower}, for any $\epsilon>0$, there exists a constant $x_0>0$ such that for any $n\geq 1$ and $x\geq x_0$,
 \beqnn
 \mathbf{E}^{(\Lambda )}\big[e^{-\Lambda  S_n}|F(I_n)-K_0I_n^{-\Lambda }|\,; I_n\geq x\big]
 \ar\leq\ar
 \epsilon\cdot \mathbf{E}^{(\Lambda )} \big[e^{-\Lambda  S_n} I_n^{-\Lambda }\big]  \leq \epsilon.
 \eeqnn
 Here the last inequality follows from the fact that $e^{-\Lambda  S_n} I_n^{-\Lambda } \leq 1$. 
 On the other hand, by Assumption~\ref{Con.F.UpperBound}, we have $e^{-\Lambda  S_n}|F(I_n)-K_0 I_n^{-\Lambda }|\leq C e^{-\Lambda  S_n} I_n^{-\Lambda } \leq C$ for some $C>0$ independent of $S$. Hence for any $K>0$,
 \beqnn
 \mathbf{E}^{(\Lambda )}\big[e^{-\Lambda  S_n}|F(I_n)-K_0 I_n^{-\Lambda }|; I_n < x\big]
 \ar\leq\ar
 C  \mathbf{P}^{(\Lambda )}(S_n\geq -K) +  C \mathbf{P}^{(\Lambda )}( S_n< -K, I_n < x).
 \eeqnn
 Since $\mathbf{P}^{(\Lambda )}(S_n\to -\infty)=1$, the first probability on the right side of this inequality vanishes as $n\to\infty$.
 By Chebyshev's inequality and $ I_n^{-1}\leq e^{S_n}$, the second one can be bounded by
 \beqnn
 x\cdot \mathbf{E}^{(\Lambda )} \big[ I_n^{ -1};S_n< - K  \big]
 \ar\leq\ar  x\cdot \mathbf{E}^{(\Lambda )} \big[e^{S_n}; S_n< - K  \big]  \leq x e^{-K},
 \eeqnn
 which vanishes as $K\to\infty$ and the desired result follows by putting all estimates above together.
 
 \subsection{Proof for Theorem~\ref{MainThm.04}}
 \setcounter{equation}{0}
 
 As we have mentioned before, since the slow decreasing of the local minimum of $S$ usually results from an early large step,  the contribution of sample paths with late arrival of the first large step, $\mathbf{E}[F(I_n); \mathcal{T}^{an}>T]$ for large $T$, can be asymptotically ignored.
 By Lemma~\ref{Lemma.BigStepMini}, the local minimum is usually not far from the first large step.
 Thus for any $\epsilon>0$ and $K\geq 1$, there exist two integers $t_0,n_0>0$ such that for any $T>t_0$ and $n\geq n_0$,
 \beqlb\label{eqn.5.21}
 \mathbf{E} \big[F(I_n); \mathcal{T}^{an}>T, \sigma_n^-\leq K \big]  \leq  C \cdot \mathbf{P}(\mathcal{T}^{an}>T, \sigma_n^-\leq K)  \leq \epsilon \cdot \mathbf{P}(\tau_0^->n ).
 \eeqlb
 Here these two inequalities follow from the boundedness of $F$ and Lemma~\ref{Lemma.BigStepMini} respectively.  
 On the other hand,  in the next proposition we show that the contribution of sample paths with local minimum attained late,   $\mathbf{E}[F(I_n);\mathcal{T}^{an}>T, \sigma_n^-> K] \leq  \mathbf{E}[F(I_n) ; \sigma_n^-> K]$ for large $K$, also can be asymptotically ignored.
 
 \begin{proposition} \label{Prop.MainThm.04.01}
 	For every $\epsilon>0$, there exist two integers $k_0,n_0\geq 1$ such that for any $K\geq k_0$ and $n\geq n_0$,
 	\beqnn
 	\mathbf{E}\big[F(I_n);\sigma^-_n\geq K \big] \ar\leq\ar \epsilon\cdot \mathbf{P}(\tau_0^->n ).
 	\eeqnn
 	Moreover, there exists a constant $C>0$ such that for any $n\geq 1$,
 	\beqnn
 	\mathbf{E} \big[F(I_n) \big] \leq C\cdot \mathbf{P}(\tau_0^->n ).
 	\eeqnn
 \end{proposition}
 \proof  By Assumption~\ref{Con.F.UpperBound},  we first have $F(I_n) \leq   C \exp\{\theta_F L_n\}$ for some $C>0$ independent of $n$ and hence
 \beqlb\label{eqn.3.11}
 \mathbf{E} \big[F(I_n);\sigma_n^-\geq K \big]
 \leq  C \sum_{k=K}^n\mathbf{E} \big[e^{\theta_F S_k};\sigma_n^-=k \big]  \ar=\ar C \sum_{k=K}^n\mathbf{E} \big[e^{\theta_F S_k};\sigma_k^-=k \big] \mathbf{P}(\tau_0^->n-k ). 
 \eeqlb
 Here the last  equality  follows from the Markov property of $S$. 
 By the duality lemma and Remark~\ref{Remark.PassageT.01}, we also have $\mathbf{E}[e^{\theta_F S_n};\sigma_n^-=n] =\mathbf{E}[e^{\theta_F S_n};M_n<0] \sim \mathbf{E}[e^{\theta_F S_n};\tau_0^+>n]$ as $n\to\infty$. 
 Moreover, by Lemma~\ref{Lemma.Passage.NegD} with $x=0$ and Lemma~\ref{Lemma.LargeDev.NegD} with $y=\infty$, we have $\mathbf{E}[e^{\theta_F S_n};\tau_0^+>n]=o(\mathbf{P}(\tau_0^->n))$ as $n\to\infty$. 
 Applying Lemma~\ref{AuxiliaryLemma} to the last sum in (\ref{eqn.3.11}) and then using  Lemma~\ref{Lemma.LargeDev.NegD} with $y=\infty$ again,  there exists  a constant $C>0$ such that as $n\to\infty$, 
 \beqnn
 \sum_{k=K}^n\mathbf{E} \big[e^{\theta_F S_k};\sigma_k^-=k \big] \mathbf{P}(\tau_0^->n-k )
 \ar\sim\ar \mathbf{P}( \tau_0^->n )\cdot  \sum_{k=K}^\infty\mathbf{E} \big[e^{-\theta_F S_k};\tau_0^->k \big] \cr
 \ar\leq\ar \mathbf{P}(\tau_0^->n )\cdot \sum_{k=K}^\infty \frac{C}{k}\mathbf{P}( X\geq ak).
 \eeqnn 
 Notice that the sequence $ k^{-1}\cdot\mathbf{P}(X\geq ak ) $ is $(-\beta-1)$-regularly varying and summable, the two desired two upper bounds follow immediately. 
 \qed

 We now turn to consider the contribution of sample paths with early large step to the expectation (\ref{ExpRW}), i.e. $ \mathbf{E}[F(I_n);\mathcal{T}^{an}=k]$  for each $k\geq 1$. By (\ref{Equvi.LargeJ}) and Lemma~\ref{Lemma.Passage.NegD},   we have for large $n$,
 \beqlb\label{eqn.5.23}
 \frac{\mathbf{E}[F(I_n);\mathcal{T}^{an}=k]}{\mathbf{P}(\tau_0^->n)} \sim  \frac{\mathbf{E}[F(I_n);X_k\geq an]}{\mathbf{P}(\tau_0^->n)}.
 \eeqlb
 By the Markov property and the independent increments of $S$, we have 
 \beqlb\label{eqn.6.83}
 \frac{\mathbf{E}[F(I_n);X_k\geq an]}{\mathbf{P}(X\geq an)}\ar=\ar \frac{1}{\mathbf{P}(X\geq an)}\cdot \mathbf{E}\big[ \mathbf{E}[F(I_{k-1}+ e^{-S_{k-1}-X_k}\tilde{I}_{n-k}); X_k\geq an|\mathscr{F}^S_{k-1}]\big]\cr
 \ar\ar\cr
 \ar=\ar \mathbf{E}\big[ \mathbf{E}[F(I_{k-1}+ e^{-S_{k-1}-X}\tilde{I}_{n-k})|X\geq an,\mathscr{F}^S_{k-1}]\big]  .
 \eeqlb
 Hence we need to consider the asymptotics of the last conditional expectation at first.
 \begin{proposition}\label{Proposition6.02}
 	Let $H$ be a bounded and Lipschitz continuous function on $(0,\infty)$.
 	For every $\epsilon>0$, there exist two integers $k_0,n_0\geq 1$ such that for any $K\geq k_0$ and $n\geq n_0$,
 	\beqnn
 	\mathbf{E}\big[|H(e^{-X}I_n)-H(e^{-X}I_K)|; X\geq an\big] \leq \epsilon \cdot \mathbf{P}(X\geq an).
 	\eeqnn
 	
 \end{proposition}
 \proof Let $I_{K,n}:=I_n-I_K$ for $0\leq K\leq n$. The boundedness and Lipschitz continuity of $H$ induce that there exists a constant $C>0$ such that $|H(x)-H(y)|\leq C (1\wedge |x-y|)$ for any $x,y>0$. Hence uniformly in $n\geq K\geq 1$, 
 \beqnn
 \mathbf{E}\big[|H(e^{-X}I_n)-H(e^{-X}I_K)| ; X\geq an\big]
 \ar\leq\ar C\cdot \mathbf{E}\big[1\wedge \big(e^{-X}I_{K,n}\big) ; X\geq an\big].
 \eeqnn
 For any $b>a$, it is obvious that the foregoing quantities can bounded by the sum of
 \beqnn
 \varepsilon(b,n):=  C \cdot\mathbf{P}( X\in[an,bn)  )
 \quad\mbox{and}\quad
 \varepsilon(b,K,n):=  C\cdot \mathbf{E}\big[ 1\wedge \big(e^{-X}I_{K,n}\big) ; X\geq bn\big] .
 \eeqnn
 Assumption~\ref{Con.RegularV} induces that $ \varepsilon(b,n)/ \mathbf{P}(X\geq an) \to C\cdot (1-(a/b)^{\beta}) \to 0$ as $n\to\infty$ and then $b\to a+$.
 Let $S_n^b:= S_n+bn$, which drifts to $\infty$ as $n\to\infty$. 
 Let $I^b$ be the exponential functional of $S^b$ and  $I_{K,n}^b:= I_{n}^b-I_{K}^b$. 
 It is easy to see that $e^{-bn}I_{K,n} \leq I_{K,n}^b$. 
 By the independence between $X$ and $I_{K,n}$,
 \beqlb\label{eqn.2.11}
 \varepsilon(b,K,n)
 \leq C\cdot\mathbf{E}\big[ 1\wedge \big(e^{-bn}I_{K,n}\big)\big]  \cdot  \mathbf{P}(X\geq bn)
 \leq C\cdot\mathbf{E}\big[1\wedge I_{K,n}^b\big] \cdot  \mathbf{P}(X\geq bn). 
 \eeqlb 
 By Lemma~\ref{FirstLemma}, we have $I_{\infty}^b<\infty$ a.s. and hence $I_{K,n}^b\leq I_{K,\infty}^b \to 0$ a.s. as $K\to\infty$.
 Applying the dominated convergence theorem to (\ref{eqn.2.11}), we have $  \varepsilon(b,K,n) /\mathbf{P}(X\geq an) \to 0$ as $n\to\infty$ and then $K\to\infty$.
 The desired result follows by putting all preceding estimates together.
 \qed

 \begin{proposition}\label{Proposition6.03}
 	Let $H$ be a continuous, positive, bounded and non-increasing function on $(0,\infty)$. For any $K>0$, we have as $n\to\infty$,
 	\beqnn
 	\mathbf{E}\big[ H(e^{-X}I_K) \,\big| \, X\geq an\big] \to C_{H,K}>0.
 	\eeqnn
 	Moreover, the sequence $C_{H,K}  $ decreases to a limit $C_H>0$ as $K\to\infty$ and
 	\beqnn
 	\lim_{n\to\infty}\mathbf{E}\big[ H(e^{-X}I_n)\,\big| \, X\geq an\big] =C_H.
 	\eeqnn
 	
 \end{proposition}
 \proof The first result follows directly from our observation that $\mathbf{E}\big[ H(e^{-X}I_K) \,\big| \, X\geq an\big]$ is non-increasing in $n$, i.e., for any $z>0$ a sample calculation show that the expectation $\mathbf{E}\big[ H(ze^{-X}) \,\big| \, X\geq an\big]$ can written as the sum of $\mathbf{E}\big[ H(ze^{-X}) \,\big| \, X\geq a(n+1)\big]$ and
 $$
 \frac{\mathbf{P}( X\in [an,a(n+1)) )}{\mathbf{P}(X\geq an)} \Big( \mathbf{E}\big[ H(ze^{-X})\,\big|\, X\in [an,a(n+1))\big]-  \mathbf{E}\big[ H(ze^{-X}) \,\big| \, X\geq a(n+1)\big]\Big).
 $$
 Since $H$ is non-increasing, the foregoing difference is non-positive and hence $\mathbf{E}\big[ H(e^{-X}I_K) \,\big| \, X\geq an\big]\geq \mathbf{E}\big[ H(e^{-X}I_K) \,\big| \, X\geq a(n+1)\big]$.
 We now turn to prove the second claim.
 Because of the monotonicity of $H$ and $I_K$,  the sequence $C_{H,K}$ is non-increasing and converges to a limit $C_H\in [0,\infty)$ as $K\to\infty$.
 Moreover, by Proposition~\ref{Proposition6.02} and the first result, 
 \beqnn
 \lim_{n\to\infty}\mathbf{E}\big[ H(e^{-X}I_n)\,\big| \, X\geq an\big] 
 \ar=\ar \lim_{K\to \infty}\lim_{n\to\infty}\mathbf{E}\big[ \big( H(e^{-X}I_n)- H(e^{-X}I_K)\big)\,\big| \, X\geq an\big] \cr
 \ar\ar +  \lim_{K\to \infty}\lim_{n\to\infty}\mathbf{E}\big[ H(e^{-X}I_K)\,\big| \, X\geq an\big] , 
 \eeqnn 
 which equals to $C_H$. 
 It remains to prove $C_H>0$.
 Recall the modified random walk $S^b_n:= S_n+bn$ for $b>a$. Since $H$ is non-increasing, we have 
 \beqnn
 \mathbf{E}\big[ H(e^{-X}I_n)\,\big| \, X\geq an\big] \geq \frac{\mathbf{E}\big[ H(e^{-X}I_n); X\geq bn\big]}{\mathbf{P}( X\geq an )}\geq \mathbf{E}\big[ H(I_n^b) \big]\cdot \frac{\mathbf{P}( X\geq bn )}{\mathbf{P}( X\geq an )}.
 \eeqnn
 Assumption~\ref{Con.RegularV} induces that $\mathbf{P}( X\geq bn )/\mathbf{P}( X\geq an ) \to (a/n)^\beta>0$ as $n\to\infty$. 
 From Lemma~\ref{FirstLemma}, the dominated convergence theorem, the boundedness and continuity of $H$, we have $\mathbf{E}\big[ H(I_n^b) \big]\to \mathbf{E}\big[ H(I_\infty^b) \big]>0$ as $n\to\infty$ and hence $C_H=\liminf_{n\to\infty} \mathbf{E}\big[ H(e^{-X}I_n)\,\big| \, X\geq an\big]>0$.
 \qed
 
 \begin{corollary}\label{Corollary.MainThm.04.01}
 	For each $k\geq1$, we have $\mathbf{E}[F(I_n); X_k\geq an]/\mathbf{P}(X\geq an) \to \mathbf{E}[ C_{F,4}(k) ]>0$ as $n\to\infty$, where $C_{F,4}(k)$ is a random variable defined in Theorem~\ref{MainThm.04}. 
 \end{corollary}
 \proof Applying Proposition~\ref{Proposition6.02} and \ref{Proposition6.03} with $H(x):= \mathbf{E}[F(I_{k-1}+  e^{-S_{k-1}}\cdot x)\,|\, \mathscr{F}^S_{k-1}]$ to  the second conditional expectation in (\ref{eqn.6.83}),  it converges to $C_{F,4}(k)>0$ a.s.  as $n\to\infty$.  Taking this back into (\ref{eqn.6.83}) we can get the desired result immediately.
 \qed

 {\it Proof for Theorem~\ref{MainThm.04}.} We first have 
 \beqnn
 \lim_{n\to\infty}  \frac{\mathbf{E}[F(I_n)]}{\mathbf{P}(X\geq an)} 
 \ar=\ar \lim_{T\to\infty} 
 \lim_{n\to\infty}  \frac{\mathbf{E}[F(I_n);\mathcal{T}^{an}\leq T]}{\mathbf{P}(X\geq an)} \cr
 \ar\ar + \lim_{K\to\infty} \lim_{T\to\infty} 
 \lim_{n\to\infty}  \frac{\mathbf{E}[F(I_n);\mathcal{T}^{an}> T, \sigma_n^-\leq K]}{\mathbf{P}(X\geq an)}\cr
 \ar\ar + \lim_{K\to\infty} \lim_{T\to\infty} 
 \lim_{n\to\infty}  \frac{\mathbf{E}[F(I_n);\mathcal{T}^{an}> T, \sigma_n^-> K]}{\mathbf{P}(X\geq an)}.
 \eeqnn
 By (\ref{eqn.5.21}) and Proposition~\ref{Prop.MainThm.04.01}, both of the last two terms on the right side of this equality equal to $0$. Hence by (\ref{eqn.5.23}) and Lemma~\ref{Lemma.Passage.NegD},
 \beqnn
 \lim_{n\to\infty}  \frac{\mathbf{E}[F(I_n)]}{\mathbf{P}(X\geq an)} \ar=\ar \sum_{k=1}^\infty  \lim_{n\to\infty}  \frac{\mathbf{E}[F(I_n);\mathcal{T}^{an} =k]}{\mathbf{P}(X\geq an)} 
 =  \sum_{k=1}^\infty  \lim_{n\to\infty}  \frac{\mathbf{E}[F(I_n); X_k\geq an]}{\mathbf{P}(X\geq an)}.
 \eeqnn
 The asymptotic result (\ref{Eqn.MainThm.04.01}) follows directly from Corollary~\ref{Corollary.MainThm.04.01}. The finiteness of the coefficient $C_{F,4}$ follows from Proposition~\ref{Prop.MainThm.04.01}.
 \qed

 \subsection{Proof for Theorem~\ref{MainThm.05}}

 Recall the sequence $B_n$ defined in (\ref{SeqB}).
 By (\ref{MeasureChange}) we first show  in the next proposition that the impact of sample paths with fast decreasing local minimum on the expectation $\mathbf{E}^{(\Lambda )}\big[e^{-\Lambda  S_n}F(I_n)\big]$ can be asymptotically ignored.
 
 \begin{proposition}\label{Prop.5.5}
 	For every $\epsilon>0$, there exist two  integers $y_0,n_0\geq1$ such that for any $Y\geq y_0$ and $n\geq n_0$,
 	\beqnn
 	\mathbf{E}^{(\Lambda )} \big[F(I_n)e^{-\Lambda S_n};L_n\leq -Y\big] \leq \epsilon \cdot B_n.
 	\eeqnn
 \end{proposition}
 \proof Assumption~\ref{Con.F.UpperBound} induces that $F(I_n)\leq C \cdot e^{\theta_F  L_n}$. By Lemma~\ref{Lemma.LaplaceTrans} with $x,y=\infty$, there exists a constant $C>0$ such that for any $n\geq 1$, 
 \beqnn
  \mathbf{E}^{(\Lambda )}[F(I_n)e^{-\Lambda S_n};L_n\leq -Y]
 \ar\leq\ar  C\cdot e^{-(\theta_F-\Lambda)Y/2} \cdot   \mathbf{E}^{(\Lambda )} \big[e^{(\theta_F+\Lambda)/2\cdot L_n-\Lambda S_n} \big]\cr
 \ar=\ar  C\cdot e^{-(\theta_F-\Lambda)Y/2} \cdot   \mathbf{E}^{(\Lambda )} \big[e^{(\theta_F-\Lambda)/2\cdot L_n+\Lambda(L_n- S_n)} \big]
 \leq  C e^{-(\theta_F-\Lambda) Y/2} \cdot B_n
 \eeqnn
 and hence the desired result follows as $Y\to\infty$.
 \qed
 
 %
 %
 %
 %

 \begin{proposition}\label{Proposition.5.6}
 	For every $\epsilon>0$, there exist two  integers $k_0,n_0\geq1$ such that for any $K\geq k_0$ and $n\geq n_0$,
 	\beqnn
 	\mathbf{E}^{(\Lambda )}\big[e^{-\Lambda  S_n}F(I_n) ; \sigma_n^-\in[K,n-K]\big] \leq \epsilon\cdot B_n.
 	\eeqnn
 	Moreover, there exists a constant $C>0$ such that for any $n\geq 1$,
 	\beqnn
 	\mathbf{E}^{(\Lambda )}\big[e^{-\Lambda  S_n}F(I_n) \big] \leq C\cdot B_n.
 	\eeqnn
 \end{proposition}
 \proof 
 By the inequality $F(I_n) \leq C e^{\theta_F L_n }$ and the independent increments of $S$,  we have for any $k\geq 1$,
 \beqnn
 \mathbf{E}^{(\Lambda )}  \big[ e^{-\Lambda  S_n} F(I_n); \sigma_n^-=k\big]
 \ar\leq \ar C\cdot \mathbf{E}^{(\Lambda )}  \big[ e^{\theta_F S_k-\Lambda  S_n }; \sigma_n^- =k\big]\cr
 \ar\ar\cr
 \ar=\ar C\cdot  \mathbf{E}^{(\Lambda )} \big[ e^{(\theta_F-\Lambda ) S_k}; \sigma_k^- =k\big] \cdot \mathbf{E}^{(\Lambda )} \big[ e^{ -\Lambda S_{n-k}}; \tau_0^- >n-k \big].
 \eeqnn
 By the duality lemma, Remark~\ref{Remark.PassageT.01} and Lemma~\ref{Lemma.LargeDev.NegD} with $(x,y)=(0,\infty)$,  we have as $n\to\infty$,
 \beqlb\label{eqn.3.12}
 \mathbf{E}^{(\Lambda )} \big[ e^{(\theta_F-\Lambda ) S_n}; \sigma_n^- =n\big] 
 \ar=\ar \mathbf{E}^{(\Lambda )} \big[ e^{(\theta_F-\Lambda ) S_n}; M_n\leq 0\big]\cr
 \ar\ar\cr
 \ar\sim\ar \mathbf{E}^{(\Lambda )} \big[ e^{(\theta_F-\Lambda ) S_n}; \tau_0^+>n\big]\sim  \mathcal{L}^{(\theta_F-\Lambda)}_{V}(\infty)\cdot B_n. 
 \eeqlb
 From these  two results and Lemma~\ref{AuxiliaryLemma}, we have for large $n$, 
 \beqnn
 \mathbf{E}^{(\Lambda )}  \big[ e^{-\Lambda  S_n} F(I_n); \sigma_n^-\in[K,n-K] \big] 
 \ar\leq\ar C\sum_{k=K}^{n-K}\mathbf{E}^{(\Lambda )} [ e^{(\theta_F-\Lambda ) S_k}; \sigma_k^- =k] \mathbf{E}^{(\Lambda )} [ e^{ -\Lambda S_{n-k}}; \tau_0^- >n-k]\cr
 \ar\sim \ar C \sum_{k=K}^\infty \big( \mathbf{E}^{(\Lambda )} \big[ e^{(\theta_F-\Lambda ) S_k}; \tau_0^+>k\big]+ \mathbf{E}^{(\Lambda )} [ e^{ -\Lambda S_{k}}; \tau_0^- >k]\big)\cdot  B_n.
 \eeqnn
 The two desired inequalities follows   from   $\mathbf{E}^{(\Lambda )}  \big[ e^{-\Lambda  S_n} F(I_n); \sigma_n^-=0\big] \leq C\cdot  \mathbf{E}^{(\Lambda )}  \big[ e^{-\Lambda  S_n} ; \tau_0^->n\big] \leq   C\cdot B_n $ and the summability of the two  sequences $\mathbf{E}^{(\Lambda )} [ e^{(\theta_F-\Lambda ) S_k}; \tau_0^+ >k] $ and $  \mathbf{E}^{(\Lambda )} [ e^{ -\Lambda S_{k}}; \tau_0^- >k]$; see Lemma~\ref{Lemma.LargeDev.NegD} {\red and (\ref{SeqB})}.
 \qed
 
 As we mentioned before, the contribution of sample paths with large final value on the expectation $\mathbf{E}[F(I_n)]$ also can be ignored; see the following two propositions.

 \begin{proposition}\label{Propo.6.9}
 	For every $\epsilon>0$ and integer $k\geq 0$, there exist two integers $n_0,y_0\geq 1$ such that for any $n\geq  n_0$ and $Y\geq y_0$,
 	\beqnn
 	\mathbf{E}^{(\Lambda )} \big[e^{-\Lambda  S_n} F(I_n); S_n\geq Y, \sigma_n^-=k \big] \leq \epsilon\cdot B_n.
 	\eeqnn
 \end{proposition}
 \proof When $k=0$, the boundedness of $F$ and Lemma~\ref{Lemma.LargeDev.NegD} yield that there exists a constant $C>0$ independent of $Y,n$ such that
 \beqnn
 \mathbf{E}^{(\Lambda )} \big[e^{-\Lambda  S_n} F(I_n); S_n\geq Y, \sigma_n^-=0 \big]
 \ar\leq\ar C \cdot\mathbf{E}^{(\Lambda )} \big[e^{-\Lambda  S_n} ; S_n\geq Y,\tau_0^->n \big] \cr
 \ar\leq\ar  C \cdot \big|\mathcal{L}^{(\Lambda)}_{\hat{V}}(\infty)-\mathcal{L}^{(\Lambda)}_{\hat{V}}(Y)\big|  \cdot B_n,
 \eeqnn
 which is $o(B_n)$ as  $Y\to\infty$.
 For $k\geq 1$, by Assumption~\ref{Con.F.UpperBound} we have
 \beqnn
 \mathbf{E}^{(\Lambda )} \big[e^{-\Lambda  S_n} F(I_n); S_n\geq Y, \sigma_n^- =k \big]
 \ar\leq\ar C \cdot \mathbf{E}^{(\Lambda )} \big[e^{\theta_F S_k-\Lambda  S_n}; S_n\geq Y, \sigma_n^- =k \big] .
 \eeqnn
 Conditioned on $\sigma_n^-=k$, we notice that $\{ S_n\geq Y  \}\subset\{ S_n-S_k \geq Y \}$ and hence
 \beqnn
 \mathbf{E}^{(\Lambda )}  \big[e^{-\Lambda  S_n} F(I_n); S_n\geq Y, \sigma_n^- =k  \big]
 \ar\leq\ar C \cdot \mathbf{E}^{(\Lambda )} \big[e^{\theta_F S_k -\Lambda  S_n}; S_n-S_k \geq Y, \sigma_n^- =k \big].
 \eeqnn
 The independent increments of $S$ induces that the expectation on  the right side of this inequality equals to
 \beqnn
 \mathbf{E}^{(\Lambda )} \big[e^{(\theta_F-\Lambda ) S_k}; \sigma_k^- =k \big] \cdot \mathbf{E}^{(\Lambda )} \big[e^{-\Lambda  S_{n-k} }; S_{n-k}\geq Y; \tau_0^- >n-k \big].
 \eeqnn
 By (\ref{eqn.3.12}), there exists a constant $C>0$ such that $\mathbf{E}^{(\Lambda )}[e^{(\theta_F-\Lambda ) S_k}; \sigma_k^- =k]\leq C $ for any $k\geq 1$ and  
 \beqnn
 \mathbf{E}^{(\Lambda )} \big[e^{-\Lambda  S_n} F(I_n); S_n\geq Y, \sigma_n^- =k \big]
 \leq  C\cdot \mathbf{E}^{(\Lambda )} \big[e^{-\Lambda  S_{n-k} }; S_{n-k}\geq Y, \tau_0^- >n-k \big],
 \eeqnn
 which can be bounded by $ C \big|\mathcal{L}^{(\Lambda)}_{\hat{V}}(\infty)-\mathcal{L}^{(\Lambda)}_{\hat{V}}(Y)\big| \cdot B_n=o(B_n)$ for large $Y$; see Lemma~\ref{Lemma.LargeDev.NegD}. 
 \qed


 \begin{proposition}\label{Proposition.5.8}
 	For every $\epsilon>0$ and integer $k\geq 0$, there exist two integers $n_0,y_0\geq 1$ such that for any $n\geq  n_0$ and $Y\geq y_0$,
 	\beqnn
 	\mathbf{E}^{(\Lambda )} \big[e^{-\Lambda  S_n} F(I_n); S_n\geq Y, \sigma_n^-=n-k \big] \leq \epsilon\cdot B_n.
 	\eeqnn
 \end{proposition}
 \proof An argument similar to that in the proof for Proposition~\ref{Propo.6.9}  shows that the expectation in the above inequality can be bounded by
 \beqnn
 C \cdot \mathbf{E}^{(\Lambda )} \big[e^{( \theta_F -\Lambda )S_{n-k}};   \sigma_{n-k}^-= n-k \big]   \cdot \mathbf{E}^{(\Lambda )} \big[e^{-\Lambda  S_k}; S_k\geq Y, \tau_0^-> k \big].
 \eeqnn
 Dividing it by $B_n$ and then using (\ref{eqn.3.12}), we see that for large $n$, it can be bounded by $C \cdot \mathbf{E}^{(\Lambda )}[e^{-\Lambda  S_k}; S_k\geq Y, \tau_0^-> k] $, which vanishes as $Y\to\infty$.
 \qed
 
 Recall the drift parameter $a$ defined in (\ref{Constant.a}).  We now turn to consider the impact of sample paths with late arrival of the first large step, slowly decreasing local minimum and the final value below some fixed level, i.e. for $b\in(0,a]$ and $Y,T\geq 1$,
 \beqnn
 \mathbf{E}^{(\Lambda )}[e^{-\Lambda  S_n}F(I_n); -Y\leq L_n\leq S_n\leq Y, \mathcal{T}^{bn}\geq T].
 \eeqnn
 The next two propositions show that the late arrival of the first large jump will cause the random walk to drift to a low level.
 It extends Lemma~9 in \cite{BansayeVatutin2014}, which considered the case of $\beta>2$.
 \begin{proposition}\label{Prop.5.9}
 	Let $r > \frac{\beta+1}{\beta-1}$, $b\in(0,a/r) $ and $Y\geq 0$.
 	For any $\epsilon>0$, there exists an integer $n_0\geq 1$ such that for any  $n\geq n_0$,
 	\beqnn
 	\mathbf{P}^{(\Lambda )}(S_n\geq -Y, \mathcal{T}^{bn}\geq n) \leq  \epsilon \cdot B_n.
 	\eeqnn
 	
 \end{proposition}
 \proof Here we just prove this proposition with $\beta\in(1,2)$. For the case of $\beta>2$, it can be proved similarly with the help of Theorem 4.1.2(i) in \cite[p.183]{BorovkovBorovkov2008}.
 We first assume that for some $C>0$, 
 \beqlb\label{eqn.682}
 \mathbf{P}^{(\Lambda )}(X\leq -x) \leq C\cdot  \mathbf{P}^{(\Lambda )}(X\geq x) ,\quad x\geq 0.
 \eeqlb
 By Lemma~3.1 in \cite{BorovkovBoxma2003} and Assumption~\ref{Con.RegularV}, there exists a constant $C>0$ such that for large $n$,
 \beqnn
 \mathbf{P}^{(\Lambda )}( S_n+ an>rbn, \mathcal{T}^{bn-a}>n) \leq C \cdot \big|n \mathbf{P}^{(\Lambda )}(X\geq bn-a)\big|^{r}
 \leq C \cdot n^{r(1-\beta)}(\ell_3(n))^r .
 \eeqnn
 From this and the assumption that $rb<a$ and $r(\beta-1)> \beta+1$, we have for large  $n$,
 \beqnn
 \mathbf{P}^{(\Lambda )}( S_n\geq -Y, \mathcal{T}^{bn}>n)\leq C n^{r(1-\beta)}(\ell_3(n))^r =o(B_n) . 
 \eeqnn
 We now consider the general case with $X$ satisfying Assumption~\ref{Con.RegularV}.
 For some $R<0$, we define $X^R:=X\vee R$ and $X^{R}_i= X_i\vee R$ such that $a_R:=-\mathbf{E}^{(\Lambda )}[X^{R}]\in(rb, a)$.  
 One can identify that (\ref{eqn.682}) holds for $X^R$ with $x>-R$, since $\mathbf{P}^{(\Lambda )}(X^R\leq -x) =0$ and $\mathbf{P}^{(\Lambda )}(X^R\geq x)=\mathbf{P}^{(\Lambda )}(X\geq x) >0 $. 
 Moreover,	for $x\in[0,-R]$ we also have  
 	\beqnn
 	\mathbf{P}^{(\Lambda )}(X^R\leq -x)
 	\ar=\ar \frac{ \mathbf{P}^{(\Lambda )}(X^R\leq -x)}{ \mathbf{P}^{(\Lambda )}(X\geq -R)} \cdot\mathbf{P}^{(\Lambda )}(X\geq -R)  \cr
 	\ar\leq\ar  \frac{ \mathbf{P}^{(\Lambda )}(X^R\leq 0)}{ \mathbf{P}^{(\Lambda )}(X\geq -R)} \cdot\mathbf{P}^{(\Lambda )}(X\geq -R)
 	\leq \frac{ \mathbf{P}^{(\Lambda )}(X\leq 0)}{ \mathbf{P}^{(\Lambda )}(X\geq -R)} \cdot\mathbf{P}^{(\Lambda )}(X\geq x).
 	\eeqnn
 	Thus the random variable $X^R$ satisfies the inequality (\ref{eqn.682}) with $C= 1\vee \frac{ \mathbf{P}^{(\Lambda )}(X\leq 0)}{ \mathbf{P}^{(\Lambda )}(X\geq -R)}$. 
 Let $S^R_n$ be the random walk generated by the sequence $\{ X_i^R: i=1,2,\cdots\}$.   
 Noting that $S \leq S^R$ a.s., we have  $\mathbf{P}^{(\Lambda )}(S_n\geq -Y, \mathcal{T}^{bn}\geq n)
 \leq 	\mathbf{P}^{(\Lambda )}(S_n^R\geq -Y, \mathcal{T}^{bn}\geq n) $. 
The previous result induces that $	\mathbf{P}^{(\Lambda )}(S_n\geq -Y, \mathcal{T}^{bn}\geq n) = o(B_n)$ as $n\to\infty$.
 \qed
 
 \begin{proposition}\label{Proposition.6.12}
 	Let  $b>0 $ and $Y\geq 0$.
 	For any $\epsilon>0$, there exist two integers $t_0,n_0>1$ such that for any $T\geq t_0$ and $n\geq n_0$,
 	\beqnn
 	\mathbf{P}^{(\Lambda )}(-Y\leq L_n\leq S_n\leq Y, T\leq \mathcal{T}^{bn}\leq n) \leq  \epsilon \cdot B_n.
 	\eeqnn
 \end{proposition}
 \proof By the definition of $\mathcal{T}^{bn}$, we have $\{T\leq \mathcal{T}^{bn}\leq n\}\subset \cup_{k=T}^n \{ X_k\geq bn \}$ and hence 
 \beqlb\label{eqn.2.12}
 \mathbf{P}^{(\Lambda )}(-Y\leq L_n\leq S_n\leq Y, T\leq \mathcal{T}^{bn}\leq n) 
 \leq \sum_{k=T}^n 	\mathbf{P}^{(\Lambda )}(-Y\leq L_n\leq S_n\leq Y,  X_k\geq bn ).
 \eeqlb
 For $k\leq n$, by the fact that $L_{k-1}\geq L_n$ and the Markov property of $S$ we have
 \beqlb\label{eqn.6.841}
 \lefteqn{ \mathbf{P}^{(\Lambda )}(-Y\leq L_n\leq S_n\leq Y, X_k\geq bn)}\ar\ar\cr
 \ar\ar\cr
 \ar\leq\ar \mathbf{P}^{(\Lambda )}(L_{k-1}\geq -Y, S_{k-1}\geq -Y, S_{k-1}+X_k+(S_n-S_k)\in[-Y,Y] , X_k\geq bn) \cr
 \ar\ar\cr
 \ar\leq\ar \int_{-Y}^\infty \mathbf{P}^{(\Lambda )}( L_{k-1}\geq -Y, S_{k-1}\in dx) \cdot \int_{bn}^\infty \mathbf{P}^{(\Lambda )}( x+y+ \tilde{S}_{n-k} \in[-Y,Y])\mathbf{P}^{(\Lambda )}( X_k\in dy)\cr
 \ar\leq\ar \mathbf{P}^{(\Lambda )}( L_{k-1}\geq -Y )
 \cdot \sup_{x\geq -Y}\int_{bn}^\infty \mathbf{P}^{(\Lambda )}( x+y+ \tilde{S}_{n-k} \in[-Y,Y])\mathbf{P}^{(\Lambda )}( X_k\in dy).
 \eeqlb
 By Assumption~\ref{Con.LocalRegularV}, there exist constants $C_0,n_0>0$ such that  $ \mathbf{P}^{(\Lambda )}( X_k\in[bn+i,bn+i+1)) \leq C_0\cdot B_n$ for any $i\geq 1$ and $n \geq n_0$. Thus 
 \beqnn
 \lefteqn{\int_{bn}^\infty \mathbf{P}^{(\Lambda )}( x+y+ \tilde{S}_{n-k} \in[-Y,Y])\mathbf{P}^{(\Lambda )}( X_k\in dy)}\ar\ar\cr
 \ar\leq\ar C_0\cdot B_n\cdot \sum_{i=0}^\infty \mathbf{P}^{(\Lambda )}( x+bn+i+ \tilde{S}_{n-k} \in[-Y-1,Y+1]) \cr
 \ar\leq\ar C_0\cdot B_n\cdot \int_{bn}^\infty \mathbf{P}^{(\Lambda )}( x+y+ \tilde{S}_{n-k} \in[-Y-2 ,Y+2 ])dy .
 \eeqnn
 Applying Fubini's theorem to the last integral, it can be bounded by
 \beqnn
 \int_{-\infty}^\infty \int_{-x-y-Y-2 }^{-x-y+Y+2 } \mathbf{P}^{(\Lambda )}(\tilde{S}_{n-k}\in dz)dy
 =\int_{-\infty}^\infty \mathbf{P}^{(\Lambda )}(\tilde{S}_{n-k}\in dy)\int_{-x-y-Y-2}^{-x-y+Y+2}dz \leq 2Y+4.
 \eeqnn
 Taking these back into (\ref{eqn.6.841}), we have
 \beqnn
 \mathbf{P}^{(\Lambda )}(-Y\leq L_n\leq S_n\leq Y, X_k\geq bn) \ar\leq\ar 2C_0\cdot (Y+2) \cdot \mathbf{P}^{(\Lambda )}( L_{k-1}\geq -Y ) \cdot B_n .
 \eeqnn
 By Lemma~\ref{Lemma.Passage.NegD}, there exists a constant $C>0$ such that $\mathbf{P}^{(\Lambda )}( L_{k-1}\geq -Y )=\mathbf{P}^{(\Lambda )}_Y( \tau_0^->k-1 ) \leq C\cdot \mathbf{P}^{(\Lambda )}( X\geq a(k-1) )$ for any $k>1$. Taking these two upper bound estimates back into (\ref{eqn.2.12}), we have for some constant $C>0$ independent of $n$ and $T$,
 \beqnn
 \frac{1}{B_n} \cdot \mathbf{P}^{(\Lambda )}(-Y\leq L_n\leq S_n\leq Y, T\leq \mathcal{T}^{bn}\leq n)
 \ar\leq\ar   C  \sum_{k=T}^\infty \mathbf{P}^{(\Lambda )}(X\geq a(k-1)),
 \eeqnn
 which vanishes as $T\to\infty$, since $ \sum_{k=1}^\infty \mathbf{P}^{(\Lambda )}(X\geq a(k-1))<\infty$.  The desired result follows.
 \qed

 Putting all preceding estimates together, we see that the main contribution to the expectation $\mathbf{E}^{(\Lambda )}[F(I_n)e^{-\Lambda S_n}]$ is made by the sample paths with an early large step. 
 In detail, for $T,n\in\mathbb{Z}$ with $0<T<n$ we have 
 \beqlb\label{eqn.5.101}
  \mathbf{E}^{(\Lambda )} \big[F(I_n)e^{-\Lambda S_n} \big]
  \ar=\ar \mathbf{E}^{(\Lambda )} \big[F(I_n)e^{-\Lambda S_n}; \mathcal{T}^{bn}<T \big] +  \mathbf{E}^{(\Lambda )} \big[F(I_n)e^{-\Lambda S_n} ;  \mathcal{T}^{bn}>n\big]\cr
  \ar\ar +  \mathbf{E}^{(\Lambda )} \big[F(I_n)e^{-\Lambda S_n} ; T\leq  \mathcal{T}^{bn} \leq n\big] ,
 \eeqlb
 where  $b\in(0,a/r) $ and $r>\frac{\beta+1}{\beta-1}$. 
 For any $Y>0$, the second expectation on the right side of this equality is equal to
 \beqnn
 \ar\ar \mathbf{E}^{(\Lambda )} \big[F(I_n)e^{-\Lambda S_n} ; S_n\geq -Y, \mathcal{T}^{bn}>n\big] +  \mathbf{E}^{(\Lambda )} \big[F(I_n)e^{-\Lambda S_n} ;  S_n\leq -Y, \mathcal{T}^{bn}>n\big]\cr 
 \ar\ar \leq C \cdot e^{\Lambda Y} \cdot \mathbf{E}^{(\Lambda )} \big[S_n\geq -Y, \mathcal{T}^{bn}>n\big] +  \mathbf{E}^{(\Lambda )} \big[F(I_n)e^{-\Lambda S_n} ;  L_n\leq -Y \big],
 \eeqnn
 which is $o(B_n)$ as $n\to\infty$; see Proposition~\ref{Prop.5.5} and \ref{Prop.5.9}. 
 The third expectation on the right side of (\ref{eqn.5.101}) is equal to
 \beqnn
\ar\ar \mathbf{E}^{(\Lambda )} \big[F(I_n)e^{-\Lambda S_n} ; -Y\leq L_n\leq S_n\leq Y, T\leq  \mathcal{T}^{bn} \leq n\big]\cr
 \ar\ar + \mathbf{E}^{(\Lambda )} \big[F(I_n)e^{-\Lambda S_n} ; L_n<-Y,  T\leq  \mathcal{T}^{bn} \leq n\big] + \mathbf{E}^{(\Lambda )} \big[F(I_n)e^{-\Lambda S_n} ; S_n> Y, T\leq  \mathcal{T}^{bn} \leq n\big].
 \eeqnn
 By the boundedness of $F$ and Proposition~\ref{Proposition.6.12}, the first term can be bounded by $C\cdot e^{\Lambda Y}\cdot \mathbf{E}^{(\Lambda )} \big[-Y\leq L_n\leq S_n\leq Y, T\leq  \mathcal{T}^{bn} \leq n\big]=o(B_n)$ as $n\to\infty$. 
 The second term is smaller than $\mathbf{E}^{(\Lambda )} \big[F(I_n)e^{-\Lambda S_n} ; L_n<-Y\big] $, which is $o(B_n) $ as $n,Y\to\infty$; see Proposition~\ref{Prop.5.5}. 
 For any $0<K< n$, the third term can be bounded by
 \beqnn
 \mathbf{E}^{(\Lambda )} \big[F(I_n)e^{-\Lambda S_n} ; S_n> Y\big] 
\ar = \ar \mathbf{E}^{(\Lambda )} \big[F(I_n)e^{-\Lambda S_n} ; S_n> Y; \sigma_n^-\in[K,n-K]\big]\cr
\ar\ar +  \sum_{k=0}^{K-1}\mathbf{E}^{(\Lambda )} \big[F(I_n)e^{-\Lambda S_n} ; S_n> Y; \sigma_n^-=k\big] \cr
\ar\ar  + \sum_{k=0}^{K-1} \mathbf{E}^{(\Lambda )} \big[F(I_n)e^{-\Lambda S_n} ; S_n> Y; \sigma_n^-=n-k\big] , 
 \eeqnn 
 which is $o(B_n)$ as $n,Y\to\infty$ and then $K\to \infty$; see Proposition~\ref{Proposition.5.6}, \ref{Propo.6.9} and \ref{Proposition.5.8}. 
Taking these estimates back into (\ref{eqn.5.101}), we have as $n\to\infty$,
 \beqnn
 \mathbf{E}^{(\Lambda )} \big[F(I_n)e^{-\Lambda S_n} \big]  \sim \sum_{k=1}^\infty  \mathbf{E}^{(\Lambda )} \big[F(I_n)e^{-\Lambda S_n}; \mathcal{T}^{bn}=k \big]. 
 \eeqnn
 By (\ref{Equvi.LargeJ}) and Lemma~\ref{Lemma.Passage.NegD}, for each $k\geq 1$ we have
 \beqlb\label{eqn5.38.1}
 \mathbf{E}^{(\Lambda )} \big[F(I_n)e^{-\Lambda S_n}; \mathcal{T}^{bn}=k \big] \sim
 \mathbf{E}^{(\Lambda )} \big[F(I_n)e^{-\Lambda S_n}; X_k\geq bn \big],
 \eeqlb
 as $n\to\infty$.  Using Proposition~\ref{Prop.5.5} again,  it is also asymptotically equivalent to
 \beqlb\label{eqn5.39}
 \mathbf{E}^{(\Lambda )} \big[F(I_n)e^{-\Lambda S_n}; S_{k-1}\geq -Y, X_k\geq bn \big]
 \eeqlb
 for large $Y$. Moreover, as we have mentioned before, the local minimum is usually attained around the early large step. 
 For sample paths that stay above a high level before the early large step, their impacts on the expectation $\mathbf{E}^{(\Lambda )}[F(I_n)e^{-\Lambda S_n}]$  can be ignored; see the following proposition.

 \begin{proposition}\label{Proposition6.13}
 	Let $b>0$ and  $k\in \mathbb{Z}_+$.	For each $\epsilon>0$, there exist two integers $n_0,y_0\geq 1$ such that for any $Y\geq y_0$ and $n\geq  n_0$,
 	\beqnn
 	\mathbf{E}^{(\Lambda )} \big[F(I_n)e^{-\Lambda S_n};S_{k-1}\geq Y, X_k\geq bn \big]
 	\leq \epsilon \cdot B_n.
 	\eeqnn
 \end{proposition}
 \proof  
 By Assumption~\ref{Con.F.UpperBound},  we have $F(I_n)\leq C\cdot e^{\theta_F L_n}$  for some constant $C>0$ depending only on $F$ and hence the expectation in the desired inequality can be bounded by $ C\cdot\mathbf{E}^{(\Lambda )}\big[e^{\theta_F L_n-\Lambda S_n};S_{k-1}\geq Y , X_k\geq bn \big] $.
 For a constant $\gamma\in(1/\beta,1)$, we write this expectation into
 \beqlb\label{eqn5.40}
 \mathbf{E}^{(\Lambda )}\big[e^{\theta_F L_n-\Lambda S_n};S_{k-1}\geq n^{\gamma} , X_k\geq bn \big] + \mathbf{E}^{(\Lambda )}\big[e^{\theta_F L_n-\Lambda S_n};Y\leq S_{k-1}\leq n^{\gamma} , X_k\geq bn \big] .
 \eeqlb 
 By H\"older's inequality  and then  the independence between $S_{k-1}$ and $X_k$,  the first expectation can be bounded by
 \beqnn 
 \big(\mathbf{E}^{(\Lambda )}[e^{2\theta_F L_n-2\Lambda S_n }]\big)^{1/2} \cdot \big(\mathbf{P}^{(\Lambda )}(S_{k-1}\geq n^{\gamma}) \cdot \mathbf{P}^{(\Lambda )}( X_k\geq bn)\big)^{1/2} .
 \eeqnn
 Using Lemma~\ref{Lemma.LaplaceTrans} with $x,y=\infty$ and (\ref{SeqB}), we see that $\mathbf{E}^{(\Lambda )}[e^{2\theta_F L_n-2\Lambda S_n }] \leq C \cdot B_n$ and $\mathbf{P}^{(\Lambda )}( X_k\geq bn) \leq C\cdot n\cdot B_n$ for any $n\geq 1$ and some constant $C>0$. Hence 
 \beqnn
 \mathbf{E}^{(\Lambda )}\big[e^{\theta_F L_n-\Lambda S_n};S_{k-1}\geq n^{\gamma} , X_k\geq bn \big] 
 \leq C \cdot B_n \cdot \big(n\cdot\mathbf{P}^{(\Lambda )}(S_{k-1}\geq n^{\gamma})\big)^{1/2}.
 \eeqnn 
 The convolution property of regular varying distributions together with Assumption~\ref{Con.RegularV} induces that  $\mathbf{P}^{(\Lambda )}(S_{k-1}\geq n^{\gamma})\sim (k-1)\mathbf{P}^{(\Lambda )}(X\geq n^{\gamma})=o(1/n)$.
 Taking this back into the preceding inequality, we have as $n\to\infty$,
 \beqlb\label{eqn5.400}
 \mathbf{E}^{(\Lambda )}\big[e^{\theta_F L_n-\Lambda S_n};S_{k-1}\geq n^{\gamma} , X_k\geq bn \big]  = o(B_n). 
 \eeqlb
 On the other hand, by the Markov property of $S$ and  
 $L_n\leq \min\{S_k,S_{k+1},\cdots,S_n\}$, 
 the second expectation in (\ref{eqn5.40}) can be bounded by
 \beqlb\label{eqn6.107}
 \lefteqn{ \int_Y^{n^{\gamma}}   \mathbf{P}^{(\Lambda )} (S_{k-1}\in dy )  \int_{bn}^\infty  \mathbf{E}^{(\Lambda )}_{y+x}\big[e^{\theta_F\tilde{L}_{n-k}-\Lambda  \tilde{S}_{n-k}}\big]\mathbf{P}^{(\Lambda )}(X_k\in dx) }\ar\ar\cr
 \ar=\ar \int_Y^{n^{\gamma}}   \mathbf{P}^{(\Lambda )} (S_{k-1}\in dy )  \int_{bn +y}^\infty  \mathbf{E}^{(\Lambda )}_{x}\big[e^{\theta_F\tilde{L}_{n-k}-\Lambda  \tilde{S}_{n-k}}\big]\mathbf{P}^{(\Lambda )}(X_k+y\in dx) .
 \eeqlb
 We first consider the inner integral  on the right-hand side of  this equality. 
 It is easy to see that  
 \beqnn
 \lefteqn{  \int_{bn +y}^\infty  \mathbf{E}^{(\Lambda )}_{x}\big[e^{\theta_F\tilde{L}_{n-k}-\Lambda  \tilde{S}_{n-k}}\big]\mathbf{P}^{(\Lambda )}(X_k+y\in dx) }\ar\ar\cr
 \ar\leq\ar  \sum_{i=0}^\infty   \mathbf{E}^{(\Lambda )}_{bn+i}[e^{\theta_F(\tilde{L}_{n-k}+1)-\Lambda  (\tilde{S}_{n-k}-1)}]\cdot \mathbf{P}^{(\Lambda )}(X_k+y\in [bn+i,bn+i+1)).
 \eeqnn
 Using Assumption~\ref{Con.RegularV}, Potter's theorem; see \cite[Theorem~1.5.6(iii)]{BinghamGoldieTeugels1987} and then  Assumption~\ref{Con.LocalRegularV}, there exist constants $C,C_1,C_2,n_0>0$ such that for any $i\geq 0$, $n\ge n_0$  and $y\in [Y,n^\gamma]$,
 \beqnn
 \mathbf{P}^{(\Lambda )}(X_k+y\in [bn+i,bn+i+1)) 
 \ar\leq\ar C_1\cdot \beta (bn-y+i)^{-1-\beta} \ell_3(bn-y+i) \cr 
 \ar\ar\cr
 \ar\leq\ar C_2\cdot \beta (bn+i)^{-1-\beta} \ell_3(bn+i) \cr 
 \ar\ar\cr
 \ar \leq\ar C \cdot  \mathbf{P}^{(\Lambda )}(X\in [bn+i,bn+i+1))
 \eeqnn  and hence 
 \beqnn
 \lefteqn{  \int_{bn +y}^\infty  \mathbf{E}^{(\Lambda )}_{x}\big[e^{\theta_F\tilde{L}_{n-k}-\Lambda  \tilde{S}_{n-k}}\big]\mathbf{P}^{(\Lambda )}(X_k+y\in dx)  }\ar\ar\cr
 \ar\leq \ar C \sum_{i=0}^\infty   \mathbf{E}^{(\Lambda )}_{bn+i}[e^{\theta_F\tilde{L}_{n-k}-\Lambda  \tilde{S}_{n-k}}]\cdot \mathbf{P}^{(\Lambda )}(X\in [bn+i,bn+i+1)) \cr
 \ar\leq \ar C \int_{bn}^\infty \mathbf{E}^{(\Lambda )}_x[e^{\theta_F\tilde{L}_{n-k}-\Lambda  \tilde{S}_{n-k}}]\cdot \mathbf{P}^{(\Lambda )}(X\in dx)
 = C\cdot \mathbf{E}^{(\Lambda )} \big[  e^{\theta_F\tilde{L}_{n-k+1}-\Lambda  \tilde{S}_{n-k+1}}; \tilde{X}_1\geq bn \big].
 \eeqnn
 By Lemma~\ref{Lemma.LaplaceTrans} with $x,y=\infty$, it can be bounded by $ C\cdot \mathbf{E}^{(\Lambda )} \big[  e^{\theta_F\tilde{L}_{n-k+1}-\Lambda  \tilde{S}_{n-k+1}}  \big] \leq C\cdot B_n$.  
 Taking this  back into (\ref{eqn6.107}), there exists a constant $C$ such that  for large $n$ and $Y$,
 \beqnn
 \mathbf{E}^{(\Lambda )} \big[e^{\theta_F L_n-\Lambda S_n};Y\leq S_{k-1}\leq n^{\gamma}, X_k\geq bn\big]
 \leq C \cdot B_n\cdot \mathbf{P}^{(\Lambda )} (S_{k-1}\geq Y) = o(B_n).
 \eeqnn  
 Taking this and (\ref{eqn5.400}) back into (\ref{eqn5.40}), we can get the desired result immediately. 
 \qed
 
 Taking the estimate in Proposition~\ref{Proposition6.13} back into (\ref{eqn5.39}) and then (\ref{eqn5.38.1}),   we have 
 \beqnn
 \mathbf{E}^{(\Lambda )} \big[F(I_n)e^{-\Lambda S_n}; \mathcal{T}^{bn}=k \big] \sim  \mathbf{E}^{(\Lambda )} \big[e^{-\Lambda S_{k-1}}\cdot F(I_n)\cdot e^{-\Lambda (S_n-S_{k-1})}; |S_{k-1}|\leq Y, X_k\geq bn \big]
 \eeqnn
 for large $n$ and $Y$. The next proposition shows that for large $Z>0$, the contribution of sample path with $|S_n-S_{k-1}|>Z$ to the forgoing expectations can be asymptotically ignored.


 \begin{proposition}\label{Proposition.5.12}
 	Let $Y\geq 0$ and $k\in \mathbb{Z}_+$.	For any $\epsilon>0$, there exist two integers $z_0,n_0\geq 1$ such that for any $Z\geq z_0$ and $n\geq n_0$,
 	\beqnn
 	\mathbf{E}^{(\Lambda )} \big[F(I_n)e^{-\Lambda S_n};|S_{k-1}|\leq Y,|S_n-S_{k-1}|\geq Z \big] \leq \epsilon\cdot B_n.
 	\eeqnn
 	
 \end{proposition}
 \proof By Assumption~\ref{Con.F.UpperBound}, it suffices to prove this equality with $F(I_n)$ replaced by $e^{\theta_F L_n}$.
 For $1\leq k\leq n$, notice that $L_n \leq  S_{k-1}+ L_{k,n}$ with $L_{k,n}:=\min\{ S_i-S_{k-1}; i=k-1,\cdots,n \}$.
 By the independent increments of $S$, 
 \beqnn
 \lefteqn{\mathbf{E}^{(\Lambda )} \big[e^{\theta_F L_n-\Lambda S_n};|S_{k-1}|\leq Y,|S_n-S_{k-1}|\geq Z \big]}\ar\ar\cr
 \ar\ar\cr
 \ar\leq\ar  \mathbf{E}^{(\Lambda )} \big[e^{\theta_F (S_{k-1}+ L_{k,n})-\Lambda S_n};|S_{k-1}|\leq Y,|S_n-S_{k-1}|\geq Z \big] \cr
 \ar\ar\cr
 \ar=\ar  \mathbf{E}^{(\Lambda )}\big[e^{(\theta_F -\Lambda )S_{k-1} };|S_{k-1}|\leq Y \big] \cdot \mathbf{E}^{(\Lambda )} \big[e^{  \theta_F L_{n-k+1}-\Lambda S_{n-k+1}}; |S_{n-k+1}|\geq Z \big] .
 \eeqnn
 The first expectation on the right side of the last equality can be bounded by $e^{(\theta_F -\Lambda )Y } $. 
 Moreover, notice that $\{S_{n-k+1} \leq -Z\}\subset \{L_{n-k+1}\leq -Z\}$ and $\{S_{n-k+1}\geq Z\}\subset \{ S_{n-k+1} -L_{n-k+1}\geq Z\}$,   
 the second expectation can be bounded by
 \beqlb\label{eqn.5.49}
 \mathbf{E}^{(\Lambda )} \big[e^{  \theta_F L_{n-k+1}-\Lambda S_{n-k+1}};  L_{n-k+1} \leq -Z \big]
 +  \mathbf{E}^{(\Lambda )} \big[e^{  \theta_F L_{n-k+1}-\Lambda S_{n-k+1}}; S_{n-k+1}- L_{n-k+1}\geq Z \big] .
 \eeqlb
 The first term can be written into
 $$ \mathbf{E}^{(\Lambda )} \big[e^{  \theta_F L_{n-k+1}-\Lambda S_{n-k+1}}\big]- \mathbf{E}^{(\Lambda )} \big[e^{  \theta_F L_{n-k+1}-\Lambda S_{n-k+1}}; - L_{n-k+1} < Z \big].$$
 Applying Lemma~\ref{Lemma.LaplaceTrans} to the first expectation with $\lambda_1=\theta_F-\Lambda$, $\lambda_2=\Lambda$, $(x,y)=(\infty,\infty)$ and then to the second expectation with $(x,y)=(Z,\infty)$, 
 we have for some constants $C,n_0>0$, the first expectation in (\ref{eqn.5.49}) can be bounded by
 \beqnn 
 C  \int_{Z}^\infty e^{-(\theta_F-\Lambda) z}\big(V(dz) + V(z)dz\big) 
 \cdot B_n,
 \eeqnn 
 uniformly in $n\geq n_0$. 
 By the two facts that $\theta_F-\Lambda>0$ and $V(z)=O(z)$ as $z\to\infty$, the preceding integral vanishes as $Z\to\infty$. 
 Similarly,  the second expectation in (\ref{eqn.5.49}) also can be uniformly bounded by
 \beqnn 
 C  \int_{Z}^\infty e^{-\Lambda z}\big(\hat{V}(dz) + \hat{V}(z)dz\big)
 \cdot B_n
 \eeqnn
 and this integral also vanishes as $Z\to\infty$. 
 Hence the desired result follows.
 \qed
 %
 
 Let $r>\frac{\beta+1}{\beta-1}$ and  $b\in(0,a/r) $.
 Combining all estimates above together, we have 
 \beqnn
 \mathbf{E}^{(\Lambda )} [e^{-\Lambda  S_n}F(I_n)]
 \ar\sim\ar \sum_{k=1}^\infty \mathbf{E}^{(\Lambda )} \big[e^{-\Lambda  S_n}F(I_n);|S_{k-1}|\leq Y, |S_n-S_{k-1}|\leq Z, X_k\geq bn \big]
 \eeqnn
 for large $Y$, $Z$ and $n$.
 By the Markov property of $S$, the summand can be written into
 \beqlb\label{eqn.5.12}
 \mathbf{E}^{(\Lambda )} \big[e^{-\Lambda  S_{k-1}} F_{\mathscr{F}^S}^{b,Z}(k,n) ;|S_{k-1}|\leq Y  \big]
 \eeqlb
 with
 \beqnn
 F_{\mathscr{F}^S}^{b,Z}(k,n):= \mathbf{E}^{(\Lambda )} \big[e^{-\Lambda  \tilde{S}_{n-k+1}} F\big(I_{k-1}+ e^{-S_{k-1}} \tilde{I}_{n-k+1}  \big); |\tilde{S}_{n-k+1}|\leq Z, \tilde{X}_1\geq bn \big| \mathscr{F}^{S}_{k-1} \big] .
 \eeqnn
 It is obvious that $ F_{\mathscr{F}^S}^{b,Z}(k,n)$ plays the main role in the asymptotic analysis of the expectation $\mathbf{E}^{(\Lambda )} [e^{-\Lambda  S_n}F(I_n)]$.
 Recall the constant $\kappa\in(1,2)$ defined in Condition~\ref{MomentCon}. 
 The next proposition states that $\{ |\tilde{S}_{n-k+1}|\leq Z , \tilde{X}_1\geq bn \}$ is asymptotically equivalent to $\{ |\tilde{S}_{n-k+1}|\leq Z, |\tilde{X}_1-an|\geq n^{1/\kappa} \}$.

 \begin{proposition}
 	For any $b\in(0,a)$, $Z\geq 1$ and $\epsilon>0$, there exists an integer $n_0$ such that for any $n\geq n_0$,
 	\beqnn
 	\mathbf{P}^{(\Lambda )} \big(X_1\in[b n , an- n^{1/\kappa}] \cup [an+ n^{1/\kappa},\infty); |S_n| \leq Z \big)  \leq \epsilon \cdot B_n.
 	\eeqnn
 \end{proposition}
 \proof 
 For $\delta\in(0,1)$ and large $n\geq 1$, we first have 
 \beqlb\label{eqn.1005}
  \lefteqn{\mathbf{P}^{(\Lambda )}\big(X_1\in[b n , an-n^{1/\kappa}]; |S_n| \leq Z\big)
 =\int_{b n}^{an-n^{1/\kappa}} \mathbf{P}^{(\Lambda )}(|x+S_{n-1}|\leq Z) \mathbf{P}^{(\Lambda )}(X\in dx) }\qquad\qquad\qquad\qquad\  \ar\ar \cr
 \ar\leq \ar \sum_{k=0}^{[((a-b)n-n^{1/\kappa})/\delta]}\int_{b n+ k\delta}^{bn+(k+1)\delta} \mathbf{P}^{(\Lambda )}(|x+S_{n-1}|\leq Z) \mathbf{P}^{(\Lambda )}(X\in dx) .
 \eeqlb
Notice that $
 \{  |x+S_{n-1}|\leq Z \} \subset\{ |b n+ k\delta+S_{n-1}|\leq Z+1 \}$ uniformly in $x
 \in [b n+ k\delta,b n+ (k+1)\delta)$. Hence the preceding summand can be bounded by
 \beqnn
  \mathbf{P}^{(\Lambda )}\big(|b n+ k\delta+S_{n-1}|\leq Z+1\big)\cdot \mathbf{P}^{(\Lambda )}\big(X\in (b n+ k\delta,b n+ (k+1)\delta]\big). 
 \eeqnn
 By Assumption~\ref{Con.LocalRegularV} and (\ref{SeqB}), there exist constants $C>0$ and $n_0>0$ such that for any $n\geq n_0$ and $0\leq k\leq [((a-b)n-n^{1/\kappa})/\delta]$, 
 \beqnn
 \mathbf{P}^{(\Lambda )}\big(X\in (b n+ k\delta,b n+ (k+1)\delta]\big) \leq C \cdot  \frac{\beta}{bn+ k\delta }\cdot \mathbf{P}^{(\Lambda )}(X\geq b n)\cdot \delta \leq C\cdot B_n\cdot \delta.
 \eeqnn
 Taking these two estimates back into  (\ref{eqn.1005}), 
 we have 
  \beqnn
 \mathbf{P}^{(\Lambda )} \big(X_1\in[b n , an-n^{1/\kappa}]; |S_n| \leq Z\big) 
 \leq  C\cdot B_n  \sum_{k=0}^{[((a-b)n-n^{1/\kappa})/\delta]}\mathbf{P}^{(\Lambda )}\big(|b n+ k\delta+S_{n-1}|\leq Z+1\big) \cdot  \delta . 
 \eeqnn
 Notice that $\{ |b n+ k\delta+S_{n-1}|\leq Z+1 \} \subset \{ |x+S_{n-1}|\leq Z+2 \}$ uniformly in $x \in (bn+k\delta, bn+(k+1)\delta]$ and $k\geq 0$. Then for large $n\geq 1$, 
 \beqnn
 \mathbf{P}^{(\Lambda )}\big(X_1\in[b n , an-n^{1/\kappa}]; |S_n| \leq Z\big) 
 \ar \leq \ar C\cdot B_n  \cdot \int_{b n}^{an-n^{1/\kappa}+1} \mathbf{P}^{(\Lambda )}(   |x+S_{n-1}|\leq Z+2) dx \cr
 \ar = \ar  C \cdot B_n \cdot \int_{b n}^{an-n^{1/\kappa}+1} \int_{-Z-2-x}^{Z+2-x} \mathbf{P}^{(\Lambda )}(    S_{n-1}\in dz) dx \cr
 \ar\leq\ar C(Z+2)  \cdot B_n \cdot    \int_{-Z-3-an+n^{1/\kappa}}^{Z+2-b n} \mathbf{P}^{(\Lambda )}(    S_{n-1}\in dz) \cr
 \ar\ar\cr
 \ar\leq\ar C(Z+2)  \cdot B_n \cdot   \mathbf{P}^{(\Lambda )}(  S_{n-1}\geq -Z-3-an+n^{1/\kappa})\cr  
 \ar\ar\cr
 \ar\leq \ar C (Z+2)\cdot B_n \cdot  \mathbf{P}^{(\Lambda )} \big(    n^{-1/\kappa}\cdot \big( S_{n-1}+a(n-1) \big)\geq 1/2 \big).
 \eeqnn 
 Here  the constant $C>0$ is independent of $n,Z$ and may change from line to line. 
 The statement below Condition~\ref{MomentCon} shows that $n^{-1/\kappa}\cdot \big( S_{n-1}+a(n-1) \big)\to 0$ in probability and hence $\mathbf{P}^{(\Lambda )}(X_1\in[b n , an-n^{1/\kappa}]; |S_n| \leq Z) \to 0$ as $n\to\infty$. Similarly, we also have for large $n$,
 \beqnn
 \mathbf{P}^{(\Lambda )} \big(X_1 > an+n^{1/\kappa} ; |S_n| \leq Z\big)
 \ar\leq\ar C \cdot B_n\cdot  \int_{an+ n^{1/\kappa}}^\infty  \mathbf{P}^{(\Lambda )}( |x+S_{n-1}| \leq  Z+2) dx\cr
 \ar\ar\cr
 \ar\leq\ar C \cdot B_n\cdot  \mathbf{P}^{(\Lambda )}(   S_{n-1}\leq Z+2-an- n^{1/\kappa})\cr
 \ar\ar\cr
 \ar\leq\ar C \cdot B_n\cdot   \mathbf{P}^{(\Lambda )} \big(   n^{-1/\kappa}\cdot \big( S_{n-1}+a(n-1) \big)\leq -1/2 \big),
 \eeqnn
 which also goes to $0$ as $n\to\infty$.
 \qed
 
 By this proposition, for large $n$ we see that  $ F_{\mathscr{F}^S}^{b,Z}(k,n)$ can be well approximated by
 \beqnn
 \mathbf{E}^{(\Lambda )} \big[e^{-\Lambda \tilde{S}_{n-k+1}} F\big(I_{k-1}+ e^{-S_{k-1}} \cdot \tilde{I}_{n-k+1}  \big);  |\tilde{S}_{n-k+1}|\leq Z, |\tilde{X}_1- an|\leq n^{1/\kappa} \big| \mathscr{F}^{S}_{k-1} \big]. 
 \eeqnn
 Let $\tilde{I}_{k,n}:=\tilde{I}_n-\tilde{I}_k$ for $0\leq k\leq n$. 
 In the next proposition, we prove that this conditional expectation will not change too much with $\tilde{I}_{n-k+1}$ replaced by $ \tilde{I}_{n-k+1-W,n-k+1}$ for large $W$.
 
 
 \begin{proposition}\label{Proposition6.16}
 	Let $H$ be a bounded and Lipschitz continuous function on $(0,\infty)$. For every $ Z,\epsilon>0$, there exist two integers $w_0, n_0 \geq 1$ such that for any $W\geq w_0$ and $n\geq n_0$,
 	\beqnn
 	\mathbf{E}^{(\Lambda )} \big[e^{-\Lambda  S_{n}} \big| H(I_n)- H(I_{n-W,n})\big|;  |S_{n}|\leq Z, |X_1- an|\leq n^{1/\kappa}  \big] \leq \epsilon \cdot B_n.
 	\eeqnn

 \end{proposition}
 \proof The Lipschitz continuity of $H$ induces that $ \big| H(I_n)- H(I_{n-W,n})\big| \leq C\cdot \big(I_{n-W} \wedge 1 \big)$ for some constant $C>0$ independent of $S,n,W$. 
 Hence the expectation in the desired inequality can be bounded by 
 \beqnn
 \lefteqn{Ce^{\Lambda Z}	\cdot \mathbf{E}^{(\Lambda )} \big[\big(I_{n-W} \wedge 1 \big) ;  |S_{n}|\leq Z, |X_1- an|\leq n^{1/\kappa}  \big]}\ar\ar\cr
 \ar\ar\cr
 \ar=\ar Ce^{\Lambda Z}\cdot	\mathbf{E}^{(\Lambda )} \big[ \big(\big(e^{-S_n}\hat{I}_{W-1,n-1} \big)\wedge 1 \big) ;  |S_{n}|\leq Z, |X_n- an|\leq n^{1/\kappa}  \big]\cr
 \ar\ar\cr
 \ar\leq\ar Ce^{\Lambda Z}\cdot	\mathbf{E}^{(\Lambda )} \big[ \big(\big(e^{Z}\hat{I}_{W-1,n-1} \big)\wedge 1 \big) ;  |S_{n-1}+X_n|\leq Z, |X_n- an|\leq n^{1/\kappa}  \big].
 \eeqnn
 Here the equality follows from  the duality lemma.  
 By the independence between $X_n$ and $\{S_k:k=0,1,\cdots,n-1\}$, the preceding quantities can be written as 
 \beqnn
 Ce^{\Lambda Z}\cdot \int_{an-n^{1/\kappa}}^{an+n^{1/\kappa}} 	\mathbf{E}^{(\Lambda )} \big[ \big((e^{Z}\hat{I}_{W-1,n-1} )\wedge 1 \big) ;  |S_{n-1}+x_n|\leq Z \big] \mathbf{P}^{(\Lambda )}(X_n\in dx_n). 
 \eeqnn
 An argument similar to that below (\ref{eqn.6.841}) induces that the foregoing quantity can be bounded by
 \beqnn
 Ce^{\Lambda Z}\cdot B_n\cdot \int_{an-n^{1/\kappa}}^{an+n^{1/\kappa}} 	\mathbf{E}^{(\Lambda )} \big[ \big((e^{Z}\hat{I}_{W-1,n-1} )\wedge 1 \big) ;  |S_{n-1}+x_n|\leq Z+2 \big] dx_n
 \eeqnn
 for large $n$ and some constant $C>0$ independent of $n$.  
 By the change of variables, it turns to be
 \beqnn
 Ce^{\Lambda Z}\cdot B_n\cdot \int_{-Z-2}^{Z+2} 	\mathbf{E}^{(\Lambda )} \big[ \big((e^{Z}\hat{I}_{W-1,n-1} )\wedge 1 \big) ;  |S_{n-1}+an-z|\leq n^{1/\kappa} \big] dz,
 \eeqnn
 which can be bounded by $Ce^{\Lambda Z}\cdot B_n\cdot \mathbf{E}^{(\Lambda )} [ ((e^{Z}\hat{I}_{W-1,n-1} )\wedge 1 )] \cdot 2(Z+1)$.
 By the dominated convergence theorem, we have $\mathbf{E}^{(\Lambda )} [ ((e^{Z}\hat{I}_{W-1,n-1} )\wedge 1 )]\to \mathbf{E}^{(\Lambda )} [ ((e^{Z}\hat{I}_{W-1,\infty} )\wedge 1 )]$ as $n\to\infty$. 
 Putting all estimates above together, there exists a constant $C>0$ independent of $W$ such that for large $n>1$,
 $$
 \mathbf{E}^{(\Lambda )} \big[e^{-\Lambda  S_{n}} \big| H(I_n)- H(I_{n-W,n})\big|;  |S_{n}|\leq Z, |X_1- an|\leq n^{1/\kappa}  \big] \leq C\cdot \mathbf{E}^{(\Lambda )} [ ((e^{Z}\hat{I}_{W-1,\infty} )\wedge 1 )] \cdot B_n.
 $$
 Since  $\hat{I}_\infty<\infty$ a.s., we have $\hat{I}_{W-1,\infty} \to 0$ a.s. as $W\to\infty$ and the desired result follows. 
 \qed

 \begin{proposition}\label{Proposition6.17}
 	Let $H$ be a nonnegative, bounded and continuous function on $(0,\infty)$. For any two integers $W,Z> 0$, we have as $n\to\infty$,
 	\beqnn
 	\frac{1}{B_n}\cdot \mathbf{E}^{(\Lambda )} \Big[e^{-\Lambda  S_{n}}   H(I_{n-W,n}) ;  |S_{n}|\leq Z, |X_1- an|\leq n^{1/\kappa}  \Big]
 	\ar \to\ar  \int_{-Z}^Z \mathbf{E}^{(\Lambda )} \big[e^{-\Lambda  z}  H(e^{-z}(1+\hat{I}_{W-1})) \big] dz.
 	\eeqnn
 	Moreover,  the limit coefficient converges as $W\to\infty$ and then $Z\to \infty$  to
 	\beqnn
 	\int_{-\infty}^\infty\mathbf{E}^{(\Lambda )} \big[e^{-\Lambda  z}  H(e^{-z}(1+\hat{I}_\infty)) \big] dz.
 	\eeqnn
 	
 \end{proposition}
 \proof By the duality lemma and the independent increments of $S$, we first have for $n>W$, 
 \beqlb\label{eqn.3.20}
 \lefteqn{ \mathbf{E}^{(\Lambda )} \big[e^{-\Lambda  S_{n}}   H(I_{n-W,n}) ;  |S_{n}|\leq Z, |X_1- an|\leq n^{1/\kappa}  \big]}
 \ar\ar\cr
 \ar\ar\cr
 \ar=\ar \mathbf{E}^{(\Lambda )} \big[e^{-\Lambda  S_{n}}   H(e^{-S_n}(1+\hat{I}_{W-1})) ;  |S_{n}|\leq Z, |X_n- an|\leq n^{1/\kappa}  \big]\cr
 \ar=\ar \int_{an- n^{1/\kappa}}^{an+ n^{1/\kappa}} \mathbf{E}^{(\Lambda )} \big[e^{-\Lambda  (S_{n-1}+x)}  H(e^{-S_{n-1}-x}(1+\hat{I}_{W-1}));  |S_{n-1}+x|\leq Z \big]  \mathbf{P}^{(\Lambda )}(X\in dx).
 \eeqlb
 For $\delta>0$, $n\geq 1$ and $i\in\mathbb{Z}$, let $H_{n,\delta}(i) :=  \sup_{i\delta\leq x <  (i+1)\delta}H(e^{-S_{n-1}-an-x}(1+\hat{I}_{W-1}))$. 
 For each integer $-[n^{1/\kappa}/\delta]-1\leq i\leq [n^{1/\kappa}/\delta]+1$, we have uniformly in $x\in [an+i\delta,an+(i+1)\delta) $,
 \beqnn
 \lefteqn{  \mathbf{E}^{(\Lambda )} \big[e^{-\Lambda  (S_{n-1}+x)}  H(e^{-S_{n-1}-x}(1+\hat{I}_{W-1}));  |S_{n-1}+x|\leq Z \big] }\ar\ar\cr
 \ar\ar\cr
 \ar\leq\ar \mathbf{E}^{(\Lambda )} \big[e^{-\Lambda  (S_{n-1}+an+i\delta)} H_{n,\delta}(i);  |S_{n-1}+an+i\delta|\leq Z+\delta \big]. 
 \eeqnn
 Taking this back into the last integral in (\ref{eqn.3.20}), we have 
 \beqnn
 \lefteqn{ \mathbf{E}^{(\Lambda )} \big[e^{-\Lambda  S_{n}}   H(I_{n-W,n}) ;  |S_{n}|\leq Z, |X_1- an|\leq n^{1/\kappa}  \big]}
 \ar\ar\cr
 \ar\leq\ar \sum_{i=-[n^{1/\kappa}/\delta]-1}^{[n^{1/\kappa}/\delta]} \Big( \mathbf{E}^{(\Lambda )} \big[e^{-\Lambda  (S_{n-1}+an+i\delta)} H_{n,\delta}(i);  |S_{n-1}+an+i\delta|\leq Z+\delta \big] \cr
 \ar\ar\cr
 \ar\ar \qquad\qquad\qquad \times  \mathbf{P}^{(\Lambda )}(X\in [an+i\delta,an+(i+1)\delta)) \Big). 
 \eeqnn
 By Assumption~\ref{Con.RegularV} and \ref{Con.LocalRegularV}, for any $ \epsilon>0$ there exists an integer $n_0>1$ such that for any $n\geq n_0$ and $-[n^{1/\kappa}/\delta]-1\leq i\leq [n^{1/\kappa}/\delta]+1$,
 \beqnn
 (1-\epsilon) \cdot B_n\cdot \delta\leq \mathbf{P}^{(\Lambda )}(X\in [an+i\delta,an+(i+1)\delta)) \leq (1+\epsilon) \cdot B_n\cdot \delta
 \eeqnn
 and hence
 \beqlb\label{eqn.2.13}
 \lefteqn{  \frac{1}{B_n}\cdot \mathbf{E}^{(\Lambda )} \Big[e^{-\Lambda  S_{n}}   H(I_{n-W,n}) ;  |S_{n}|\leq Z, |X_1- an|\leq n^{1/\kappa}  \Big]}
 \ar\ar\cr
 \ar\leq\ar    (1+\epsilon)  \sum_{i=-[n^{1/\kappa}/\delta]-1}^{[n^{1/\kappa}/\delta]}   \mathbf{E}^{(\Lambda )} \Big[e^{-\Lambda  (S_{n-1}+an+i\delta)} H_{n,\delta}(i);  |S_{n-1}+an+i\delta|\leq Z+\delta \Big]\cdot  \delta. 
 \eeqlb
 For $y\in \mathbb{R}$, let $ H_\delta(y):= \sup_{y-2\delta\leq x <   y+2\delta}H(e^{-x}(1+\hat{I}_{W-1}))$. For each $-[n^{1/\kappa}/\delta]-1\leq i\leq [n^{1/\kappa}/\delta]$ and $y\in [i\delta,(i+1)\delta)$, we have $H_{n,\delta}(i)\leq H_\delta(S_{n-1}+an+y)$ and 
 \beqnn
 \lefteqn{\mathbf{E}^{(\Lambda )} \big[e^{-\Lambda  (S_{n-1}+an+i\delta)} H_{n,\delta}(i);  |S_{n-1}+an+i\delta|\leq Z+\delta \big]}\ar\ar\cr
 \ar\ar\cr
 \ar \leq\ar  e^{\Lambda \delta}\cdot \mathbf{E}^{(\Lambda )} \big[e^{-\Lambda  (S_{n-1}+an+y)} H_{\delta}(S_{n-1}+an+y );  |S_{n-1}+an+y|\leq Z+2\delta \big] .
 \eeqnn
 Taking this back into (\ref{eqn.2.13}), we have 
 \beqnn
 \lefteqn{  \frac{1}{B_n}\cdot \mathbf{E}^{(\Lambda )} \Big[e^{-\Lambda  S_{n}}   H(I_{n-W,n}) ;  |S_{n}|\leq Z, |X_1- an|\leq n^{1/\kappa}  \Big]}
 \ar\ar\cr
 \ar\leq\ar    (1+\epsilon) e^{\Lambda \delta} \int_{-n^{1/\kappa}-\delta}^{n^{1/\kappa}+\delta} \mathbf{E}^{(\Lambda )} \big[e^{-\Lambda  (S_{n-1}+an+y)} H_{\delta}(S_{n-1}+an+y );  |S_{n-1}+an+y|\leq Z+2\delta \big]dy\cr
 \ar=\ar (1+\epsilon) e^{\Lambda \delta} \int_{-Z-2\delta}^{Z+2\delta} \mathbf{E}^{(\Lambda )} \big[e^{-\Lambda  z} H_{\delta}(z);  |S_{n-1}+an-z|\leq  n^{1/\kappa}+\delta \big]dz.
 \eeqnn
 Here the last equality follows from  the change of variables.
 The statements below Condition~\ref{MomentCon} yields that $n^{-1/\kappa}|S_{n-1}+an-z|\to0$ uniformly in $|z|\leq Z+2\delta$. 
 From this and the dominated convergence theorem, we have $\mathbf{E}^{(\Lambda )} \big[e^{-\Lambda  z} H_{\delta}(z);  |S_{n-1}+an-z|\leq  n^{1/\kappa}+\delta \big] \to \mathbf{E}^{(\Lambda )} \big[e^{-\Lambda  z} H_{\delta}(z) \big]$  as $n\to\infty$. 
 In addition, the continuity of $H$ induces that $H_{\delta}(z)\to H(e^{-z}(1+\hat{I}_{W-1}))$  as $\delta\to 0+$. 
 By Fatou's lemma (as $n\to\infty$) and the dominated convergence theorem (as $\delta, \epsilon \to 0+$), we have 
 \beqnn
 \lefteqn{\limsup_{n\to \infty} \frac{1}{B_n}\cdot \mathbf{E}^{(\Lambda )} \Big[e^{-\Lambda  S_{n}}   H(I_{n-W,n}) ;  |S_{n}|\leq Z, |X_1- an|\leq n^{1/\kappa}  \Big] }\ar\ar\cr
 \ar\leq\ar  
 \lim_{\delta, \epsilon \to 0+} (1+\epsilon) e^{\Lambda \delta} \int_{-Z-2\delta}^{Z+2\delta} \mathbf{E}^{(\Lambda )} \big[e^{-\Lambda  z} H_{\delta}(z)  \big]dz  
 = \int_{-Z}^Z \mathbf{E}^{(\Lambda )} \big[e^{-\Lambda  z}  H(e^{-z}(1+\hat{I}_{W-1})) \big] dz.
 \eeqnn 
 On the other hand, a similar argument also can yields that
 \beqnn
 \lefteqn{ \liminf_{n\to\infty}  \frac{1}{B_n}\cdot \mathbf{E}^{(\Lambda )} \Big[e^{-\Lambda  S_{n}}   H(I_{n-W,n}) ;  |S_{n}|\leq Z, |X_1- an|\leq n^{1/\kappa}  \Big] } \qquad \qquad \ar\ar\cr
 \ar \ar \qquad \geq \int_{-Z }^{Z} \mathbf{E}^{(\Lambda )} \Big[e^{-\Lambda  z} H(e^{-z}(1+\hat{I}_{W-1}))  \Big]dz.
 \eeqnn
 Thus the first claim follows. The second one can be gotten by using dominated convergence theorem as $W\to\infty$ and then the monotone convergence theorem as $Z\to\infty$.
 \qed
 %
 
 Applying the preceding three propositions to $F_{\mathscr{F}^S}^{b,Z}(k,n) $, we can get the following corollary immediately.
 \begin{corollary}\label{Corollary.5.16}
 	For every $b\in(0,a)$ and $k,Z\geq 1$, we have  as $n\to\infty$
 	\beqnn
 	\frac{ 1}{B_n}\cdot F_{\mathscr{F}^S}^{b,Z}(k,n) \to \int_{-Z}^{Z} \mathbf{E}^{(\Lambda )} \big[e^{ -\Lambda   z} F\big(I_{k-1}+ e^{-S_{k-1}-z}(1+\tilde{\hat{I}}_\infty)  \big) \big| \mathscr{F}^S_{k-1}  \big]dz.
 	\eeqnn
 	%
 	
 \end{corollary}
 
 {\it Proof for Theorem~\ref{MainThm.05}.} 
 By (\ref{MeasureChange}), 
 it suffices to prove that $\mathbf{E}^{(\Lambda )} \big[e^{-\Lambda  S_n}F(I_n) \big] \sim C_{F,5}\cdot B_n$ as $n\to\infty$. 
 For $b\in (0,a/r)$ with $r>\frac{\beta+1}{\beta-1}$, there exists a constant $C>0$ such that  
 \beqnn
 \lefteqn{\limsup_{T\to\infty} \limsup_{n\to\infty}\frac{ 1}{B_n}\cdot \mathbf{E}^{(\Lambda )} \big[e^{-\Lambda  S_n}F(I_n);\mathcal{T}^{bn}\geq T \big]}\ar\ar\cr
 \ar\leq\ar  \limsup_{Y\to\infty}\limsup_{T\to\infty} \limsup_{n\to\infty} \frac{1}{B_n}\cdot   \mathbf{E}^{(\Lambda )} \big[e^{-\Lambda  S_n}F(I_n); -Y\leq  L_n\leq S_n\leq Y, \mathcal{T}^{bn}\geq T\big] \cr
 \ar\ar +\limsup_{Y\to\infty} \limsup_{T\to\infty} \limsup_{n\to\infty}  \frac{1}{B_n}\cdot   \mathbf{E}^{(\Lambda )} \big[e^{-\Lambda  S_n}F(I_n);    L_n\leq - Y, \mathcal{T}^{bn}\geq T\big] \cr
 \ar\ar +  \limsup_{Y\to\infty} \limsup_{T\to\infty} \limsup_{n\to\infty} \frac{1}{B_n}\cdot   \mathbf{E}^{(\Lambda )} \big[e^{-\Lambda  S_n}F(I_n);  S_n\geq Y, \mathcal{T}^{bn}\geq T\big] .
 \eeqnn
 By the boundedness of $F$ and  Proposition~\ref{Prop.5.9}-\ref{Proposition.6.12}, the first limit on  the right side of this inequality equals to $0$.  The second limit can be bounded by 
 \beqnn
 \limsup_{Y\to\infty}  \limsup_{n\to\infty}  \frac{1}{B_n}\cdot   \mathbf{E}^{(\Lambda )} \big[e^{-\Lambda  S_n}F(I_n);    L_n\leq - Y \big] ,
 \eeqnn
 which equals to $0$; see  Proposition~\ref{Prop.5.5}. 
 The third limit can be bounded by
 \beqnn
 \lefteqn{\limsup_{Y\to\infty}  \limsup_{n\to\infty} \frac{1}{B_n}\cdot   \mathbf{E}^{(\Lambda )} \big[e^{-\Lambda  S_n}F(I_n);  S_n\geq Y \big] }\ar\ar\cr
 \ar\leq\ar \limsup_{K\to\infty } \limsup_{n\to\infty} \frac{1}{B_n}\cdot   \mathbf{E}^{(\Lambda )} \big[e^{-\Lambda  S_n}F(I_n);    \sigma_n^-\in [K,n-K] \big]\cr
 \ar\ar + \limsup_{K\to\infty } \sum_{k=0}^{K-1} \limsup_{Y\to\infty}\limsup_{n\to\infty} \frac{1}{B_n}\cdot   \mathbf{E}^{(\Lambda )} \big[e^{-\Lambda  S_n}F(I_n);  S_n\geq Y, \sigma_n^- =k \big]\cr
 \ar\ar + \limsup_{K\to\infty } \sum_{k=n-K+1}^{n}\limsup_{Y\to\infty}\limsup_{n\to\infty} \frac{1}{B_n}\cdot   \mathbf{E}^{(\Lambda )} \big[e^{-\Lambda  S_n}F(I_n);  S_n\geq Y, \sigma_n^- =k \big].
 \eeqnn
 By Proposition~\ref{Proposition.5.6}, \ref{Propo.6.9} and \ref{Proposition.5.8},  the three limits on  the right side of  this inequality equal to $0$. In conclusion, we have 
 \beqnn
 \limsup_{n\to\infty}\frac{ 1}{B_n}\cdot \mathbf{E}^{(\Lambda )} \big[e^{-\Lambda  S_n}F(I_n) \big]
 \ar=\ar\lim_{T\to\infty} \lim_{n\to\infty}\frac{ 1}{B_n}\cdot \mathbf{E}^{(\Lambda )} \big[e^{-\Lambda  S_n}F(I_n);\mathcal{T}^{bn}\leq T \big]\cr
 \ar=\ar   \sum_{k=1}^\infty\lim_{n\to\infty}\frac{ 1}{B_n}\cdot \mathbf{E}^{(\Lambda )} \big[e^{-\Lambda  S_n}F(I_n);\mathcal{T}^{bn}=k\big].
 \eeqnn
 By (\ref{eqn5.38.1}) and (\ref{eqn5.39}), we have 
 \beqnn
 \lefteqn{ \limsup_{n\to\infty}\frac{ 1}{B_n}\cdot \mathbf{E}^{(\Lambda )} \big[e^{-\Lambda  S_n}F(I_n) \big]}\ar\ar\cr
 \ar=\ar  \sum_{k=1}^\infty \lim_{Y\to\infty}\lim_{n\to\infty}\frac{ 1}{B_n}\cdot \mathbf{E}^{(\Lambda )} \big[e^{-\Lambda  S_n}F(I_n);S_{k-1}\geq -Y, X_k\geq bn\big]\cr
 \ar=\ar \sum_{k=1}^\infty \lim_{Y\to\infty}\lim_{n\to\infty}\frac{ 1}{B_n}\cdot \mathbf{E}^{(\Lambda )} \big[e^{-\Lambda  S_n}F(I_n);S_{k-1}> Y, X_k\geq bn\big]\cr
 \ar\ar + \sum_{k=1}^\infty \lim_{Y\to\infty}\lim_{Z\to\infty}\lim_{n\to\infty}\frac{ 1}{B_n}\cdot \mathbf{E}^{(\Lambda )} \big[e^{-\Lambda  S_n}F(I_n);|S_{k-1}|\leq Y,|S_n-S_{k-1}|> Z, X_k\geq bn\big]\cr
 \ar\ar + \sum_{k=1}^\infty \lim_{Y\to\infty}\lim_{Z\to\infty}\lim_{n\to\infty}\frac{ 1}{B_n}\cdot \mathbf{E}^{(\Lambda )} \big[e^{-\Lambda  S_n}F(I_n);|S_{k-1}|\leq Y,|S_n-S_{k-1}|\leq Z, X_k\geq bn\big] . 
 \eeqnn
 From Proposition~\ref{Proposition6.13} and \ref{Proposition.5.12}, both of the first two terms on the right side of the last equality equal to $0$.  By (\ref{eqn.5.12}), the third term equals to 
 \beqnn
 \sum_{k=1}^\infty \lim_{Y\to\infty}\lim_{Z\to\infty}\lim_{n\to\infty}  \mathbf{E}^{(\Lambda )} \big[e^{-\Lambda  S_{k-1}} \cdot B_n^{-1}\cdot F_{\mathscr{F}^S}^{b,Z}(k,n) ;|S_{k-1}|\leq Y  \big].
 \eeqnn
 For $Y,Z>0$ fixed, notice that $F_{\mathscr{F}^S}^{b,Z}(k,n)$ is uniformly bounded in $n$. By the dominated convergence theorem and Corollary~\ref{Corollary.5.16}, the preceding sum equals to 
 \beqnn
 \sum_{k=1}^\infty \lim_{Y\to\infty}\lim_{Z\to\infty}  \int_{-Z}^{Z} \mathbf{E}^{(\Lambda )} \big[e^{ -\Lambda  S_{k-1}-\Lambda   z} F\big(I_{k-1}+ e^{-S_{k-1}-z}(1+\tilde{\hat{I}}_\infty)  \big); |S_{k-1}|\leq Y \big]dz . 
 \eeqnn
 Notice that this integral is non-decreasing as $Y,Z\to\infty$.
 By the monotone convergence theorem, the limit in the foregoing limit converges to $C_{F,5}(k)$ and hence 
 \beqnn
 \limsup_{n\to\infty}\frac{ 1}{B_n}\cdot \mathbf{E}^{(\Lambda )} \big[e^{-\Lambda  S_n}F(I_n) \big] = \sum_{k=1}^\infty C_{F,5}(k)= C_{F,5},
 \eeqnn
 which is finite; see the second claim in Proposition~\ref{Proposition.5.6}.
 Here we have got the desired result.
 \qed

 {\bf Acknowledgments.} The author would like to thank Professor Mladen Savov and the two professional referees for their enlightening and helpful comments.

  \bibliographystyle{plain}

 \bibliography{Reference}

 \end{document}